\documentclass{article}
\usepackage{latexsym}
\usepackage{graphics}
\usepackage{amsmath}
\usepackage{amsfonts}
\begin{document}
\bibliographystyle{plain}
\title{The decomposition of global conformal invariants II: 
The Fefferman-Graham ambient metric and the nature of the decomposition.}
\author{Spyros Alexakis\thanks{University of Toronto, alexakis@math.toronto.edu.
\newline
This work has absorbed
the best part of the author's energy over many years. 
This research was partially conducted during 
the period the author served as a Clay Research Fellow, 
an MSRI postdoctoral fellow,
a Clay Liftoff fellow and a Procter Fellow.  
\newline
The author is immensely indebted to Charles
Fefferman for devoting twelve long months to the meticulous
proof-reading of the present paper. He also wishes to express his
gratitude to the Mathematics Department of Princeton University
for its support during his work on this project.}} \date{}
\maketitle
\newtheorem{proposition}{Proposition}
\newtheorem{theorem}{Theorem}
\newcommand{\Sum}{\sum}
\newtheorem{lemma}{Lemma}
\newtheorem{observation}{Observation}
\newtheorem{formulation}{Formulation}
\newtheorem{definition}{Definition}
\newtheorem{conjecture}{Conjecture}
\newtheorem{corollary}{Corollary}
\numberwithin{equation}{section}
\numberwithin{lemma}{section}
\numberwithin{theorem}{section}
\numberwithin{definition}{section}
\numberwithin{proposition}{section}
\begin{abstract}
This is the second in a series of papers where we  
prove a conjecture of Deser and Schwimmer
regarding the algebraic structure of ``global conformal
invariants''; these  are defined to
be conformally invariant integrals of geometric scalars.
 The conjecture asserts that the integrand of
any such  integral can be expressed as a linear
combination of a local conformal invariant, a divergence and of
the Chern-Gauss-Bonnet integrand. 

 The present paper addresses the hardest challenge in this series: It shows 
how to {\it separate} the local conformal invariant from the divergence 
term in the integrand; we make full use of the Fefferman-Graham ambient metric 
to construct the necessary local conformal invariants, as well as all 
the author's prior work \cite{a:dgciI,a:dgciII, alexakis1} 
to construct the necessary divergences. 
This result combined with \cite{alexakis1} completes the proof
of the conjecture, subject to establishing a purely algebraic result 
which is proven in \cite{alexakis4,alexakis5,alexakis6}.
\end{abstract}

\tableofcontents
\section{Introduction.}
 This is the second in a series of papers 
\cite{alexakis1,alexakis3,alexakis4,alexakis5,alexakis6} 
where we prove a conjecture of Deser-Schwmimmer \cite{ds:gccaad} regarding the algebraic 
structure of global conformal invariants. We recall that a global 
conformal invariant is an integral of a natural scalar-valued function of 
Riemannian metrics, $\int_{M^n}P(g)dV_g$, with the property that this integral 
remains invariant under conformal re-scalings of the underlying 
metric.\footnote{See the introduction of \cite{alexakis1}
for a detailed discussion of the Deser-Schwimmer 
conjecture, and for background on scalar Riemannian invariants.} 
More precisely, $P(g)$ is assumed to be a linear combination, $P(g)=\sum_{l\in L} a_l C^l(g)$, 
where each $C^l(g)$ is a complete contraction in the form:

\begin{equation}
\label{contraction} 
contr^l(\nabla^{(m_1)}R\otimes\dots\otimes\nabla^{(m_s)}R);
\end{equation}
 here each factor $\nabla^{(m)}R$ stands for the $m^{th}$ iterated 
covariant derivative of the curvature tensor $R$. $\nabla$ is the Levi-Civita 
connection of the metric $g$ and $R$ is the curvature associated to this connection. 
The contractions are taken with respect to the quadratic form $g^{ij}$.
In this series of papers we prove:

\begin{theorem}
\label{thetheorem} 
Assume that $P(g)=\sum_{l\in L} a_l C^l(g)$, where each $C^l(g)$ is a 
complete contraction in the form (\ref{contraction}), with weight $-n$. 
Assume that for every closed Riemannian manifold $(M^n,g)$ and every $\phi\in C^\infty (M^n)$:
$$\int_{M^n}P(e^{2\phi}g)dV_{e^{2\phi}g}=\int_{M^n}P(g)dV_g.$$

Then $P(g)$ can then be expressed in the form:
$$P(g)=W(g)+div_iT^i(g)+\operatorname{Pfaff}(R_{ijkl}).$$
Here $W(g)$ stands for a local conformal invariant of weight $-n$  (meaning 
that $W(e^{2\phi}g)=e^{-n\phi}W(g)$ for every $\phi\in C^\infty (M^n)$), 
$div_iT^i(g)$ is the divergence of a Riemannian vector field of 
weight $-n+1$, and   $\operatorname{Pfaff}(R_{ijkl})$ is the Pfaffian of the curvature tensor. 
\end{theorem}

 We recall from the introduction in \cite{alexakis1} 
that this entire wok can be naturally subdivided
 into two parts: Part I, consisting 
of \cite{alexakis1}, the present paper and 
\cite{alexakis3} prove Theorem \ref{thetheorem} subject to 
deriving certain purely algebraic propositions, namely the 
``main algebraic Proposition'' 5.2 
in \cite{alexakis1} and Propositions \ref{pregiade2}, \ref{tool'}
in the present paper. Part II, consisting of papers 
\cite{alexakis4,alexakis5,alexakis6} 
is devoted to proving these algebraic Propositions. 

\par In \cite{alexakis1} we explained that our proof of Theorem \ref{thetheorem} 
relied on a main inductive step 
 which asserts that given a 
$P(g)$ as in Theorem \ref{thetheorem}, if the minimum number of factors among all 
complete contractions $C^l(g)$ in $P(g)$ is $\sigma<\frac{n}{2}$, 
then we can {\it subtract} a divergence and a local conformal invariant from $P(g)$ 
so as to {\it cancel out} the terms with $\sigma$ factors in $P(g)$, 
modulo introducing new terms $\sigma+1$ factors. In conjunction with 
the results of \cite{a:dgciI, a:dgciII}, this main inductive step 
implies Theorem \ref{thetheorem}. 

This main inductive step consists of two sub-steps, the Propositions 3.1, 3.2 in section 3 
in \cite{alexakis1}. Proposition 3.1  
was proven in \cite{alexakis1} 
(subject to deriving the ``main algebraic Proposition'' 5.2 there). 
The present paper is devoted to proving the second (harder) Proposition 3.2.  

We state this second Proposition again, after recalling two main 
pieces of notation:

{\bf Conventions:} For any complete or partial contraction in $\alpha$ 
factors $T^1,\dots$,\\$T^{\alpha}$,\footnote{See the first section in \cite{a:dgciI} 
for a very rigorous definition of complete and partial contractions.}
 $contr(T^1\otimes\dots\otimes T^\alpha)$, an {\it internal contraction }
is a pair of indices in a given factor $T^\beta$  which contract against each other. 
For  each complete contraction in the form (\ref{karami}), $\delta_W$
stands for the total number of internal contractions. Also, for each complete or partial contraction, 
its ``length'' will be its number of factors. 

Proposition 3.2 in \cite{alexakis1} states: 

\begin{proposition}
\label{bigprop}
\par Consider any $P(g)$,
$P(g)=\Sum_{l\in L} a_l C^l(g)$ where each $C^l(g)$ has length
$\ge\sigma$, and each $C^l(g)$ of length $\sigma$ 
 is in the form:\footnote{$W$ below stands for the Weyl tensor 
($W_{ijkl}$ if we write out the indices), the trace-free 
part of the curvature tensor. $\nabla^{(m)}W$
is the $m^{th}$ covariant derivative of the Weyl tensor.}
\begin{equation}
\label{karami}
contr(\nabla^{(m_1)}W\otimes\dots\otimes\nabla^{(m_\sigma)}W).
\end{equation}
Assume that $\int_{M^n}P(g)dV_g$ is a global conformal invariant. 
 Denote by $L_\sigma\subset L$ 
the index set of terms with length $\sigma$.

We claim that there is a local conformal invariant $W(g)$
 of weight $-n$ and also a vector field $T^i(g)$ as in the statement
of Theorem \ref{thetheorem}, so that:

\begin{equation}
\label{karami2} \Sum_{l\in L_\sigma} a_l C^l(g)-W(g) -div_i
T^i(g)=0,
\end{equation}
modulo complete contractions in the form (\ref{contraction}) of length $\ge\sigma+1$.
\end{proposition}

{\it Note:} As explained in \cite{alexakis1}, we prove this proposition for $\sigma\ge 3$. 
The special cases $\sigma=1,\sigma=2$ are treated in section 2 in \cite{alexakis3}.
\newline

{\bf Digression: A general discussion on local 
conformal invariants.} Since the present paper deals extensively with the issue of constructing 
{\it local conformal invariants} (as required in (\ref{karami2})),
we digress slightly to discuss some background material regarding 
the history and known constructions  of such invariants:

 The theory of {\it local} invariants of Riemannian structures 
(and indeed, of more general geometries,
e.g.~conformal, projective, or CR)  has a long history. As stated above, the original foundations of this 
field were laid in the work of Hermann Weyl and \'Elie Cartan, see \cite{w:cg, cartan}. 
The task of writing out local invariants of a given geometry is intimately connected
with understanding polynomials in a space of tensors with  given symmetries, 
which remain invariant under the action of a Lie group. 
In particular, the problem of writing down all 
 local Riemannian invariants reduces to understanding 
the invariants of the orthogonal group. 

 In more recent times, a major program was laid out by C.~Fefferman in \cite{f:ma}
aimed at finding all local invariants in CR geometry. This was motivated 
by the problem of understanding the 
local invariants which appear in the asymptotic expansions of the 
Bergman and Szeg\"o kernels of CR manifolds,
 in a similar way to which Riemannian invariants appear in the asymptotic expansion  
of the heat kernel; the study of the local invariants
in the singularities of these kernels led to important breakthroughs 
in \cite{beg:itccg} and more recently by Hirachi in \cite{hirachi1}.
 This program was later extended  to conformal geometry in \cite{fg:ci}. 
Both these geometries belong to a 
broader class of structures, the
{\it parabolic geometries}; these admit a principal bundle with 
structure group a parabolic subgroup $P$ of a semi-simple 
Lie group $G$, and a Cartan connection on that principle bundle 
(see the introduction in \cite{cg1}). 
An important question in the study of these structures 
is the problem of constructing all their local invariants, which 
can be thought of as the {\it natural, intrinsic} scalars of these structures.

  In the context of conformal geometry, the first (modern) landmark 
in understanding {\it local conformal invariants} was the work of Fefferman 
and Graham in 1985 \cite{fg:ci},
where they introduced the {\it ambient metric}. This allows one to 
construct local conformal invariants of any order in odd 
dimensions, and up to order $\frac{n}{2}$ in even dimensions. 
The question is then whether {\it all} invariants arise via this construction. 

The subsequent work of Bailey-Eastwood-Graham \cite{beg:itccg} proved that
indeed in odd dimensions all conformal invariants arise 
via the Fefferman-Graham construction; in even dimensions, 
they proved that the same holds true 
when the weight in absolute value is bounded above by the dimension. The ambient metric construction 
in even dimensions was recently extended by Graham-Hirachi, \cite{grhir}; this enables them to then
indentify in a satisfactory way all local conformal invariants 
even when the weight (in absolute value) exceeds the dimension.  

 An alternative 
construction of local conformal invariants can be obtained via the {\it tractor calculus} 
introduced by Bailey-Eastwood-Gover in \cite{bego}. This construction bears a strong 
resemblance to the Cartan conformal connection, and to 
the work of T.Y.~Thomas, \cite{thomas}. The tractor 
calculus has proven to be very universal; 
tractor buncles have been constructed \cite{cg1} for an entire class of parabolic geometries. 
The relation betweeen the conformal tractor calculus and the Fefferman-Graham 
ambient metric construction has been elucidated in \cite{cg2}.

The present series of papers \cite{alexakis1}--\cite{alexakis6}, 
while pertaining to the question above
(given that it ultimately deals with the algebraic form of local 
 Riemannian and  conformal invariants), nonetheless addresses a different 
{\it type} of problem:  We here consider Riemannian invariants $P(g)$ for 
which the {\it integral} $\int_{M^n}P(g)dV_g$ remains invariant 
under conformal changes of the underlying metric; we then seek to understand 
the possible algebraic form of the {\it integrand} $P(g)$, 
ultimately proving that it can be de-composed 
in the way that Deser and Schwimmer asserted. 
It is thus not surprising that the prior work on 
 the construction and understanding of local {\it conformal} 
invariants plays a central role in this endeavor, 
in the present paper and in \cite{alexakis3}.

On the other hand, our resolution of the Deser-Scwimmer conjecture will also rely 
heavily on a deeper understanding of the algebraic properties 
of the {\it classical} local Riemannian invariants.
 The fundamental theorem of invariant theory (see Theorem B.4 in 
\cite{beg:itccg} and also Theorem 2 in \cite{a:dgciI})
is used extensively throughout this series of papers. However, the most important 
algebraic tool on which our method relies are certain ``main algebraic Propositions'' 
presented in \cite{alexakis1} and in the present paper. 
These are purely algebraic propositions that deal with {\it local Riemannian invariants}. 
While the author was led to these Propositions 
out of the strategy that he felt was necessary to 
solve the Deser-Schwimmer conjecture, they can 
be thought of as results with an independent interest. 
The {\it proof} of these Propositions, presented
 in \cite{alexakis4,alexakis5,alexakis6} is in fact 
not particularily intuitive. It is the author's 
sincere hope that deeper insight 
will be obtained in the future as to  {\it why} these algebraic 
 Propositions hold. 
\newline

{\bf Local conformal invariants in our proof:} The 
first half of this paper deals largely with the problem
 of identifying a local conformal invariant {\it in} $P(g)$. 
As explained in more detail in the ``outline'' below, we construct three {\it different
kinds} of local conformal invariants,  each of which
 will ``cancel out'' a particular {\it kind} of terms in $P(g)$.  
 The next challenge is to then cancel out the {\it remaining piece} in
 $P(g)|_\sigma:=\sum_{l\in L_\sigma} a_l C^l(g)$ 
by subtracting {\it only} a divergence of a vector field. 
As explained above, one has a powerful tool in the challenge of 
  constructing local conformal invariants; this is the Fefferman-Graham 
ambient metric, \cite{fg:ci, fg:latest}. 
A note is in order here: The roughest simplification 
of $P(g)|_\sigma$ is to ``cancel out'' terms $C^l(g)$ in $P(g)$ 
which do not involve {\it internal contractions}.\footnote{I.e.~we cancel out 
the terms in the form  (\ref{karami}) which contain 
 no factor $\nabla^{(m)}W$ 
with two indices  contracting against each 
other; we do this modulo introducing 
{\it new terms} in the form (\ref{karami})
which {\it do} contain such internal contractions. This is the content of Lemma \ref{subcinv}.}
For many reasons, the use of local conformal invariants to cancel out 
these particular terms in $P(g)$ is hardly surprising, in view (for example) 
of the proof of Proposition 3.2 in \cite{beg:itccg}. However, we 
believe that the use of local conformal invariants for 
the next two ``simplifications'' of $P(g)$ (in Lemmas 
\ref{1stclaim}, \ref{killstacks}), {\it is} somewhat surprising: 
We prove that we can start with complete contractions in the 
form (\ref{karami}) {\it which do contain internal contractions}, 
and {\it cancel out} {\it certain particular terms in this form} in   
$P(g)$ by subtracting local conformal invariants and explicitly 
constructed divergences; we do this modulo introducing {\it new terms} 
which are ``better'' (from our point of view) 
than the terms we cancelled out. The construction 
of local conformal invariants for Lemmas \ref{1stclaim}, \ref{killstacks} depends essentially on 
explicitly constructed {\it divergences
in the ambient metric},\footnote{These are 
{\it local conformal invariants}, by construction; they 
are (at least apriori) {\it not} divergences for the base metric $g$.} and the subsequent  
of cancellations that occur in these constructions.

\subsection{Outline of the argument.}
\label{outlofarg}
The proof of Proposition \ref{bigprop} addresses two main challenges:
Firstly, how to separate the local conformal invariant ``piece'' $W(g)$  
in $\sum_{l\in L_\sigma} a_l C^l(g)$ from the divergence ``piece'' $div_iT^i(g)$. 
Secondly, how to use the {\it local equation} (i.e.~the ``super divergegce formula'') 
that we have derived in \cite{a:dgciI} regarding the {\it conformal variation} of $P(g)$
to construct the divergence $div_iT^i(g)$ needed in \ref{karami2}. 

\par Our argument proceeds as follows: At a first step 
we explicitly construct a local conformal invariant
$W(g)$ and a divergence $div_iT^i(g)$ such that:
$$\sum_{l\in L_\sigma} a_l C^l(g)=W(g)+div_iT^i(g)+\sum_{l\in L_{new}} C^l(g),$$
where the complete contractions $C^l(g), l\in L_{new}$ are in the form (\ref{karami}),
but have certain {\it additionnal} algebraic properties.  
(See Lemmas 
\ref{subcinv}, \ref{1stclaim}, \ref{killstacks} below).

Thus, the first step reduces matters to proving Proposition \ref{bigprop} with 
\\$\sum_{l\in L_\sigma} a_l C^l(g)$ replaced by $\sum_{l\in L_{new}} a_l C^l(g)$.
In step 2, we then revert to studying the {\it new} $P(g)$ by focusing on 
the conformal variation $I^1_g(\phi)$ and the super divergence 
formula for $I^1_g(\phi)$.\footnote{Recall that $I^1_g(\phi):=\frac{d}{dt}|_{t=0} [e^{nt\phi}P(e^{2t\phi}g)]$, 
$\int_{M^n}I^1_g(\phi)dV_g$; we have then derived a useful local formula in \cite{alexakis1}, 
which we called the ``super divergence formula''.} We prove that 
 $\sum_{l\in L_{new}} a_l C^l(g)=div_iT^i(g)$, for some vector field $T^i(g)$. 
(This is proven in Lemma \ref{nabothewri2} below). 

A few rermarks: It is not at all clear that the local conformal invariant
we construct to prove (\ref{karami2}) is the {\it unique} 
local conformal invariant for which (\ref{karami2}) 
is true. It is also not clear that one cannot 
subtract further conformal invariants from $\sum_{l\in L_{new}} a_l C^l(g)$ 
in order to simplify it even further. At any rate, (as discussed in section 
\ref{confinvariant}) the local conformal invariants $W(g)$ that one subtracts from $P(g)$ 
in order to simplify it are all explicitly constructed using 
the Fefferman-Graham ambient metric, 
\cite{fg:ci, fg:latest}. Our construction 
is elaborate and relies on the study of linear combinations of 
complete contractions in the ambient metric with specific algebraic properties. 
It is also worth noting that the algebraic properties of the terms 
in $\sum_{l\in L_{new}} a_l C^l(g)$ are precisely what is needed in order 
to prove (by the methods in section \ref{puredivergencesec}) that 
$\sum_{l\in L_{new}} a_l C^l(g)=div_iT^i(g)$. 
\newline

We now discuss in more detail the two main steps in the proof of Proposition \ref{bigprop}.
 {\it First Step:} We show (in Lemmas \ref{subcinv},
\ref{1stclaim}, \ref{killstacks} below) that there exists a local
conformal invariant $W(g)$ and a divergence $div_iT^i(g)$ so that:

\begin{equation}
\label{bo'az9} P(g)=W(g)+div_i T^i(g)+\sum_{f\in F}
a_f C^f(g)+\sum_{j\in Junk-Terms} a_j C^j(g),
\end{equation}
where the terms indexed in  $Junk-Terms$ have at least $\sigma+1$
factors, while the terms indexed in  $F$  are in the form
(\ref{karami}) {\it and have internal contractions in at least
two different factors}. 
\newline

The second step is to show that there exists another divergence
$div_iT'^i(g)$ so that (in the notation of (\ref{bo'az9})):

\begin{equation}
\label{bo'az11} \sum_{f\in F} a_f C^f(g)=div_i T'^i(g)+\sum_{j\in
Junk-Terms} a_j C^j(g).
\end{equation}

\par The above two sub-steps combined prove Proposition \ref{bigprop}.
\newline

{\it An outline of the proof of (\ref{bo'az9}):} As explained, the proof of this step relies heavily on the 
{\it ambient metric} of Fefferman and Graham \cite{fg:ci, fg:latest}. 
This is a strong tool that 
allows one to explicitly construct all {\it local conformal invariants} of 
weight $-n$ in dimension $n$. By making 
detailed use of the precise form of the ambient metric (and in one 
instance of the super divergence formula), we are able to {\it explicitly} construct 
the local conformal invariant $W(g)$  and the divergence $div_iT^i(g)$ 
needed for (\ref{bo'az9}). The ``main algebraic Propositions'' 
\ref{pregiade2} and \ref{tool'} of the present paper are 
{\it not} used in deriving (\ref{bo'az9}).

{\it An outline  of the proof of (\ref{bo'az11}):} (\ref{bo'az11})
is proven by a new induction: {\it At a  rough level,} the new
induction can be described as follows: We denote by $j$ the {\it
minimum number of internal contractions} among the complete
contractions indexed in $F$. Denote the corresponding index set by
$F_j\subset F$. We then  prove that
we can write:

\begin{equation}
\label{bo'az12} \sum_{f\in F_j} a_f C^f(g)=div_i
T'^i(g)+\sum_{f\in F'_{j+1}} a_f C^f(g)+ \sum_{t\in Junk-Terms}
a_t C^t(g).
\end{equation}
Here the terms indexed in $F'_{j+1}$ are complete contractions in
the form (\ref{bo'az11}) with $j+1$ internal contractions in
total, of which at least two belong to different factors.

\par Observe that (\ref{bo'az11}) clearly
 follows  
 by iterative application of  (\ref{bo'az12}).\footnote{This 
is because there can be at most $\frac{n}{2}$ internal contractions
in any non-zero complete contraction in the form (\ref{karami})
with weight $-n$.}

\par Now, there are a number of difficulties in proving
(\ref{bo'az12}):
\newline

{\it How can one ``recognize'' the linear combination $\sum_{f\in F_j} a_f
C^f(g)$ (which appears in $P(g)$) in $I^1_g(\phi)$?} We observe
that in this setting, $I^1_g(\phi)$ can be expressed as:

\begin{equation}
\label{goforit}
\begin{split}
& I^1_g(\phi)=\sum_{x\in X} a_x contr^x(\nabla^{(m_1)}R\otimes
\dots\otimes\nabla^{(m_\sigma)}R\otimes\nabla^{(p)}\phi)+\sum_{j\in Junk-Terms} a_j C^j(g),
\end{split}
\end{equation}
where the complete contractions indexed in $X$
have $\sigma+1$ factors in total, while
$\sum_{j\in Junk-Terms} a_j C^j(g)$ stands for a generic linear combination of
 complete contractions with at least $\sigma+2$ factors. We denote by
 $(I^1_g(\phi))_{\nabla\phi}$ the ``piece'' in $\sum_{x\in X}\dots$
which consists of terms with a factor $\nabla\phi$.

\par Observe that $(I^1_g(\phi))_{\nabla\phi}$ 
can {\it only} arise from the terms in $P(g)$ with $\sigma$ factors in
total.\footnote{This is by virtue of the transformation law for the Levi-Civita
connection, see (\ref{levicivita})--these facts are explained in detail in subsection \ref{prelimwork}.}

 In fact (as we will show in subsection \ref{prelimwork}),
we can easily reconstruct  \\$\sum_{f\in
F_j} a_f C^f(g)$ in $P(g)$ if we are given the ``piece''
$(I^1_g(\phi))_{\nabla\phi}$ in $I^1_g(\phi)$. Thus, our aim is to
use the super divergence formula for $I^1_g(\phi)$ to express 
$(I^1_g(\phi))_{\nabla\phi}$ as ``essentially a
divergence''.
\newline

 The most important difficulty in deriving (\ref{bo'az12}) appears when we 
apply the super divergence formula to $I^1_g(\phi)$. To illustrate
why this case is harder than the case $s>0$ 
which was treated in \cite{alexakis1}, we will 
note that $I^1_g(\phi)$ can be expressed as follows:

\begin{equation}
\label{bo'az13}
I^1_g(\phi)=\big{(}I^1_g(\phi)\big{)}_{\nabla\phi}+ \sum_{q\in Q}
a_q  C^q(g)\cdot \Delta\phi+\sum_{z\in Z} a_z C^z_g(\phi)
+\sum_{j\in Junk-Terms} a_j C^j_g(\phi),
\end{equation}
where the terms $C^q(g)\cdot \Delta\phi$ are in the form:
\begin{equation}
\label{bo'az14}
contr(\nabla^{(m_1)}R\otimes\dots\otimes\nabla^{(m_\sigma)}R\otimes\Delta\phi),
\end{equation}
(notice they have $\sigma+1$ factors in total), while the terms
$C^z_g(\phi)$ are in the form:
\begin{equation}
\label{bo'az15}
contr(\nabla^{(m_1)}R\otimes\dots\otimes\nabla^{(m_\sigma)}R\otimes\nabla^{(p)}\phi),
\end{equation}
where $p\ge 2$ and moreover if $p=2$ then the two indices
${}_a,{}_b$ in $\nabla^{(2)}_{ab}\phi$ are {\it not} contracting
against each other in (\ref{bo'az15}). Furthermore, it follows
that\footnote{See subsection \ref{prelimwork}.} the terms
$C^q(g)\cdot \Delta\phi$ in (\ref{bo'az13}) arise {\it both} from
the terms with $\sigma$ factors in $P(g)$ {\it and} from terms
with $\sigma+1$ factors in $P(g)$.

\par Now, to obtain a {\it local formula} for the expression
$\big{(}I^1_g(\phi)\big{)}_{\nabla\phi}$  we consider the super
divergence formula applied to $I^1_g(\phi)$ and pick out the terms
which have $\sigma+1$ factors in total, and furthermore have a factor
$\nabla\phi$ (differentiated only once). It follows that the
linear combination of {\it those} terms in $supdiv[I^1_g(\phi)]$
must vanish separately (modulo junk terms with more that
$\sigma+1$ factors): We thus obtain a new {\it local} equation,
which we denote by:
\begin{equation}
\label{den3eroume}
supdiv_{+}[I^1_g(\phi)]=\sum_{j\in Junk-Terms} a_j C^j(g).
\end{equation}

 Now, if one
follows the algorithm for the super divergence
formula,\footnote{(See subsection 6.3 in \cite{a:dgciI}, and also the ``main 
consequence of the super divergence formula'' in
 \cite{alexakis1}).} one observes that the terms in the RHS of (\ref{bo'az13})
 which can contribute
terms to $supdiv_{+}[I^1_g(\phi)]$ are the terms in
$\big{(}I^1_g(\phi)\big{)}_{\nabla\phi}$ in (\ref{bo'az13}), but
{\it also} the terms \\$\sum_{q\in Q} a_q  C^q(g)\cdot \Delta\phi$
in (\ref{bo'az13}).\footnote{In particular,
the terms $C^q(g)\cdot \Delta\phi$ can give rise to an iterated
divergence of a tensor field with a factor $\nabla_i\phi$ (with
the index ${}_i$ being free).} Thus, whereas in $I^1_g(\phi)$
({\it before} we consider the super divergence formula), the terms
 with a factor $\nabla\phi$ (and with $\sigma+1$ factors in total) can
{\it only} arise from the ``worst piece'' of $P(g)$,\footnote{(See the statement 
of Proposition \ref{bigprop}). This
is good for our purposes, since we are trying to ``recover'' the
``piece'' $\sum_{l\in L_\sigma} a_l C^l(g)$ in $P(g)$ 
by examining $I^1_g(\phi)$.}  {\it after} we apply 
the super divergence formula we {\it also} obtain terms in
$supdiv_{+}[I^1_g(\phi)]$ which arise from $\sum_{q\in Q} a_q C^q(g)\cdot \Delta\phi$;
this is very worrisome since {\it we have no information on the algebraic form 
of $\sum_{q\in Q} a_q C^q(g)\cdot \Delta\phi$}.\footnote{This
is bad for our purposes, because it is then not obvious how the
contribution of the terms $\sum_{q\in Q} a_q  C^q(g)\cdot
\Delta\phi$ to the formula $supdiv_{+}[I^1_g(\phi)]=(Junk-Terms)$
can be {\it distinguished} from the contribution of the terms
$\big{(}I^1_g(\phi)\big{)}_{\nabla\phi}$.} At this point
 the fact that the terms with length $\sigma$ in $P(g)$ are all
in the form (\ref{karami}) and have two internal contractions
belonging to different factors is crucial. We prove (\ref{bo'az11}) 
(explaining how to overcome these difficulties)
in section \ref{puredivergencesec}.

\subsection{Divide Proposition \ref{bigprop} into four smaller claims.}

The main assumption for all four of the Lemmas below is that
 $\int_{M^n}P(g)dV_g$ is a global conformal invariant and
$P(g)=\sum_{l\in L} a_l C^l(g)$, where each
$C^l(g)$ has at least $\sigma$ factors. The contractions $C^l(g)$
that do have $\sigma$ factors will be indexed in $L_\sigma\subset
L$ and will all be in the form (\ref{karami}).
\newline

{\bf First three Lemmas: The local conformal invariant ``piece'' in $P(g)$:}
 We firstly focus on the complete contractions $C^l(g), l\in
L_\sigma$ with {\it no internal
 contractions}. We index those complete contractions in the
set $L^0_\sigma$ and consider the sublinear combination
$P(g)|_{L^0_\sigma}(:=\sum_{l\in L_\sigma^0} a_l C^l(g))$. Our first claim is the following:

\begin{lemma}
\label{subcinv} Assume that $P(g)=\sum_{l\in L} a_l C^l(g)$ 
is as in the assumption of Proposition \ref{bigprop}.
We claim that there is a  conformally invariant scalar $W(g)$ of weight $-n$ such
that:
$$P(g)-W(g)={\Sigma}_{t\in T} a_t C^t(g)+
\Sum_{u\in U} a_u C^u(g),$$ where all the complete contractions
$C^t(g)$ are in the
 form (\ref{karami}) with length $\sigma$ and $\delta_W\ge 1$.
Each $C^u(g)$ has length $\ge\sigma +1$.
\end{lemma}

\par Assuming we can prove the above Lemma, we are reduced to proving Proposition 
\ref{bigprop} under the
 extra assumption that all complete contractions $C^l(g)$,
$l\in L_\sigma$ have at least one internal contraction (in other words we may now 
assume that $L^0_\sigma=\emptyset$).

\par We then denote by $L^1_\sigma\subset L_\sigma$ the index set of the
complete contractions $C^l(g), l\in L^1_\sigma$ with one
 internal contraction. We claim:

\begin{lemma}
\label{1stclaim} Assume that $P(g)=\sum_{l\in L} a_l C^l(g)$ 
satisfies the assumptions of Proposition \ref{bigprop}; additionaly, 
assume that each $C^l(g), l\in L_\sigma$ satisfies $\delta_W\ge 1$. We claim that
there is a conformally invariant scalar $W(g)$
and a vector
 field $T^i(g)=\Sum_{r\in R} a_r C^{r,i}(g)$ of weight $-n+1$, where each
 $C^{r,i}(g)$ is in the form (\ref{karami2}) such that:

\begin{equation}
\label{kinhta} \Sum_{l\in L^1_\sigma} a_l C^l(g)-W(g)-div_iT^i(g)=
\Sum_{f\in F} a_f C^f(g),
\end{equation}
where each $C^f(g)$ is in the form (\ref{karami}) and has
$\delta_W\ge 2$. The above holds modulo complete contractions of
length $\ge\sigma+1$.
\end{lemma}

\par Assuming we can prove the above Lemma, we are reduced to proving Proposition 
\ref{bigprop} under the
 extra assumption that each
$C^l(g), l\in L_\sigma$ has $\delta_W\ge 2$.

\begin{lemma}
\label{killstacks} Assume that $P(g)=\sum_{l\in L} a_l C^l(g)$ satisfies the 
assumptions of Proposition \ref{bigprop}; assume moreover that
each $C^l(g), l\in L_\sigma$ has
$\delta_W \ge 2$. Consider the index set $L^{stack}_\sigma\subset L_\sigma$ that
consists of complete contractions in the form (\ref{karami}) that
have all their $\delta_W$ internal contractions
 belonging to the same factor $\nabla^{(m)}W_{ijkl}$. 

\par We claim that
 there is a scalar conformal invariant $W(g)$ and a vector
 field $T^i(g)=\Sum_{r\in R} a_r C^{r,i}(g)$ of weight $-n+1$
 (where each $C^{r,i}(g)$ is in the form (\ref{karami2})) so that:

\begin{equation}
\label{hadjipapares} \Sum_{l\in L^{stack}} a_l C^l(g)-W(g)-div_i
T^i(g)=\Sum_{j\in J} a_j C^j(g),
\end{equation}
where each $C^j(g)$ is in the form (\ref{karami}), has $\delta_W\ge
2$ and also at least two internal contractions
 belonging to different factors. The above holds modulo complete
 contractions of length $\ge\sigma+1$.
\end{lemma}

{\bf Next Claims: The divergence ``piece'' in $P(g)$.}

\par For our next claim we will be assuming that for  
$P(g)=\sum_{l\in L} a_l C^l(g)$ (which satisfies the assumptions of Proposition \ref{bigprop}),
all complete contractions $C^l(g)$, $l\in L_\sigma$ are in the form
(\ref{karami}) with $\delta_W\ge 2$ internal contractions, {\it and at
least two internal contractions belong to different factors}.

Let $j=min_{l\in L_\sigma}\{\delta_W[C^l(g)]\}$. We define
$L^j_\sigma\subset L_\sigma$ to stand for the index set of
 the complete contractions with $\delta_W=j$.

\begin{lemma}
\label{nabothewri2} Assume that $P(g)=\sum_{l\in L} a_l C^l(g)$
satisfies the assumptions of Proposition \ref{bigprop}; 
assume additionaly that all complete contractions
$C^l(g)$, $l\in L_\sigma$ have $\delta_W\ge 2$ internal
contractions and at least two of those internal 
contractions belong to different factors.

\par  We claim that there is a linear
combination of partial contractions, $\Sum_{h\in H} a_h
C^{h,i}(g)$, where each $C^{h,i}(g)$ is in the form (\ref{karami2})
with weight $-n+1$ and $\delta_W=j$, so that:

\begin{equation}
\label{salvation} \Sum_{l\in L^j_\sigma} a_l C^l(g)-div_i
\Sum_{h\in H} a_h C^{h,i}(g)=\Sum_{v\in V} a_v C^v(g),
\end{equation}
where each $C^v(g)$ is in the form (\ref{karami}) with $\delta_W\ge
j+1$, and with at least two internal contractions
 in different factors. The above equation holds modulo
 complete contractions of length $\ge\sigma +1$.
\end{lemma}

\par Clearly, if we can show the above four Lemmas,
Proposition \ref{bigprop} will follow: First 
applying the first three Lemmas to $P(g)$ to derive that there exists a 
local conformal invariant $W(g)$ and a divergence $div_iT^i(g)$  
as claimed in Proposition \ref{bigprop} such that:

$$P(g)=W(g)+div_iT^i(g)+\sum_{l\in L'} a_l C^l(g)+\sum_{j\in J} a_j C^j(g);$$
here each $C^j(g)$ has length $\ge\sigma+1$; each $C^l(g), l\in L'$ 
has length $\sigma$ and is in the form (\ref{karami}) 
and moreover has $\delta_W\ge 2$ internal contractions, 
at least two of which belong to different factors.  Then, iteratively applying 
(\ref{nabothewri2}) to $P'(g):=P(g)-W(g)-div_iT^i(g)$,\footnote{Notice 
that $\int_{M^n}P'(g)dV_g$ is {\it also} a global conformal invariant.} 
we derive that there exists a divergence $div_iT'^i(g)$ 
as required by Proposition \ref{bigprop} such that:
$$\sum_{l\in L'} a_l C^l(g)=div_iT'^i(g)+\sum_{j\in J} a_j C^j(g),$$
where each $C^j(g)$ has length $\ge\sigma+1$. 
Therefore, the above four Lemmas indeed imply Proposition \ref{bigprop}. 
 We present the proofs of these four Lemmas in the remainder of this paper.

\section{The locally conformally invariant ``piece'' in $P(g)$: A proof of Lemmas 
\ref{subcinv}, \ref{1stclaim}, \ref{killstacks}.}
\label{confinvariant}

\subsection{The Fefferman-Graham construction of local conformal invariants 
and an algorithm for computations.}
\label{algorithm}

\par We start with a brief discussion of the ambient metric of Fefferman and Graham, 
\cite{fg:ci, fg:latest} which we use to construct local conformal invariants. 
 The ambient metric construction provides a canonical 
embedding of a Riemannian manifold $(M^n,g)$ into an 
ambient Ricci-flat Lorentzian manifold 
$(\tilde{G}^{n+2},\tilde{g}^{n+2})$. We refer the reader to the papers  
 \cite{fg:ci}, \cite{fg:latest} for a detailed exposition of this construction.
We recall here a few features of this construction which will be useful to us: 
 Recall that given  coordinates $\{x^1,\dots ,x^n\}$ for $(M^n,g)$,
 then the ambient metric embedding 
provides a special coordinate system $\{x^0,x^1,\dots ,x^n,x^{n+1}\}$
for the ambient manifold $(\tilde{G}^{n+2},\tilde{g}^{n+2})$. More precisely, 
any given point $x_0\in M^n$ is mapped to $\tilde{x}_0$ in $\tilde{G}^{n+2}$ with $x^0=1$ and
$x^{n+1}=0$. We denote the vectors in $T\tilde{G}|_{\tilde{x}_0}$
that correspond to the directions of $x^0,\dots ,x^n,x^{n+1}$ by
$X^0,\dots , X^n, X^{\infty}$ respectively. In what follows, 
we will often use the notion of  ``assigning of values'' 
${}_0,{}_1,\dots,{}_n,{}_\infty$ to the (lower) indices of 
tensors. By this we will mean that we evaluate those 
(covariant) tensors against the vectors $X^{0},\dots,X^\infty$.

\par Now, let us furthermore recall that in the coordinate system $\{x^0,\dots
x^{n+1}\}$ the ambient metric at
$\tilde{x}_0$ is of the form:

\begin{equation}
 \label{ambform}
\tilde{g}^{n+2}_{IJ}dx^Idx^J =2dx^0dx^{n+1}+ g_{ij}dx^idx^j,
\end{equation}
where $0\le I,J\le n+1$ and $1\le i,j\le n$.

The Christoffel symbols of the ambient metric 
(evaluated with respect to this special
coordinate system at a point with coordinates $(t,x^1,\dots,x^n,0)$) are:

\begin{equation}
\label{cnr}
\begin{split}
&\tilde{\Gamma}_{IJ}^0 = 
\left(
\begin{matrix}
0&0&0\\
0&- tP_{ij}&0\\
0&0&0
\end{matrix}
\right),
\\&\tilde{\Gamma}_{IJ}^k = 
\left(
\begin{matrix}
0&t^{-1}\delta_j{}^k&0\\
t^{-1}\delta_i{}^k&\Gamma_{ij}^k& g^{kl}P_{il}\\ 
0& g^{kl}P_{jl}&0 
\end{matrix}
\right),
\\&\tilde{\Gamma}_{IJ}^\infty = 
\left(
\begin{matrix}
0&0&t^{-1}\\
0&-g_{ij} &0\\
t^{-1}&0&0
\end{matrix}
\right).
\end{split}
\end{equation}

\par Now, as proven in \cite{fg:ci}, \cite{fg:latest}, 
one can construct local conformal 
invariants by considering complete contractions involving the 
{\it ambient} curvature tensor $\tilde{R}$ and its {\it ambient} covariant derivatives
(with respect to the Levi-Civita connection $\tilde{\nabla}$ of the ambient metric). 
In other words, Fefferman and Graham proved that any linear combination  of 
complete contractions (of weight $-n$) in the form:

\begin{equation}
\label{ambient} contr({\tilde{\nabla}}_{r_1\dots
r_{m_1}}^{(m_1)}{\tilde{R}}_{i_1j_1k_1l_1} \otimes \dots \otimes
{\tilde{\nabla}}_{v_1\dots
v_{m_s}}^{(m_s)}{\tilde{R}}_{i_sj_sk_sl_s}),
\end{equation}
say $L(\tilde{g})=\sum_{h\in H} a_h C^h(\tilde{g})$, 
is {\it by construction} a local conformal 
invariant of weight $-n$.\footnote{When $n$ is even, 
the jet of the ambient metric at $\tilde{x}^0$ is only 
defined up to order $\frac{n}{2}-1$. In our constructions, 
this restriction will always be fulfilled.} 
In particular, 
$L(\tilde{g})$ can {\it also} be expressed a linear combination 
of complete contractions in the form (\ref{karami}), involving 
covariant derivatives of the curvature tensor of 
the metric $g$, so $L(\tilde{g})=F(g)$, where 
$F(g):=\sum_{h\in H} a_h contr^h(\nabla^{(m)}R\otimes\dots\otimes\nabla^{(m')}R)$. 
Furthermore, for any function $\lambda>0$, we will 
have $F(\lambda^2\cdot g)=\lambda^{-n}F(g)$.  
\newline

In the rest of this subsection we will seek to understand how 
a given complete contraction in the form (\ref{ambient}) can be 
expressed as a linear combination of complete contractions 
of the form (\ref{karami}), involving 
covariant derivatives of the curvature tensor of the metric $g$.
\newline

{\bf From
contractions in the ambient metric $\tilde{g}$ to 
contractions in the base metric $g$:} 

We do this
in steps. Consider any complete contraction 
$C(\tilde{g}^{n+2})$ in the form (\ref{ambient}), of weight $-n$. 
We denote by $\tilde{\nabla}$ the Levi-Civita connection of $\tilde{g}$ 
and by $\nabla$ the Levi-Civita connection of $g$.

Our aim is to write
$C(\tilde{g}^{n+2})$ as a linear combination
of complete contractions ({\it with respect to the metric $g$}) of the form:

\begin{equation}
\label{prapo}
contr_g(T^{\alpha_1}\otimes\dots\otimes T^{\alpha_y}\otimes F^{\beta_1}\otimes F^{\beta_w}),
\end{equation}
where the factors $T^\alpha$ are all in the form $\nabla^{(m)}W$ (so each $T^\alpha$ is an 
iterated covariant
derivative of the Weyl tensor), all of whose $(m+4)$ indices are being contracted against another index
in $contr_g(\dots)$. Each factor $F^\beta$ is in the form
$\nabla^{(a)}\tilde{\nabla}^{(q)}\tilde{R}$,\footnote{$\nabla$ is
the Levi-Civita connection of $g$ and $\tilde{\nabla}$
 is the Levi-Civita connection of $\tilde{g}$.} where
 $\delta>0$ of its indices have a (fixed) value ${}_\infty$, $\epsilon\ge 0$ of its indices
have a value ${}_0$ and the rest of its indices are being contracted
({\it with respect to the metric $g$})
 against some other index in (\ref{prapo}).
\newline

\par We discuss how any complete contraction $contr_{\tilde{g}}
(\tilde{\nabla}^{(m_1)}\tilde{R}\otimes\dots\otimes\tilde{\nabla}^{(m_p)}
\tilde{R})$ can be written as a linear combination of contractions
in the form (\ref{prapo}). This is a two step-procedure:

 Consider any complete contraction
$C(\tilde{g}^{n+2})$ in the ambient metric. Firstly we define the
set of {\it assignments}, $ASSIGN$:

\begin{definition}
\label{assign}
 An element $ass\in ASSIGN$ is
a rule that acts by picking out each  particular contraction between indices 
$({}_a,{}_b)$ and assigning to that pair of indices
 either the values $({}_\infty,{}_0)$ or
the values $({}_0,{}_\infty)$ or repeatedly assigning any pair of
numbers ${}_i,{}_j$, $1\le i,j\le n$ and then multiplying by
$g^{ij}$ and summing over all such pairs. For
each assignment $ass\in ASSIGN$ we obtain a complete contraction (in the metric $g$)
involving tensors $\tilde{\nabla}^{(m)}\tilde{R}$;  we  denote this complete contraction 
 $ass[C(\tilde{g}^{n+2})]$.
\end{definition}

\par Thus, for each element $ass\in ASSIGN$  
$ass[C(\tilde{g}^{n+2})]$ is a complete contraction in the
quadratic form $g^{ij}$ (denote it by $contr_g(T^1\otimes\dots\otimes T^s)$)
where the factors $T^i$ are tensors
$\tilde{\nabla}^{(m)}\tilde{R}$ with some indices having the
(fixed) value ${}_{\infty}$, some having the (fixed) value ${}_0$
and all the rest of the indices are being contracted against each other with respect to the metric $g^{ij}$.
(Note: when we separately consider a tensor $T^i$ in the contraction $contr_g(T^1\otimes\dots\otimes T^s)$
we will call the indices of the third kind above ``free'' indices; 
in $contr_g(T^1\otimes\dots\otimes T^s)$
these can only be assigned values ${}_1,\dots,{}_n$). 
 It follows that we can write:

\begin{equation}
\label{expect}
C(\tilde{g}^{n+2})=\sum_{ass\in ASSIGN} ass[C(\tilde{g}^{n+2})].
\end{equation}

\par Now, the next step is to write out each $ass[C(\tilde{g}^{n+2})]$
 as a linear combination of contractions in the form (\ref{prapo}) (modulo a linear
combination of contractions with length $\ge\sigma +1$ which
we do not care about). In order to do that, we will pick out each
 tensor $T=\tilde{\nabla}^{(m)}_{r_1\dots r_m}\tilde{R}_{ijkl}$
in $ass[C(\tilde{g}^{n+2})]$ (recall that each of the indices
 ${}_{r_1},\dots ,{}_{r_m},{}_i,{}_j,{}_k,{}_l$  either has a (fixed) value ${}_0$ or ${}_\infty$, or is a
free index that may take any values between
 ${}_1,\dots ,{}_n$). We denote by $\{{}_{u_1},\dots ,{}_{u_t}\}\subset
\{{}_{r_1},\dots ,{}_{r_m},{}_i,{}_j,{}_k,{}_l\}$ the set of free indices in $T$ (i.e. the indices that
 may take values between ${}_1,\dots ,{}_n$).
We will write $T_{u_1\dots u_t}$ to stress the fact that
 the free indices in $T$ that are being contracted (with respect to the
 metric $g^{ij}$) against other indices in $ass[C(\tilde{g}^{n+2})]$
are precisely the indices ${}_{u_1},\dots ,{}_{u_t}$. 
\newline
{\it Step 2:} We claim that $T_{u_1\dots u_t}$ can be expressed in the form:

\begin{equation}
\label{breakT} T_{u_1\dots u_t}=\Sum_{f=1}^F a_f T^f_{u_1\dots
u_t}(g),
\end{equation}
where each of the  terms $T^f_{u_1\dots u_t}(g)$\footnote{Here ${}^f$ is 
a label, while ${}_{u_1},\dots,{}_{u_t}$ are free indices.} is a
 tensor product with free indices ${}_{u_1},\dots ,{}_{u_t}$ in one of two
forms: Either $T^f_{u_1\dots u_t}(g)$ will be a tensor 
product:\footnote{Recall $\nabla$ is the Levi-Civita connection of $g$, and
$\tilde{\nabla}$ is the Levi-Civita connection of $\tilde{g}^{n+2}$}

\begin{equation}
\label{kefalo} \nabla^{(b)}_{c_1\dots c_b}\tilde{\nabla}^{(q)}_{
f_1\dots f_q} \tilde{R}_{zxcv}\otimes g_{ab}\otimes\dots\otimes
g_{a'b'}
\end{equation}
{\it for which at least one index in the tensor
$\tilde{\nabla}^{(q)}\tilde{R}_{ijkl}$ has the value ${}_\infty$},
 or it will just be a tensor of the form:

\begin{equation}
\label{kefalo2} \nabla^{(d)}_{f_1\dots f_d} W_{zxcv}.
\end{equation}

{\it Proof of (\ref{breakT}):} Firstly, if $T_{u_1\dots u_t}$ has
an index with value ${}_\infty$, there is nothing to prove, since we are
in the first of the two desired forms. Thus, we only have to study the
case where all the indices ${}_{r_1},\dots, ,{}_{r_{m+4}}$ have
values between ${}_0$ and ${}_n$. We will then show by induction how
$T_{u_1\dots u_t}$ can be written as a linear combination in the
form (\ref{breakT}): We assume that we are able to write out any
factor $T_{u_1\dots u_t} (=\tilde{\nabla}^{(m)}\tilde{R})$ as in
(\ref{breakT}) provided $m\le d$. We now show how we can do
 this for $m=d+1$:

\par We consider $T_{u_1\dots u_t}=\tilde{\nabla}^{(d+1)}_{r_1\dots r_{d+1}}
\tilde{R}_{r_{d+2}r_{d+3}r_{d+4}r_{d+5}}$ and we
 distinguish cases based on the {\it value} of the index ${}_{r_1}$:
 If ${}_{r_1}$ has the value ${}_0$, we then just denote by $d^{*}$ the number of the indices
${}_{r_2},\dots ,{}_{r_{d+5}}$  that do {\it not} have the value ${}_0$. It then
follows that $T_{u_1\dots
u_t}=(-d^{*}+2)\tilde{\nabla}^{(d)}_{r_2\dots
r_d}\tilde{R}_{r_{d+2}r_{d+3}r_{d+4}r_{d+5}}$. Thus we are done,
 by our inductive hypothesis.

 \par On the other hand,
if the index ${}_{r_1}$ is a free index (allowed to take values
${}_1,\dots ,{}_n$), we then denote by $\{ a_1,\dots a_z\}\subset
\{2,\dots d+5\}$ the set of numbers for which the index
${}_{r_k}$, $k\in \{a_1,\dots a_z\}$ has been assigned the value
${}_0$. We also denote by $\{b_1\dots ,b_x \}\subset \{2,\dots d+5\}$
the set of numbers for which ${}_{r_k}$, $k\in \{b_1\dots ,b_x \}$
 is a free index (taking values between ${}_1,\dots,{}_n$). We then have:

\begin{equation}
\label{breakTwell} \begin{split} & \tilde{\nabla}^{(d+1)}_{r_1\dots
r_{d+1}}\tilde{R}_{r_{d+2}r_{d+3}r_{d+4}r_{d+5}}=\nabla_{r_1}(\tilde{\nabla}^{(d)}_{r_2\dots
r_{d+1}}\tilde{R}_{r_{d+2}r_{d+3}r_{d+4}r_{d+5}})
\\&-\Sum_{k=1}^z\tilde{\nabla}^{(d)}_{r_2\dots
r_{a_k-1}r_1r_{a_k+1} r_{d+1}}\tilde{R}_{r_{d+2}r_{d+3}r_{d+4}r_{d+5}}
\\&+\Sum_{k=1}^x\tilde{\nabla}^{(d)}_{r_2\dots r_{b_k-1}\infty
r_{b_k+1} r_{d+1}}\tilde{R}_{r_{d+2}r_{d+3}r_{d+4}r_{d+5}}\otimes
g_{r_1r_{b_k}}+Q(R).
\end{split}
\end{equation}
 (\ref{breakTwell}) follows virtue of the
Christoffel symbols  $\tilde{\Gamma}_{0i}^k=\delta^k_i$ and
$\tilde{\Gamma}_{ij}^\infty=-g_{ij}$.

\par Now, observe that the tensors $\tilde{\nabla}^{(d)}_{r_2\dots r_{b_k-1}\infty
r_{b_k+1} r_{d+1}}\tilde{R}_{r_{d+2}r_{d+3}r_{d+4}r_{d+5}}$ are in our desired
form (\ref{kefalo}) because they contain an index ${}_\infty$.
 Moreover, the tensors $(\tilde{\nabla}^{(d)}_{r_2\dots
r_{d+1}}\tilde{R}_{r_{d+2}r_{d+3}r_{d+4}r_{d+5}})$ and
$\tilde{\nabla}^{(d)}_{r_2\dots r_{a_k-1}r_1r_{a_k+1}
r_{d+1}}\tilde{R}_{r_{d+2}r_{d+3}r_{d+4}r_{d+5}}$ fall under our
inductive hypothesis. Therefore, expressing these tensors  as
linear combinations of tensors in the form (\ref{kefalo}),
(\ref{kefalo2}) we derive that (\ref{breakT}) follows by induction. $\Box$
\newline

\par This analysis of each of the tensors $T_{u_1\dots u_t}$
allows us to analogously write each complete contraction
$C(\tilde{g}^{n+2})$  as a sum of complete
contractions in the form (\ref{prapo}): We pick out each factor
$T^y$ ($y=1,\dots ,\sigma$) and replace it by one of the summands
$T^y_r$ (times a coefficient) in (\ref{breakT}). We then take the
complete contraction (with respect to $g^{ij}$) of these new tensors $T^y_f$. The
formal sum over all these substitutions will be denoted by
$\Sum_{\gamma\in BR_{ass}} \gamma\{ C(\tilde{g}^{n+2})\}$.
Therefore, using the above and (\ref{expect}) we derive a formula:

\begin{equation}
\label{definebreak} C(\tilde{g}^{n+2})=\sum_{ass\in ASSIGN}
\sum_{\gamma\in BR_{ass}} \gamma\{ ass[C(\tilde{g}^{n+2})]\};
\end{equation}
(the symbol $BR_{ass}$ serves to illustrate that we have first
picked a particular assignment $ass\in ASSIGN$).
\newline

Finally, we introduce a definition. 

.

\begin{definition}
\label{amban}
\par For each complete contraction $C^l(g)$ in the form:
\begin{equation}
\label{noq} contr({\nabla}_{r_1\dots
r_{m_1}}^{(m_1)}W_{i_1j_1k_1l_1}\otimes \dots \otimes
{\nabla}_{v_1\dots v_{m_s}}^{(m_s)}W_{i_sj_sk_sl_s}), 
\end{equation}
with each $m_a\le \frac{n-4}{2}$,
 we construct a  complete contraction $Amb[C^l(g)]$ in the
ambient metric:

\begin{equation}
\label{ambient} contr({\tilde{\nabla}}_{r_1\dots
r_{m_1}}^{(m_1)}{\tilde{R}}_{i_1j_1k_1l_1} \otimes \dots \otimes
{\tilde{\nabla}}_{v_1\dots
v_{m_s}}^{(m_s)}{\tilde{R}}_{i_sj_sk_sl_s});
\end{equation}
 this is related to $C^l(g)$ in the following way: Consider
 any two
factors $T_a={\nabla}^{(m_w)}_{x_1\dots x_{m_w}} W_{i_wj_wk_wl_w}$
and $T_b={\nabla}^{(m_z)}_{y_1\dots y_{m_z}}W_{i_zj_zk_zl_z}$ in
$C^l(g)$ (the $a^{th}$ and $b^{th}$ factors)
and suppose that the $k^{th}$ index in $T_1$
contracts against the $l^{th}$ index in $T_2$.
Then, we require that in $Amb[C^l]$ we will have two factors
 ${\tilde{T}}_a={\tilde{\nabla}}^{(m_w)}_{x_1\dots x_{m_w}}
{\tilde{R}}_{i_wj_wk_wl_w}$ and
${\tilde{T}}_b={\tilde{\nabla}}^{(m_z)}_{y_1\dots y_{m_z}}
{\tilde{R}}_{i_zj_zk_zl_z}$ in $C^l({\tilde{g}}^{n+2})$ and furthermore
 the $k^{th}$ index
in ${\tilde{T}}_1$ will contract against the $l^{th}$ index in
${\tilde{T}}_2$. 

We will call the local conformal invariant $Amb[C(g)]$
 the ambient analogue of $C(g)$. 
\end{definition}
{\it Note:} Since we are assuming $m_a\le \frac{n-4}{2}$, we know by
\cite{fg:ci}, \cite{fg:latest}
 that the above ambient complete contraction is
well-defined.

\subsection{Proof of Lemma \ref{subcinv}:}

 In order to prove this Lemma,  let us consider the
 linear combination
${\Sigma}_{l\in L^0_\sigma} a_l C^l(g)$. 
We claim that the local conformal invariant 
needed for the proof of Lemma \ref{subcinv} is precisely 
$W(g):=\sum_{l\in L} a_l Amb[C^l(g)]$.\footnote{See Definition \ref{amban}.}
 This will follow by virtue of the following Lemma:

\begin{lemma}
\label{cancelweyl} Given any $C^l(g)$ in the form (\ref{karami}) with no internal contractions, we claim:
$$C^l(g)- Amb[C^l]=\sum_{t\in T} a_t C^t(g),$$
where each $C^t(g)$ either  has length $\ge
\sigma +1$ (and is in the form (\ref{contraction})), 
or is in the form (\ref{karami}) and has length $\sigma$
and $\delta_W\ge 1$.
\end{lemma}

{\it Proof:} Recall the algorithm
 from the previous subsection. Recall the equation (\ref{definebreak}). 
Let  $ass^{*}\in ASSIGN$ be the 
(unique) assignment where {\it no} index in the complete contraction $Amb[C^l(g)]$
is assigned the value $\infty$.\footnote{Equivalently, $ass^{*}$ is the unique 
assignment where no particular contraction $({}_a,{}_b)$ is 
assigned the values $({}_\infty,{}_0)$ or $({}_0,{}_\infty)$.}
Firstly we claim that for every $ass\in ASSIGN\setminus \{ass^{*}\}$:

\begin{equation}
\label{fair}
\sum_{\gamma\in BR_{ass}} \gamma\{ ass[C(\tilde{g}^{n+2})]\} =\sum_{t\in T} a_t C^t(g),
\end{equation}
where the RHS is as in the claim of Lemma \ref{cancelweyl}. 
To see this, recall the algorithm from the previous subsection. Observe that 
$ass[C(\tilde{g}^{n+2})]$ stands for a complete contraction 
(with respect to the metric $g^{ij}$) of factors $\tilde{\nabla}^{(m)}_{r_1\dots r_m}\tilde{R}_{ijkl}$,
 {\it and at least one of those factors has an index with the (fixed) value ${}_\infty$};
it then follows from the higher-order Taylor expansion of the ambient metric $\tilde{g}$ 
(see \cite{fg:ci}, \cite{fg:latest})
 that such a tensor $T_{i_1\dots i_f}$,
  can be written out in terms of the metric $g$ as follows: 

$$T_{i_1\dots i_f}= \sum{\nabla}^{(m)}W_{ijkl} +
{\Sum} T'_{i_1\dots i_f},$$ where each tensor ${\nabla}^{(m)}
W_{ijkl}$ has $f$ free lower
 indices, $m+4>f$ and the rest of the indices are internally
 contracting. Moreover, {\it there will be} at least one internal contraction
  in $\nabla^{(m)}W_{ijkl}$. 
$\sum T'_{i_1\dots i_f}$ stands for a linear combination
 of partial contractions in the form:

\begin{equation}
\label{pcontr} pcontr({\nabla}^{(m_1)}R_{ijkl}\otimes\dots\otimes
{\nabla}^{(m_w)} R_{i'j'k'l'}\otimes g_{ab}\otimes\dots\otimes g_{a'b'})
\end{equation}
with $w\ge 2$ (i.e. with at least two factors $\nabla^{(m)}R$). 
Substituting this expression into $\sum_{\gamma\in BR_{ass}} \gamma\{ ass[C(\tilde{g}^{n+2})]\}$
we derive (\ref{fair}).

\par So we are reduced to considering 
$\sum_{\gamma\in BR_{ass^{*}}} \gamma\{ ass^{*}[C(\tilde{g}^{n+2})]\}$.
 That is, we are reduced to considering the case where {\it no index} in the 
complete contraction $C^l(\tilde{g}^{n+2})$ is assigned the value ${}_0$ or ${}_\infty$. 
 But then recall the Christoffel symbols
$\tilde{\Gamma}^k_{ij}$ for the ambient metric $\tilde{g}^{n+2}$:

\par If $1\le i,j,k\le n$ we have $\tilde{\Gamma}^k_{ij}=
{\Gamma}^k_{ij}$. On the other hand, $\tilde{\Gamma}^0_{ij}=
P_{ij}$ and $\tilde{\Gamma}^{n+1}_{ij}=-g_{ij}$.
Thus we derive that for values ${}_{r_1},\dots,{}_l\in \{{}_1,\dots,{}_n\}$:

\begin{equation}
\label{dichoamb} \tilde{\nabla}^{(m)}_{r_1\dots
r_m}\tilde{R}_{ijkl}= {\nabla}^{(m)}_{r_1\dots r_m}W_{ijkl}+
T_{r_1\dots r_mijkl}+ T'_{r_1\dots r_mijkl},
\end{equation}
where $ T_{r_1\dots r_mijkl}$ is a linear combination of tensor
fields of the form $g_{ab}\otimes\dots\otimes g_{cd}\otimes
\tilde{\nabla}^{(p)}\tilde{R}_{ijkl}$, where the factor
$\tilde{\nabla}^{(p)}\tilde{R}_{ijkl}$ contains at least one index with a value
${}_\infty$. $T'_{r_1\dots r_mijkl}$ stands for a linear combination
of tensor fields of the form (\ref{pcontr}) with $w\ge 2$.

\par This completes the proof of Lemma \ref{cancelweyl}. $\Box$

\subsection{Proof of Lemma \ref{1stclaim}: Proof of half the Lemma.}

\par We prove Lemma \ref{1stclaim} 
 in two steps. In order to explain this proof, we need
one piece of notation:

\begin{definition}
\label{splitlapl} 
 Consider a complete contraction $C(g)$
in the form (\ref{karami}) {\it with weight $-n+2$}. We define
$\Delta_r^k[C(g)]$, for any $k, 1\le k\le \sigma$ to be the
complete contraction that arises from $C(g)$ by picking out its
$k^{th}$ factor $F_k$ and replacing it by $\Delta F_k$. Then, we
define:

$$\Delta_r[C(g)]=\Sum_{k=1}^\sigma \Delta_r^k[C(g)]$$
This operation extends to linear combinations.
\end{definition}

\par Lemma \ref{1stclaim} will follow by the next two claims:

\begin{lemma}
\label{palrights1}
 Under the assumptions of Lemma \ref{1stclaim}
 we claim that there is a vector field $T^i(g)=\Sum_{r\in R}
 a_r C^{r,i}(g)$ (of the form required by Lemma \ref{1stclaim}) so that:

 \begin{equation}
 \label{ekratos}
\Sum_{l\in L^1_\sigma} a_l C^l(g)-div_i\Sum_{r\in R}
 a_r C^{r,i}(g)=\Delta_r[\Sum_{v\in V} a_v C^v(g)]+\Sum_{f\in
 F}a_f C^f(g),
 \end{equation}
 modulo complete contractions of length $\ge\sigma +1$. Here each $C^f(g)$
 is in the form (\ref{karami}) with length $\sigma$ and $\delta_W\ge 2$.
  On the other hand, each $C^v(g)$ is in the
 form (\ref{karami}) with length $\sigma$, weight $-n+2$ and
with no internal contractions. 
\end{lemma}

\par The next claim starts with the conclusion of the previous
one:
\begin{lemma}
\label{definit} 
Assume that $P(g)$ is as in the assumption of Lemma 
\ref{1stclaim}; assume  furthermore that 
the sublinear combination $\sum_{l\in L^1} a_l C^l(g)$
 of contractions with length $\sigma$ and precisely
  one internal contraction can be expressed in the form
$\sum_{l\in L^1} a_l C^l(g)=\sum_{v\in V} a_v \Delta_r[C^v(g)]$.

We claim that  there
is a local conformal invariant $W(g)$ and a vector field
$\Sum_{r\in R}  a_r C^{r,i}(g)$ of weight $-n+1$ so that:

\begin{equation}
\label{karter} \Delta_r[\Sum_{v\in V} a_v
C^v(g)]-W(g)-div_i\Sum_{r\in R}  a_r C^{r,i}(g)=\Sum_{f\in
 F}a_f C^f(g).
\end{equation}
The above holds modulo complete contractions of length $\sigma
+1$. Here again each $C^f(g)$
 is in the form (\ref{karami}) with length $\sigma$ and $\delta_W\ge 2$.
\end{lemma}

{\it Proof of Lemma \ref{palrights1}:} In order to prove this
Lemma, we first slightly manipulate the sublinear 
combination $\sum_{l\in L_\sigma^1} a_l C^l(g)$ in $P(g)$:

\par We first prove that by subtracting divergences from
$P(g)$, we may assume that every complete contraction $C^l(g)$
with $l\in L^1_\sigma$ will have its internal contraction in a
factor of the form ${\nabla}_l{W_{ijk}}^l$:

\par This is done as follows: For each complete
 contraction $C^l(g)$, $l\in L^1_\sigma$, we isolate
the one factor ${\nabla}^{(m_a)}_{r_1\dots r_{m_a}}W_{ijkl}$ which
contains the internal contraction (there can be only one such
factor, by the definition of $L^1_\sigma$). Then, if the internal
contraction is between two indices ${}_{r_s}, {}_{r_t}$, we bring
them to the positions ${}_{r_{m_a-1}}, {}_{r_{m_a}}$ by repeatedly applying
the curvature identity. That way we introduce
 correction terms that have length at least $\sigma +1$.
Then, we apply the identity
\begin{equation} \label{in&out}
\nabla_rW_{ijkl}+\nabla_jW_{rikl}+\nabla_iW_{jrkl}=\sum
(\nabla^sW_{srty}\otimes g),
\end{equation}
(see subsection 2.3 \cite{alexakis1}--recall that  the symbol $\sum(\nabla^sW_{srty}\otimes g)$ stands for a
linear combination of a tensor product of the three-tensor
$\nabla^sW_{sqty}$ with an
un-contracted metric tensor).

 Thus (modulo introducing correction terms that are allowed in
Lemma \ref{palrights1}),
 we may assume that in each complete contraction $C^l(g)$, $l\in L^1$ the internally
 contracting factor is of the form
${{\nabla}^{(m)}_{r_1\dots r_{m-1}l}}{W_{ijk}}^l$. Finally, we
 subtract divergences from $C^l(g)$ as in the proof of the
silly divergence formula in \cite{a:dgciI}, in order to arrange that all complete
contractions $C^l(g), l\in L^1_\sigma$ have the internal
 contraction in a factor of the form ${\nabla}_l{W_{ijk}}^l$ (this can be
  done modulo introducing complete contractions with more than one internal contractions--but these
  are allowed in the conclusion of our Lemma).

  \par We note the
transformation law for the factor ${\nabla}^lW_{ijkl}$ under
 the re-scaling $\hat{g}=e^{2\phi (x)} g$.

\begin{equation}
\label{trangw}
({\nabla}^lW_{ijkl})_{\hat{g}}=({\nabla}^lW_{ijkl})_{g} + (n-3)
[ W_{ijkl}{\nabla}^l\phi]_{g}.
\end{equation}

\par Now, consider $I^1_{g}(\phi)(:=Image^1_\phi[P(g)]$).\footnote{Recall
 that $Image^1_\phi[P(g)]:=\frac{d}{dt}|_{t=0}e^{nt\phi}P(e^{2t\phi}g)$ 
and that $\int_{M^n}I^1_g(\phi)dV_g=0$.} 
We firstly study
the sublinear combination $Image^1_{\phi}[P(g)]|_\sigma]$.\footnote{Recall 
that $P(g)|_\sigma$ stands for the sublinear combination of terms in $P(g)$ 
with length $\sigma$. Those terms are indexed in $L_\sigma\subset L$. Recall 
that  $Image^1[P(g)_\sigma](:=)\frac{d}{dt}|_{t=0}e^{nt\phi}\sum_{l\in L_\sigma} a_l C^l(e^{2t\phi}g)$.}

Initially, we write $Image^1_{\phi}[P(g)|_\sigma]$ out as a linear combination of
complete contractions in the form:

\begin{equation}
\label{normcontrphi}
contr(\nabla^{(m_1)}W_{ijkl}\otimes\dots\otimes\nabla^{(m_\sigma)}
W_{i'j'k'l'}\otimes\nabla^{(\nu)}\phi).
\end{equation}

 Given that
$P(g)|_\sigma$ is a linear combination of complete contractions in
the from (\ref{karami}), we immediately obtain such an expression
for $Image^1_{\phi}[P(g)|_\sigma]$, by applying the
 transformation laws from the Weyl tensor and the Levi-Civita connection, 
see the subsection 2.3 in \cite{alexakis1}. Thus
we write

$$Image^1_{\phi}[P(g)|_\sigma]={\Sigma}_{k\in K} a_k C^k_{g}(\phi),$$
where each $C^k_{g}(\phi)$ is in the form (\ref{normcontrphi})
with length $\sigma+1$.

 \par Now, we
 re-express $Image^1_{\phi}[P(g)|_\sigma]$ as a linear
 combination of complete contractions in the form:

 \begin{equation}
 \label{linisymric'}
contr(\nabla^{(m_1)}R\otimes\dots\otimes\nabla^{(m_\sigma)}R\otimes\nabla^{(\nu)}\phi)
\end{equation}
($\nabla^{(m)}R$ above stands for the differentiated curvature
tensor, where we do not write out the indices). This is done by
picking out each complete contraction $C^k_{g}(\phi)$ in the above
equation and decomposing
  each Weyl tensor according to $W_{ijkl}=R_{ijkl}+[P\wedge g]$ (see 
 subsection 2.3 in \cite{alexakis1}).
 Hence, we write out:

$$Image^1_{\phi}[P(g)|_\sigma]={\Sigma}_{k\in K} a_k
[{\Sigma}_{w\in W_k}a_wC^{k,w}_{g}(\phi)].$$

Each complete contraction $C^{k,w}_{g}(\phi)$ is in the form
(\ref{linisymric'}).
\newline

\par Now, we introduce some definitions regarding complete contractions
$C_{g}(\phi)$ of the form (\ref{linisymric'}).

\begin{definition}
\label{deltaq}
 For any complete contraction $C_{g}(\phi)$,
$\delta$ will stand for the number of internal contractions
 among its factors
$\nabla^{(m)}R_{ijkl}$, $\nabla^{(p)}Ric$, $\nabla^{(y)}R$
(including the one (resp. two) internal contractions in the term
$Ric_{ij}={R^a}_{iaj}$, $R={R^{ab}}_{ab}$). $q$ will stand for the
number of factors $\nabla^{(p)}Ric$, $R$ (scalar curvature).
\end{definition}

\begin{definition}
\label{targetc}

\par Consider complete contractions $C_{g}(\phi)$ in the
 form (\ref{linisymric'}). $\nu$ will stand for the number of
 derivatives on the factor $\nabla^{(\nu)}\phi$.

\par Any such complete
 contraction with $q=\delta=0$ and $\nu =1$
will be called a target. Any complete contraction with $\nu \ge 2$
and ${\nabla}^{(\nu)}_{r_1\dots r_{\nu}}\phi\ne \Delta \phi$ will
be called irrelevant. Any complete contraction with $\nu =1$ and
 $q+\delta>0$ will be called a contributor.

\par Now, we consider complete contractions in the
 form (\ref{linisymric'}) of length $\sigma +1$ with a factor
$\Delta \phi$.  If $q=\delta =0$, we will
 call it dangerous. If $\delta +q>0$ we will call it a
 contributor. Finally, any complete contraction 
of length $\ge \sigma +2$ will be
called irrelevant.
\end{definition}

\par In the rest of this subsection, $J_{g}(\phi)$ will stand for a
 generic linear combination of contributors and irrelevant
 complete contractions.

\par Our next aim is to understand $Image^1_\phi[C^l(g)]$ for $l\in L^1_\sigma$.
  We will need some definitions:

\begin{definition}
\label{skeletons}
\par Consider a complete contraction $C^l_{g}(\phi)$,
$l\in L^1_\sigma$ which is in the form (\ref{karami}) and has its
internal contraction in a factor ${\nabla}^lW_{ijkl}$. We define
the skeleton $C^{l,\iota}_{g}(\phi)$ of $C^l(g)$ to be the
complete contraction which is obtained from $C^l(g)$ by
substituting the factor ${\nabla}_l{W_{ijk}}^l$ by
 ${\nabla}_l\phi{R_{ijk}}^l$ and every other factor
${\nabla}^{(m)}_{r_1\dots r_m}W_{ijkl}$ by ${\nabla}^{(m)}_{r_1\dots
r_m}R_{ijkl}$.

\par Also, consider any complete contraction $C^l(g)$ which is in the form (\ref{karami})
(not necessarily of weight $-n$) with $q=0$ and $\delta=0$. We
then define its skeleton $C^{l,\iota}(g)$ to be the
 complete contraction which is obtained from it by
 substituting each factor
${\nabla}^{(m)}_{r_1\dots r_m}W_{ijkl}$ by ${\nabla}^{(m)}_{r_1\dots
r_m}R_{ijkl}$.
\end{definition}

\begin{lemma}
\label{y1targ} For any complete contraction $C^l(g)$ with $l\in
L^1_\sigma$, we claim that $Image_{\phi}^1[C^l(g)]$ can be expressed as follows:

\begin{equation}
\label{specyitarg} Image_{\phi}^1[C^l(g)]=
(n-3)C^{l,\iota}_{g}(\phi)+ J_{g}(\phi).
\end{equation}
\end{lemma}

{\it Proof:} This is proven in two steps:

\par We first consider the complete contraction
in $Image_{\phi}^1[C^l(g)]$ which arises as follows:
 We first substitute the factor
${\nabla}^lW_{ijkl}$ in $C^l(g)$ by $(n-3)W_{ijkl}{\nabla}^l\phi$.
Let us denote the complete contraction that we obtain thus by
$C^l_{g}(\phi)$. Then, we write $C^l_{g}(\phi)$ as a linear
combination of complete contractions in the form
 (\ref{linisymric'}):

$$C^l_{g}(\phi)= {\Sigma}^F_{t=1} a_t C^{l,t}_{g}(\phi)$$

Let us assume that $C^{l,1}_{g}(\phi)$ is obtained from
$C^l_{g}(\phi)$ by substituting each factor
${\nabla}^{(m)}_{r_1\dots r_m}W_{ijkl}$ by a factor
 ${\nabla}^{(m)}_{r_1\dots r_m}R_{ijkl}$.
The complete contractions $C^{l,t}_{g}(\phi)$ with $t\ge 2$ arise
by substituting at least one factor ${\nabla}^{(m)}W_{ijkl}$ in
$C^l_{g}(\phi)$ by a factor ${\nabla}^{(m)}[Ric\otimes g]$ 
or $\nabla^{(m)}[R\otimes g\dots \otimes g]$. Hence,
each complete contraction $C^{l,t}_{g}(\phi)$ with $t\ge 2$
 will either have $q>0$ or $\delta>0$,
so it will be a contributor.
We see that the complete contraction $C^{l,1}_{g}(\phi)$ above is
$(n-3)\cdot C^{l,\iota}_{g} (\phi)$.

\par Let us now consider any 
complete contraction $C^{l'}_{g}(\phi)$ (in the form)
(\ref{normcontrphi}) in $Image_{\phi}^1[C^l(g)]$ other than $C^{l,1}_g(\phi)$. Then necessarily 
$C^{l'}_{g}(\phi)$ has arisen by
 applying the transformation laws for the Levi-Civita connection or the Weyl tensor 
to any indices in $C^l(g)$ other than the internal
 contraction in the factor
${\nabla}^lW_{ijkl}$. Hence $C^{l'}_{g}(\phi)$ will contain a
factor ${\nabla}^lW_{ijkl}$, which still has an internal contraction.
Therefore, writing $C^{l'}_{g}(\phi)$ as a linear combination of
complete contractions in the form (\ref{linisymric'}), as below:

$$C^{l'}_{g}(\phi)= {\Sigma}_{r\in R} C^{l',r}_{g}(\phi)$$
we have that each $C^{l',r}_{g}(\phi)$ will either have a factor
${\nabla}^{(m)}R_{ijkl}$
 with an internal contraction or a
 factor ${\nabla}^{(p)}Ric$ or a factor $R$ (scalar curvature).
Therefore, each such complete
 contraction $C^{l',r}_{g}(\phi)$
 is either a contributor or irrelevant. $\Box$

\par By the same calculations we also derive:

\begin{lemma}
\label{ylnotarg} For any complete contraction $C^l(g)$ with $l\in
L_\sigma\setminus L^1_\sigma$:

\begin{equation}
\label{specylnotarg} Image_{\phi}^1[C^l(g)]= J_{g}(\phi).
\end{equation}
\end{lemma}

{\it Proof:} This fact follows by the same proof as for the
previous Lemma. $\Box$

\begin{lemma}
\label{badtoir} Consider the sublinear combination $P(g)|_{\sigma
+1}$ in $P(g)$.  Then:

\begin{equation}
\label{dangtarg} Image_{\phi}^1[P(g)|_{\sigma +1}]= \Sum_{u\in U}
a_u C^u(g)\Delta\phi+ J_{g}(\phi),
\end{equation}
where $\Sum_{u\in U} a_u C^u(g)\Delta\phi$ stands for a generic
 linear combination of dangerous complete contractions.
\end{lemma}

{\it Proof:} Straightforward from the transforation laws of 
the Levi-Civita connection and the curvature tensor under 
conformal re-scalings. (These can be found in  subsection 2.3 
in \cite{alexakis1}). $\Box$
\newline

\par We are now ready to get to the main part of proving
 Lemma \ref{palrights1}. It will be useful to recall a few facts about the super divergence
 formula from \cite{a:dgciI}. This is the second and last instance 
in this series of papers where we make use of the super divergence 
formula in full strength. In all other instances we use the 
``main consequence'' of the super divergence formula, as codified in subsection 2.2 in \cite{alexakis1}. 

 {\it A few facts about the super divergence formula:} We apply the super divergence formula
to the operator $I^1_g(\phi)=\sum_{l\in L} a_l C^l(g)$. Recall 
that $\int_{M^n} I^1_g(\phi)dV_g=0$ for every 
compact $(M^n,g)$, $\phi\in C^\infty (M^n)$ and each $C^l(g)$ is in the form (\ref{linisymric'}). 
The super divergence formula applied to $I^1_g(\phi)$ provides a local formula which expresses 
$I^1_g(\phi)$ as a divergence of a vector field. We recall 
that there is a process by which
 each term $C^l(g)$ in $I^1_g(\phi)$ gives rise to  divergences in
the super divergence formula. In the end of \cite{a:dgciI} we 
 summarized the conclusion of the super divergence formula as follows: 
Given $I^1_g(\phi)=\sum_{l\in L} a_l C^l_g(\phi)$ we introduces 
the notion of  ``descendents'' of 
each $C^l(g),l\in L$; these are 
 complete contractions involving factors $\vec{\xi}$ (or covariant
 derivatives thereof). Such complete contractions are divided up into categories 
(e.g. ``good'', ``hard'', ``undecided'') and we then proceed by 
integrating by parts the factors $\vec{\xi}$ in those complete contractions. 
At each stage,  we {\it discard} 
complete contractions which are ``bad'', ``hard'' or ``stigmatized''. In the end 
each complete contraction $C^l_g(\phi)$ in 
$I^1_g(\phi)$ contributes itself plus a linear combination of 
divergences to the super divergence formula; we have 
denoted this sum by $Tail[C^l(g)]$. In other words, the 
super diveregnce formula can be summarized as:
\begin{equation}
\label{prosfyges}
\sum_{l\in L} a_l Tail[C^l(g)]=0. 
\end{equation}
The only further remarks we wish to make is that 
if $C^l_g(\phi)$ has length $\ge\sigma+2$ 
then $Tail[C^l_g(\phi)$ consists of complete contractions 
with length $\ge\sigma+2$, and that  for every complete contraction $C^l_g(\phi)$ of length $\sigma+1$
 the only way that the factor $\nabla^{(\nu)}\phi$ 
in a complete contraction $C^l_g(\phi)$ can give rise to a factor $\vec{\xi}$ is by replacing 
a pair of indices ${}_{r_s},{}_{r_t}$ in the factor $\nabla^{(\nu)}_{r_1\dots r_\nu}\phi$ that 
contract against each other by an expression ${}_{r_t}\vec{\xi}^{r_t}$. 
This factor $\vec{\xi}^{r_t}$ will 
then be integrated by parts, giving rise to a 
divergence with respect to the index ${}_{r_t}$. 

\par Our Lemma will follow  by picking out a particular ``piece'' 
in the super divergence formula {\it which vanishes separately} and then applying the operation
 $Weylify$ (from subsection 5.1 in \cite{alexakis1}) to the resulting equation:

\begin{lemma}
\label{cancirrel} Consider any irrelevant complete
contraction $C^l_{g}(\phi)$. Then:

$$Tail[C^l_{g}(\phi)]={\Sigma}_{w\in W} C^w_{g}(\phi),$$
where each complete contraction $C^w_{g}(\phi)$ has length $\ge
\sigma +1$, and if $C^w_{g}(\phi)$ does have length $\sigma +1$ then
it is irrelevant.
\end{lemma}

{\it Proof:} If $C^l_g(\phi)$ has length $\ge\sigma+2$ then $Tail[C^l_g(\phi)]$ 
consists of terms with length $\ge\sigma+2$, so we are done. Now, 
the case where $C^l_g(\phi)$ has length $\sigma+1$: 
 By our definition, a complete contraction of length
$\sigma +1$ is irrelevant if it has a factor
${\nabla}^{(m)}_{r_1\dots r_m}\phi\ne\Delta\phi$. But then, by
 virtue of the Lemma 16 in \cite{a:dgciI} and the iterative
 integration by parts,  it follows that each complete contraction
 in $Tail[C^l_{g}(\phi)]$ of length
$\sigma +1$ will have a factor ${\nabla}^{(p)}_{r_1\dots r_p}\phi,
p\ge 2$. Hence each complete contraction of length $\sigma +1$ in
$Tail[C^l_{g}(\phi)]$ is irrelevant. $\Box$
\newline

\par We next consider $Tail[C^l_g(\phi)]$ when $C^l_g(\phi)$ is a contributor. 
(By definition if $C^l_{g}(\phi)$ is a
contributor, then its factor ${\nabla}^{(p)}\phi$ is either of the
form $\nabla\phi$ or $\Delta\phi$).

We claim the following:

\begin{lemma}
\label{tailcontr} Consider any contributor $C^l_{g}(\phi)$.
 Modulo complete contractions of length $\ge \sigma +2$, we
 can write $Tail[C^l_{g}(\phi)]$ as a linear combination:

\begin{equation}
\label{tailana1} Tail[C^l_{g}(\phi)] = C^l_{g}(\phi)+
{\Sigma}_{r\in R_l} a_r div_{j_r} C^{j_r}_{g} (\phi) +
{\Sigma}_{w\in W} a_w C^w_{g}(\phi),
\end{equation}
where each vector field $C^{l,j_r}_{g} (\phi)$ in the above
equation is a partial contraction in the form (\ref{linisymric'})
with length $\sigma +1$, with one free index,
 and furthermore
$\nu =1$ (that is, it has a factor ${\nabla}_a\phi$), and the
index ${}_a$ in $\nabla_a\phi$ is not the free index, and
furthermore $q=\delta =0$; each complete contraction
$C^w_{g}(\phi)$  is in the form (\ref{linisymric'}) with
 either $\nu\ge 2$ or $q+\delta >0$.
\end{lemma}

{\it Proof of Lemma \ref{tailcontr}:} First consider
the case where $C^l(g)$ contains a factor $\nabla\phi$: By virtue
of the algorithm for the super divergence
 formula (see the concluding remarks in \cite{a:dgciI}),
we derive that, modulo complete contractions of
 length $\ge \sigma +2$:

$$Tail[C^l_{g}(\phi)]= C^l_{g}(\phi) -
 {\Sigma}_{f\in F} div_{i_f}C^{l,i_f}_{g}(\phi),$$
where ${\Sigma}_{f\in F} C^{l,i_f}_{g}(\phi)$ is a linear combination 
of vector fields in the form (\ref{linisymric'}) with $\nu\ge 1$; if $\nu=1$ then
 the index ${}_a$ in $\nabla_a\phi$ {\it is not} 
the free index ${}_{i_f}$ because $C^l_{g}(\phi)$
contains a factor $\nabla\phi$, therefore no descendent of
$C^l_{g}(\phi)$ has a factor $\vec{\xi}$ contracting against
$\nabla\phi$.
 Now, if a vector field $C^{l,i_f}_{g}(\phi)$
has $\nu =1$ and no internal contractions,
 we  index it in $R_l$. Otherwise,
 we place $div_{i_f}C^{l,i_f}_{g}(\phi)$ into
 the sum ${\Sigma}_{w\in W} a_w C^w_{g}(\phi)$. $\Box$
\newline

\par We now consider the case where $C^l_{g}(\phi)$ has a
 factor $\Delta \phi$. In that case, recall that by definition
$\delta+ q>0$. We denote the set of good or undecided
descendants of $C^l_{g}(\phi)$ by $\{
C^{l,b}_{g}(\phi,\vec{\xi})\}_{b\in B}$. Also recall Lemma
16 from \cite{a:dgciI}. Now, if $C^{l,b}_{g}(\phi,
\vec{\xi})$ contains a factor $\Delta\phi$, it follows from Lemma
20 in \cite{a:dgciI} that, modulo
 complete contractions of length $\ge \sigma +2$:

$$PO[C^{l,b}_{g}(\phi,\vec{\xi})]={\Sigma}_{t\in T} a_t
C^t_{g}(\phi),$$
 where each $C^t_{g}(\phi)$ is in the form (\ref{linisymric'})
 and has a factor
${\nabla}^{(\nu)}\phi,\nu \ge 2$. We then place the complete
contractions   $C^t_{g}(\phi)$ into the sum 
${\Sigma}_{w\in W} a_w C^w_{g}(\phi)$.

\par If the $\vec{\xi}$-contraction $C^{l,b}_{g}(\phi,
\vec{\xi})$ has a factor ${\nabla}_i\phi$, it follows from Lemma
16 in \cite{a:dgciI} that it will have an expression
${\nabla}_i\phi\vec{\xi}^i$. We then decree
that the factor $\vec{\xi}^i$ will be the first to be integrated
by parts.\footnote{Note that by the remark made in
the subsection ``Conclusion'' in \cite{a:dgciI}
we are free to impose this restriction.}
 Notice the following: If $\vec{\xi}^i$ was the only
$\vec{\xi}$-factor in $C^{l,b}_{g}(\phi,\vec{\xi})$ is
$\vec{\xi}^i$, then $C^{l,b}_{g}(\phi,\vec{\xi})$ is in the form
(\ref{linisymric'}) with at least one factor 
$\nabla^{(p)}Ric$ or $R$ (of the scalar curvature).
 Hence, in that case we place
$PO[C^{l,b}_{g}(\phi,\vec{\xi})]$ into the sum ${\Sigma}_{w\in W}
a_w C^w_{g}(\phi)$.

\par If $C^{l,b}_{g}(\phi,\vec{\xi})$
has more $\vec{\xi}$-factors, then after the integration by
 parts of $\vec{\xi}^i$, we are reduced to the previous case.
Thus our Lemma follows. $\Box$
\newline

\par Finally, let us consider a dangerous complete contraction $C^l_{g}(\phi)$. It will be in the form:

$$contr({\nabla}^{(m_1)}R_{ijkl}\otimes\dots\otimes
{\nabla}^{(m_{\sigma})}R_{ijkl}\otimes\Delta\phi),$$
where none of the factors ${\nabla}^{(m)}R_{ijkl}$ have internal
 contractions. Hence, we derive:

\begin{equation}
\label{taildang}
Tail[C^l_{g}(\phi)]=-contr({\nabla}^i[{\nabla}^{(m_1)}
R_{ijkl}\otimes\dots\otimes
{\nabla}^{(m_{\sigma})}R_{ijkl}]\otimes{\nabla}_i\phi).
\end{equation}
We denote a linear combination as in the right hand side of the
above by $\sum_{h\in H} a_h \nabla^j[C^h(g)]\nabla_j\phi$.
\newline

\par In view of the above Lemmas, we see that by applying the super divergence
formula to $I^1_g(\phi)$, and pick out the sublinear
combination with $\sigma+1$ factors, with one factor $\nabla\phi$
and without internal contractions,\footnote{Notice that this sublinear combination must vanish separately.}
 we derive a {\it new} local equation:

\begin{equation}
\label{simonsjim} \sum_{l\in L^1_\sigma} a_l
C^{l,\iota}_g(\phi)+\sum_{r\in R} a_r
Xdiv_{j}C^{r,j}_g(\phi)=\sum_{h\in H} a_h
\nabla^j[C^h(g)]\nabla_j\phi
\end{equation}

\par Now, applying the operation $Weylify$ to the above, we
derive our Lemma \ref{palrights1} (by virtue of the discussion 
on the operation $Weylify$ in subsection 5.1 in \cite{alexakis1}). $\Box$

\par We will show Lemma \ref{definit} in the next section, after proving
Lemma \ref{killstacks}.

\subsection{Proof of Lemmas \ref{killstacks} and 
\ref{definit} (the second half of Lemma \ref{1stclaim}).}
 This subsection contains  elaborate constructions and calculations 
of local conformal invariants. All the divergences of vector fields 
that appear in the proofs of these two Lemmas are constructed explicitly; 
there is no recourse to the ``super divergnece formula''. The key 
to these constructions and to the calculations below is 
the ambient metric of Fefferman and Graham, see 
\cite{fg:ci}, \cite{fg:latest}. We will be using the algorithm 
that was presented in subsection \ref{algorithm}.
\newline

{\bf Proof of  Lemma \ref{killstacks}:} For this section, for
each complete contraction $C^l(g)$, $l\in L^{stack}$ we will call
the factor to which all the internal contractions belong the {\it
important factor}.

Firstly, we observe (easily) that by applying the fake second
Bianchi identities from \cite{alexakis1} (see (\ref{in&out}) in this paper)
  we can write modulo terms of length $\ge\sigma+1$:

$$\Sum_{l\in \tilde{L}^{stack}} a_l C^l(g)=
\Sum_{l\in \overline{L}^{stack}} a_l C^l(g)+\sum_{j\in J} a_j C^j_g(\phi),$$ 
where each $C^l(g)$
$l\in \overline{L}^{stack}$ is in the form (\ref{karami}) 
and has two of its internal contractions
involving the indices ${}_i,{}_k$ in the important factor
$\nabla^{(m)}W_{ijkl}$;\footnote{In other words, we have two derivative
indices $\nabla^i,\nabla^k$ contracting against the indices
${}_i,{}_k$ in the important factor} also each $C^j(g)$ is in the 
form (\ref{karami}) and has at least two internal
contractions belonging to {\it different} factors $\nabla^{(m)}W_{ijkl}$.

 Thus we are reduced
 to proving our claim under the assumption that
the important factor of each $C^l(g)$, $l\in L^{stack}$ is
 of the form $\nabla^{t_1\dots t_{\delta-2}ik}
\nabla^{(m)}_{r_1\dots r_m}W_{ijkl}$.

\par We then observe that by intergating by parts the 
indices ${}_{r_1},\dots, {}_{r_m}$ in the factor 
${\nabla_{r_1\dots r_m}}^{ik}W_{ijkl}$ for each $l\in L^{stack}$ we
can explicitly construct a vector field $T^i(g)=\Sum_{h\in H} a_h C^{h,i}(g)$
so that:

\begin{equation}
\label{wraiagynh} C^l(g)-div_iT^i(g)=\Sum C^{*}(g)+\Sum_{j\in J}
a_j C^j(g),
\end{equation}
modulo complete contractions of length $\sigma +1$. Here $\Sum
C^{*}(g)$ stands for a linear combination of complete contractions
in the form (\ref{karami}) with $\delta_W=2$, where both the
internal
 contractions belong to a factor in the form
$\nabla^{ik}W_{ijkl}$.  Each $C^j(g)$ is a complete
 contraction in the form (\ref{karami}) with at least two
internal contractions belonging to different factors. Now, abusing
notation we will again denote $\Sum_{l\in \overline{L}^{stack}}
a_l C^l(g)$ by $\Sum_{l\in L^{stack}} a_l C^l(g)$. Therefore, we
are reduced to showing our Lemma in the case where each $C^l(g)$,
 $l\in L^{stack}$ has $\delta_W=2$ and the two internal
 contractions belong to an important factor in the form
$\nabla^{ik}W_{ijkl}$.
\newline

\par In order to state our claim,
we recall the {\it ambient analogue} $Amb[C(g)]$,
of any complete contraction in the from (\ref{karami})
 (see Definition \ref{amban}); we also introduce a new definition:

\begin{definition}
\label{elab.constr}
Consider any $C^l(g)$, $l\in L^{stack}$ in
the form (\ref{karami}), where there are precisely two internal
contractions in $C^l(g)$, in some factor $\nabla^{ik}W_{ijkl}$.
We let $C^{l,i_1i_2}(g)$ be the tensor field  that arises from
$C^l(g)$ by making the two internal
 contractions into free indices. 

We also let $C^{l,i_1}(g)$ be the vector field that arises
from $C^l(g)$ by replacing the important factor
$\nabla^{ik}W_{ijkl}$ by $\nabla^kW_{i_1jkl}$ and $C^{l,i_2}(g)$
be the vector field that arises from $C^l(g)$ by replacing the
important factor $\nabla^{ik}W_{ijkl}$ by $\nabla^iW_{iji_2l}$.
\end{definition}

\begin{lemma}
\label{xtizw} 
Consider the tensor field $C^{l,i_1i_2}(g)$; consider 
$Xdiv_{i_1}Xdiv_{i_2}C^{l,i_1i_2}(g)$ and construct
$Amb[Xdiv_{i_1}Xdiv_{i_2}C^{l,i_1i_2}(g)]$. 
We claim that modulo terms of length $\ge\sigma+1$:

\begin{equation}
\label{loggos}
\begin{split}
&Amb[Xdiv_{i_1}Xdiv_{i_2}C^{l,i_1i_2}(g)]=
Xdiv_{i_1}Xdiv_{i_2}C^{l,i_1i_2}(g)+\frac{n-4}{n-3}Xdiv_{i_1}C^{l,i_1}(g)
\\&+\frac{n-4}{n-3}Xdiv_{i_2}C^{l,i_2}(g)+\frac{(n-4)^2-2(n-4)+2}{(n-3)(n-4)}C^l(g)+
 \Sum_{j\in J} a_j C^j(g)+
\\&\Sum_{t\in T^\sharp} a_t
Xdiv_{i_1}\dots Xdiv_{i_{a_t}}
 C^{t,i_1\dots i_{a_t}}(g)+\Sum_{t\in T^{\sharp\sharp}}
a_t Xdiv_{i_1}\dots Xdiv_{i_{a_t}} C^{t,i_1\dots i_{a_t}}(g).
\end{split}
\end{equation}
Here each tensor field $C^{t,i_1\dots i_{a_t}}(g)$, $t\in
T^\sharp$ is in the form (\ref{karami}), has $\delta_W>0$ and at
least one of the free indices belongs to a
 factor $T_1$, and at least one internal contraction belongs to a
 factor $T_2$ with $T_1\ne T_2$. Each tensor field
$C^{t,i_1\dots i_{a_t}}(g)$, $t\in T^{\sharp\sharp}$ also has
$\delta_W>0$ and in addition has a factor $T=\nabla^{(m)}W_{ijkl}$
with one internal contraction between a
 derivative index and an index ${}_i$ or ${}_k$ and then one of the
 free indices ${}_{i_1},\dots ,{}_{i_{a_t}}$ is the index ${}_j$ or ${}_l$ in $T$,
respectively. 
 $\Sum_{j\in J} a_j\dots$ is as in the
conclusion of Lemma \ref{killstacks}.
\end{lemma}

{\it Lemma \ref{xtizw} implies Lemma \ref{killstacks}:}
 Firstly, for each $t\in T^\sharp$ we suppose with no loss of
generality that ${}_{i_1}$ belongs to $T_1$ (see the definition
above).
 We then observe
that for each $t\in T^\sharp$:

\begin{equation}
\label{idf} \begin{split} &Xdiv_{i_1}\dots
Xdiv_{i_{a_t}}C^{t,i_1\dots i_{a_t}}(g)- div_{i_1}Xdiv_{i_2}\dots
Xdiv_{i_{a_t}}C^{t,i_1\dots i_{a_t}}(g)=
\\& \Sum_{j\in J} a_j C^j(g),
\end{split}
\end{equation}
modulo complete contractions of length $\ge\sigma +1$.

\par For each
$t\in T^{\sharp\sharp}$ we assume with no loss of generality that
${}_{i_1}$ is the index ${}_j$ or ${}_l$ (see the definition
above). Then, for each $t\in T^{\sharp\sharp}$:

\begin{equation}
\label{idf} Xdiv_{i_1}\dots Xdiv_{i_{a_t}}C^{t,i_1\dots
i_{a_t}}(g)- div_{i_1}Xdiv_{i_2}\dots Xdiv_{i_{a_t}}C^{t,i_1\dots
i_{a_t}}(g) =0,
\end{equation}
modulo complete contractions of length $\ge\sigma +1$.

\par By the same reasoning
 we explicitly construct  a linear
combination of vector fields $\Sum_{r\in R} a_r C^{r,i}(g)$ so
that:

\begin{equation}
\label{pousthdes}
\begin{split}
&Xdiv_{i_1}Xdiv_{i_2}C^{l,i_1i_2}(g)+\frac{n-4}{n-3}Xdiv_{i_1}
C^{l,i_1}(g)+\frac{n-4}{n-3}Xdiv_{i_2}C^{l,i_2}(g)+
\\&\frac{(n-4)^2-2(n-4)+2}{(n-3)(n-4)}C^l(g) -div_i \Sum_{r\in R}
a_r C^{r,i}(g)\\&= \frac{(n-4)(n-3)-2(n-4)^2+(n-4)^2-2(n-4)+2}{(n-3)(n-4)}
C^l(g).
\end{split}
\end{equation}
We note that the constant $C=\frac{-n+6}{(n-3)(n-4)}$ on the right
hand side is not zero, since for $n=6$ there can not be any complete contraction in the form:
\begin{equation}
\label{onlyW}
contr(\nabla^{(m_1)}W_{ijkl}\otimes\dots\otimes\nabla^{(m_\sigma)}W_{i'j'k'l'})
\end{equation}
 with a factor $\nabla^{ik}W_{ijkl}$.

\par Therefore, if we can prove Lemma \ref{xtizw}, our Lemma
 \ref{killstacks} will follow.
\newline

{\it Proof of Lemma \ref{xtizw}:}
\newline

Our proof relies on a careful calculation of $Amb[Xdiv_{i_1}Xdiv_{i_2}C^{l,i_1i_2}(g)]$,
based on the algorithm presented in subsection \ref{algorithm}.
We write $Xdiv_{i_1}Xdiv_{i_2}C^{l,i_1i_2}(g)$
 out as a sum in the obvious way,
 $Xdiv_{i_1}Xdiv_{i_2}C^{l,i_1i_2}(g)=
 \sum_{t=1}^{(\sigma-1)^2}C^t(g)$; then we write:

\begin{equation}
\label{definebreak2}
Amb[Xdiv_{i_1}Xdiv_{i_2}C^{l,i_1i_2}(g)]=\sum_{t=1}^{(\sigma-1)^2}\sum_{ass\in
ASSIGN^t}\sum_{\gamma\in BR_{ass}} \gamma\{ass\{
Amb[C^t(g)]\}\}.
\end{equation}

\par Now, our Lemma will follow by a careful analysis of the right
 hand side of the above. To perform this analysis we must divide the right hand side into
further sublinear combinations. 

\par We call the factor
$W_{i_1ji_2l}$ in $C^{l,i_1i_2}(g)$ to which the two free indices
belong the {\it important factor}. For the purposes of the discussion below,
each complete contraction $C^t(g)$ above, the indices
 ${}_{i_1}$, ${}_{i_2}$ in the important factor
 (which are contracting against $\nabla^{i_1},\nabla^{i_2}$)
  will still be called the {\it free indices}.

\begin{definition}
\label{karamanal}
Refer to (\ref{definebreak2}). Consider any given term $Amb\{C^t(g)\}$. We will write $ASSIGN$ 
instead of $ASSIGN^t$ for simplicity.

\par Define $ASSIGN^{-}\subset ASSIGN$ to stand for the index set of assignments that assign
the values $({}_\infty,{}_0)$ to two particular contractions
$({}_a,{}_b), ({}_c,{}_d)$ where ${}_a,{}_c$ belong to
different factors.

Define $ASSIGN^1\subset (ASSIGN\setminus ASSIGN^{-})$ to stand for
the index set of assignments in $(ASSIGN\setminus ASSIGN^{-})$
that assign to at least one particular contraction $({}_a,{}_b)$, where
${}_a$ belongs to the important factor, the values
$({}_0,{}_\infty)$.

Define $ASSIGN^\sharp \subset [ASSIGN\setminus (ASSIGN^{-}\bigcup
ASSIGN^1)]$ to stand for the index set of assignments in
$(ASSIGN\setminus ASSIGN^{-})$ that 
 either do not assign the values
 $({}_\infty,{}_0)$ to any particular contraction, or
 do assign such values
$({}_\infty,{}_0)$ but subject to the restriction that 
 all the ${}_\infty$'s
are assigned to the same factor, which is not the important factor.

 Define $ASSIGN^{*} \subset [ASSIGN\setminus (ASSIGN^{-}\bigcup ASSIGN^1)]$
to stand for the index set of assignments in $(ASSIGN\setminus
ASSIGN^{-})$ for which the value ${}_\infty$ is assigned 
to at least one index,  and moreover 
all the indices ${}_\infty$ are assigned to the important factor,
 but not to the free indices ${}_{i_1},{}_{i_2}$.

 Define $ASSIGN^{+} \subset [ASSIGN\setminus (ASSIGN^{-}\bigcup ASSIGN^1)]$
to stand for the index set of assignments in $(ASSIGN\setminus
ASSIGN^{-})$ for which all the values ${}_\infty$ are assigned to indices in 
the important factor,
 and at least one  of the free indices ${}_{i_1},{}_{i_2}$ is assigned the value ${}_\infty$.
\end{definition}

\par Let us firstly observe that for each $t, 1\le t\le (\sigma-1)^2$:

\begin{equation}
\label{eokab}
ASSIGN=ASSIGN^{-}\bigcup ASSIGN^1\bigcup ASSIGN^\sharp\bigcup ASSIGN^{*}\bigcup ASSIGN^{+};
\end{equation}
(in the above $\bigcup$ is a disjoint union).
We define:
$$\sum_{ASSIGN^{-}}ass\{ Amb[Xdiv_{i_1}Xdiv_{i_2}C^{l,i_1i_2}(g)]\}
=\sum_{t=1}^{(\sigma-1)^2}\sum_{ASSIGN^{-}}ass\{Amb[C^t(g)]\},$$
and also use the same definition for the other subsets of $ASSIGN$.

\par We will use (\ref{eokab}), along with (\ref{definebreak2}). Before doing so,
 let us recall a few more facts about the ambient metric:

\par A key observation is that any tensor in the form
$\tilde{\nabla}^{(q)}_{f_1\dots f_d}
\tilde{R}_{zxcv}$ as in (\ref{kefalo}) (i.e. with at least one index have the value ${}_\infty$)
can be written out as a linear combination of tensors
$\nabla^{(y)}W_{ijkl}(\otimes
 g\dots \otimes g)$ {\it with} an internal contraction in
$\nabla^{(y)}W$, modulo quadratic terms in the curvature.\footnote{Here
$(\otimes g\dots \otimes g)$ means that there may be no such
factors.}

\par We only need to be more precise in the case of the factor
$T=\tilde{R}_{ijkl}$: If one of the indices ${}_i,{}_j,{}_k,{}_l$
is given the value ${}_0$ then $T=0$. On the other hand, if $1\le {}_j,{}_k,{}_l\le n$
 then $\tilde{R}_{\infty jkl}=-\frac{1}{n-3}\nabla^i
W_{ijkl}$. If at least one of the pairs ${}_i,{}_j$ or ${}_k,{}_l$
are both assigned the value ${}_\infty$ then $T=0$. Finally, if $1\le j,l\le n$ then
$\tilde{R}_{\infty j\infty l}=\frac{1}{(n-3)(n-4)}\nabla^{ik}W_{ijkl}+Q(R)$ 
(see \cite{fg:ci, fg:latest}).
\newline

\par We will now calculate 
$\sum_{ass\in ASSIGN^{-}}
ass[Xdiv_{i_1}Xdiv_{i_2}C^{l,i_1i_2}(g)],$ $\dots,$ \\$\sum_{ass\in ASSIGN^{+}}
ass[Xdiv_{i_1}Xdiv_{i_2}C^{l,i_1i_2}(g)]$. All equations 
below hold  modulo terms of length $\ge\sigma+1$.

 It immediately follows that:

\begin{equation}
\label{agapil} \sum_{ass\in ASSIGN^1}
ass[Xdiv_{i_1}Xdiv_{i_2}C^{l,i_1i_2}(g)]=0,
\end{equation}

\begin{equation}
\label{russland} \sum_{ass\in ASSIGN^{-}}
ass[Xdiv_{i_1}Xdiv_{i_2}C^{l,i_1i_2}(g)]= \Sum_{j\in J} a_j C^j(g).
\end{equation}

\par We next seek to understand $\sum_{ass\in ASSIGN^\sharp}
 ass[Xdiv_{i_1}Xdiv_{i_2}C^{l,i_1i_2}(g)]$. We claim that:

 \begin{equation}
\label{sabbatob}
\begin{split}
&\sum_{ass\in ASSIGN^\sharp}
ass[Xdiv_{i_1}Xdiv_{i_2}C^{l,i_1i_2}(g)]=
Xdiv_{i_1}Xdiv_{i_2}C^{l,i_1i_2}(g)+
\\& \sum_{t\in T^\sharp} a_t Xdiv_{i_1}\dots Xdiv_{a_t}
 C^{t,i_1\dots i_{a_t}}(g)+\sum_{j\in J} a_j C^j(g).
\end{split}
 \end{equation}

{\it Proof of (\ref{sabbatob}):} Consider the {\it break-ups} of each different 
$ass[Xdiv_{i_1}Xdiv_{i_2}C^{l,i_1i_2}(g)]$:\footnote{Recall the {\it break-ups} from subsection \ref{algorithm}.}

\begin{equation}
\label{sienna} 
\begin{split}
&\sum_{ass\in ASSIGN^\sharp}
ass[Xdiv_{i_1}Xdiv_{i_2}C^{l,i_1i_2}(g)]=
\\&\sum_{ass\in ASSIGN^\sharp}\sum_{\gamma\in BR_{ass}}
\gamma\{ ass[Xdiv_{i_1}Xdiv_{i_2}C^{l,i_1i_2}(g)]\}
\end{split}
\end{equation}
Denote the RHS of the
above by $\Lambda$.
 We will break $\Lambda$ into two sublinear
combinations, $\Lambda^1,\Lambda^2$: A term 
$\gamma\{ ass[Xdiv_{i_1}Xdiv_{i_2}C^{l,i_1i_2}(g)]\}$
 will belong to $\Lambda^1$ if it arises by
applying one of the rules $\tilde{\nabla}_{i_s}X_j\rightarrow
\tilde{\Gamma}_{i_sj}^\infty X_\infty=-g_{i_sj}X_\infty$ or
$\tilde{\nabla}_{i_s}X_0\rightarrow
\tilde{\Gamma}_{i_s0}^kX_k=\delta^k_{i_s}X_k$ at least once (as in (\ref{breakTwell})).
 (Here $\tilde{\nabla}_{i_s}$ stands for a divergence index,\footnote{This means that the index 
is not one of the indices in the vector field $Amb[C^{l,i_1i_2}(g)]$ but corresponds to 
an index from a divergence $div_{i_1}$ or $div_{i_2}$.} which has
been assigned a value between ${}_1,\dots,{}_n$, and ${}_j$ is an original index which
has been assigned a value ${}_0,\dots ,{}_n$). 
A term $\gamma\{ ass[Xdiv_{i_1}Xdiv_{i_2}C^{l,i_1i_2}(g)]\}$
 belongs to $\Lambda^2$ if it 
arises without applying any of the above two Christoffel symbols.

\par In the above notation, we will show that:

\begin{equation}
\label{imamis1} \Lambda^1=0,
\end{equation}

 \begin{equation}
\label{imamis2} \Lambda^2= Xdiv_{i_1}Xdiv_{i_2}C^{l,i_1i_2}(g)+
\sum_{t\in T^\sharp} a_t Xdiv_{i_1}\dots Xdiv_{a_t} 
C^{t,i_1\dots i_{a_t}} +\sum_{j\in J} a_j C^j(g).
\end{equation}

\par In fact, (\ref{imamis2}) just follows by the definition of
$\Lambda^2$: If we consider the complete contractions that belong
to $\Lambda^2$ that arise without 
assigning any index the value ${}_\infty$, it follows that they add up to
the term $Xdiv_{i_1}Xdiv_{i_2}C^{l,i_1i_2}(g)$; if we consider the
complete contractions in $\Lambda^2$ that arise by assigning at least one index the value ${}_\infty$,
 it follows that they will add up to a linear combination
$\sum_{t\in T^\sharp} a_t Xdiv_{i_1}\dots Xdiv_{a_t} 
C^{t,i_1\dots i_{a_t}}(g)+\sum_{j\in J} a_j C^j(g)$.

\par So, matters are reduced to showing (\ref{imamis1}).
We will establish a bijection that will help us prove cancellation
in $\Lambda^1$.

\par We arbitrarily pick out one of the free indices ${}_{i_s}$ (so $s=1$ or $s=2$).
 For convenience we just set $s=1$, but the same claim will be true if we just
switch ${}_{i_1}$ and ${}_{i_2}$.
We also arbitrarily pick out a particular contraction in the
 tensor field $C^{l,i_1i_2}(g)$, say $\pi=({}_a,{}_b)$
 where ${}_a$ belongs to the factor $T^k$ and ${}_b$ belongs
 to the factor $T^l$ and $k,l\ne 1$.
 We also pick out a factor $T^r, r>1$ arbitrarily.
Then, we consider the complete contraction $C^{l,i_1i_2|k,r}(g)$
that arises  in $Xdiv_{i_1}Xdiv_{i_2} C^{l,i_1i_2}(g)$  when
$\nabla^{i_1}$ hits the $k^{th}$ factor
 $T^k$ in $C^{l,i_1i_2}(g)$ and $\nabla^{i_2}$ hits the $r^{th}$ factor $T^r$ in
$C^{l,i_1i_2}(g)$. Accordingly, we consider the complete
contraction $C^{l,i_1i_2|l,r}(g)$ that arises  in
$Xdiv_{i_1}Xdiv_{i_2} C^{l,i_1i_2}(g)$  when $\nabla^{i_1}$ hits
the $l^{th}$
 factor $T^l$ in $C^{l,i_1i_2}(g)$ and $\nabla^{i_2}$
  hits the $r^{th}$ factor $T^r$ in $C^{l,i_1i_2}(g)$.

\par We now consider the ambient analogues
$Amb[C^{l,i_1i_2|k,r}(g)]$, $Amb[C^{l,i_1i_2|k,r}
(g)]$. For $Amb[C^{l,i_1i_2|l,r}(g)]$ we
define $ASSIGN^{\sharp,\pi,A}_{i_1}\subset ASSIGN^\sharp$
 to stand for the set of assignments that assign the particular
 contraction $\pi=({}_a,{}_b)$ the values $({}_0,{}_\infty)$ and
assign the pair $({}_{i_1},\nabla^{i_1})$ base values (i.e.
${}_{i_1}\in \{{}_1,\dots ,{}_n\})$.

\par On the other hand, for $Amb[C^{l,i_1i_2|l,r}
(g)]$ we define $ASSIGN^{\sharp,\pi,B}_{i_1}$ to
 stand for the set of assignments where $\pi=({}_a,{}_b)$ is assigned base values
 and the pair $({}_{i_1},\nabla^{i_1})$ is also
 assigned base values.

\par Now, for each element $ass\in ASSIGN^{\sharp,\pi,A}_{i_1}$
we define $BR_{ass}^1$ to stand for the set of break-ups
 that replace the factor $T^k=\tilde{\nabla}^{(m+1)}_{i_sr_1\dots 0\dots r_m}
 \tilde{R}_{ijkl}$ (or analogously
when ${}_a={}_0$ is one of the internal indices) by a factor
$-\tilde{\nabla}^{(m)}_{r_1\dots i_1\dots r_m}\tilde{R}_{ijkl}$
(i.e. we apply $\tilde{\Gamma}^k_{i_s0}=\delta^k_{i_s}$ to the
indices $\tilde{\nabla}_{i_s},\tilde{\nabla}_0$, as in (\ref{breakTwell})). For
$C^{l,i_1i_2|l,r}(\tilde{g}^{n+2})$ and
 each element $ass\in ASSIGN^{\sharp,\pi,B}_{i_1}$,
we define $BR_{ass}^2$ to stand for the set of break-ups
 that replace the factor $T^l=
\tilde{\nabla}^{(m'+1)}_{i_sr'_1\dots b\dots r_{m'}}
\tilde{R}_{i'j'k'l'}$ (or analogously where $b$ is one of the
 internal indices) by a factor
$\tilde{\nabla}^{(m')}_{r'_1\dots \infty\dots r'_m}
\tilde{R}_{i'j'k'l'}g_{i_1b}$ (i.e. we apply
$\tilde{\Gamma}^\infty_{i_sb}=-g_{i_sb}$.

\par We then observe that in the above notation,
for each $\pi\in \Pi$, each $r=2,\dots ,\sigma$:

\begin{equation}
\label{broeder}
\begin{split}
& \Sum_{ass\in ASSIGN^{\sharp,\pi,A}_{i_s}}
\Sum_{\gamma\in BR^1_{ass}} \gamma\{ass
[Amb[C^{l,i_1i_2|k,r}(g)]]\}
\\&+ \Sum_{ass\in
ASSIGN^{\sharp,\pi,B}_{i_s}} \Sum_{\gamma\in BR^2_{ass}}
\gamma\{[Amb[C^{l,i_1i_2|l,r}(g)]]\}=0.
\end{split}
\end{equation}

\par In view of the above cancellation, when analyzing
$\Lambda^1$ we may {\it discard} any contraction that involves
using the symbol $\tilde{\Gamma}^\infty_{ij}$ to a pair
$\tilde{\nabla}_{i_s}X_k$ (where ${}_{i_s}$ is a divergence
 index and $X_k$ is not contracting against the important factor). We may  also
discard any complete contraction that arises when we use the
symbol $\tilde{\Gamma}_{i_s0}^k$ to a pair
$\tilde{\nabla}_{i_s}X_d$ where $X_d$ is an original index in
$C^{l,i_1i_2}(g)$ that is not contracting the important factor.

\par So, to show (\ref{imamis1}), we only have to consider the
terms in the LHS that arise by using  one of the symbols
$\tilde{\Gamma}^\infty_{i_sj}$, $\tilde{\Gamma}_{i_s0}^k$ to a
pair $\tilde{\nabla}_{i_s}X_v$ where ${}_{i_s}$ is a divergence
 index and $X_v$ {\it is} contracting against the important factor
 (and where in addition ${}_v$ has been assigned a value ${}_0,\dots,{}_n$).

 \par Now, observe that there can be no contractions arising from
 an assignment $ass\in ASSIGN^\sharp$ with an index ${}_v$ having
 been assigned the value ${}_0$, and the same index ${}_v$ contracting
 against an index ${}_b$ in the important factor in
 $C^{l,i_1i_2}(g)$: In that case, the index ${}_b$ would have been
 assigned a value ${}_\infty$, which contradicts the fact that $ass\in
 ASSIGN^\sharp$.

\par So we only need to observe that we may discard any complete
contraction that arises when we use the symbol
$\tilde{\Gamma}^\infty_{i_sk}$ to a pair of indices
$(\tilde{\nabla}_{i_s}, {}_k)$ where $\tilde{\nabla}_{i_s}$ is the
divergence index and ${}_k$ is contracting against an index in the
important factor-this follows because
$\tilde{\Gamma}^\infty_{i_sk}=-g_{i_sk}$ and the Weyl tensor is
trace-free. This proves our claim. $\Box$
\newline

\par Next, we wish to analyze

$$\sum_{ass\in ASSIGN^{*}}\sum_{\gamma\in BR_{ass}} \gamma\{ass[Xdiv_{i_1}Xdiv_{i_2}C^{l,i_1i_2}(g)]\}.$$

Using the definition of $ASSIGN^{*}$ and the Christoffel symbols for $\tilde{g}$,
we observe that:

\begin{equation}
\label{kestek3}
\begin{split}
&\Sum_{ass\in ASSIGN^{*}} \sum_{\gamma\in BR_{ass}}\gamma\{ ass
Amb [Xdiv_{i_1}Xdiv_{i_2}C^{l,i_1i_2}(g)]\}=
\\&-\frac{1}{n-3}Xdiv_{i_1}C^{l,i_1}(g)-\frac{1}{n-3}Xdiv_{i_2}C^{l,i_2}(g)
\\&+\frac{2}{(n-3)(n-4)}C^l(g)+\Sum_{t\in T^\sharp} b_t Xdiv_{i_1}\dots
 Xdiv_{i_{a_t}}C^{t,i_1\dots i_{a_t}}(g)
 \\&+\Sum_{t\in T^{\sharp\sharp}} b_t Xdiv_{i_1}\dots
 Xdiv_{i_{a_t}}C^{t,i_1\dots i_{a_t}}(g)+\Sum_{j\in J} a_j C^j(g).
\end{split}
\end{equation}
Note: The first three terms in the RHS arise by applying the Christoffel symbol 
$\tilde{\Gamma}_{i_s0}^k=\delta_{i_s}^k$ at 
least once (we also use the first Bianchi 
identity here--the details are left to the reader).\footnote{Here 
${}_{i_s}$ comes from an index $\tilde{\nabla}_{i_s}$, from one 
of the divergences $div_{i_1},div_{i_2}$.} It is easy to observe that the sublinear combination 
that arises when the symbol 
$\tilde{\Gamma}_{i_s0}^k=\delta_{i_s}^k$ is never applied will equal  
$\Sum_{t\in T^{\sharp\sharp}}\dots$. 

Finally, we must understand  the sublinear combination
$$\sum_{ass\in ASSIGN^{+}} \sum_{\gamma\in BR_{ass}}\gamma\{ ass
[Xdiv_{i_1}Xdiv_{i_2}C^{l,i_1i_2}(\tilde{g}^{n+2})]\}.$$
 We calculate that:

\begin{equation}
\label{kestek4}
\begin{split}
&\Sum_{ass\in ASSIGN^{+}}ass
Amb[Xdiv_{i_1}Xdiv_{i_2}C^{l,i_1i_2}(g)]=
\frac{(n-4)^2}{(n-3)(n-4)}C^l(g)
\\&+Xdiv_{i_1}C^{l,i_1}(g)+Xdiv_{i_2}C^{l,i_2}(g) -\frac{2}{(n-3)} C^l(g)
\\&+ \Sum_{t\in T^{\sharp\sharp}} a_t Xdiv_{i_1}\dots
Xdiv_{i_{a_t}} C^{t,i_1\dots i_{a_t}}(g)+\Sum_{j\in J} a_j C^j(g).
\end{split}
\end{equation}
Note again, the first four terms in the RHS arise 
 by applying the Christoffel symbol 
$\tilde{\Gamma}_{i_s0}^k=\delta_{i_s}^k$ at 
least once (we also again use the first Bianchi 
identity).

\par Thus, in view of (\ref{definebreak2}), (\ref{eokab}), by just adding up the equations 
 (\ref{agapil}), (\ref{russland}),
(\ref{sabbatob}), (\ref{kestek3}), (\ref{kestek4}) we derive 
Lemma \ref{killstacks}. $\Box$
\newline

{\bf Proof of Lemma \ref{definit}:}
\newline

\par Now, to show Lemma \ref{definit}, we recall the notational 
conventions from the statement of the Lemma and introduce some additionnal 
ones:\footnote{The reader should 
note that the notational conventions here are different 
from the ones in the proof of the previous Lemma.}
Recall that $\Sum_{f\in F} a_f C^f(g)$ 
stands for a generic linear
combination of complete contractions in the form (\ref{karami})
with $\delta_W\ge 2$. Also, $\Sum_{t\in T^\sharp} a_t
C^{t,i_1\dots i_{a_t}}(g)$ will now stand for a generic
linear combination of tensor fields in the form (\ref{karami}) with
$a_t>0$ (i.e. with at least one free index) and with  one internal
contraction. All equations in this proof will hold modulo complete
contractions of length $\ge\sigma+1$.

\par Now, for each $C^v(g)$ as in the statement of Lemma
\ref{definit} we consider $\nabla_iC^v(g)$ (which is thought 
of as the sum of $\sigma$ partial contractions). 
Then, for each $v\in V$ we consider the linear
combination $Xdiv_i[\nabla_i C^v(g)]$ and we construct its {\it
ambient analogue} (see the Definition \ref{amban}):

$$Amb\{Xdiv_i[\nabla_i C^v(g)] \}.$$

 \par We then claim:

\begin{lemma}
\label{kallinikos}
For each $v\in V$:

\begin{equation}
\label{eqnkalli} \begin{split} &Amb\{Xdiv_i[\nabla_i C^v(g)]
\}=2\Delta_rC^v(g)+ Xdiv_i[\nabla_iC^v(g)]+\Sum_{f\in F} a_f
C^f(g) \\&+\Sum_{t\in T^\sharp} a_tXdiv_{i_1}\dots
Xdiv_{i_{a_t}}C^{t,i_1\dots i_{a_t}}(g),
\end{split}
\end{equation}
modulo complete contractions of length $\ge\sigma +1$.
\end{lemma}

\par We observe that if we can prove the above, then Lemma
\ref{definit} follows from two easy observations: Firstly:

$$Xdiv_i[\nabla_iC^v(g)]-div_i[\nabla_iC^v(g)]=-\Delta_rC^v(g).$$
Secondly, for every $t\in T^\sharp$:
$$Xdiv_{i_1}\dots
Xdiv_{i_{a_t}}C^{t,i_1\dots i_{a_t}}(g)- div_{i_1}Xdiv_{i_2}\dots
Xdiv_{i_{a_t}}C^{t,i_1\dots i_{a_t}}(g)=\sum_{f\in F} a_f C^f(g).$$
Thus, it suffices to show Lemma \ref{kallinikos}.
\newline

{\it Proof of Lemma \ref{kallinikos}:} Firstly, we will number the
factors in $C^v(g)$ and denote them by $F_1,\dots ,F_\sigma$.
Then, we will denote by $\nabla^\tau_iC^v(g)$ the vector field
that arises from $C^v(g)$ by replacing the factor $F_\tau$ by
$\nabla_i F_\tau$. We also denote by $m_\tau$ the number of
 derivatives in the factor $F_\tau=\nabla^{(m_\tau)}W_{ijkl}$.

\par We will then show that for any $\tau$,
$1\le \tau\le \sigma$:

\begin{equation}
\label{appel}
\begin{split}
&Amb\{ Xdiv_i[\nabla_i^\tau C^v(g)]\}=\Delta_r^\tau
C^v(g)+Xdiv_i[\nabla_i^\tau C^v(g)]+ \frac{m_\tau +2}{n-2}
\Delta_r C^v(g)\\& +\Sum_{f\in F} a_f C^f(g) +\Sum_{t\in T^\sharp}
a_tXdiv_{i_1}\dots Xdiv_{i_{a_t}}C^{t,i_1\dots i_{a_t}}(g).
\end{split}
\end{equation}

\par In view of the equation $\Sum_{\tau=1}^\sigma
(m_\tau+2)=n-2$, (\ref{appel}) implies 
 Lemma \ref{kallinikos}.
\newline

{\it Proof of (\ref{appel}):} We will show this equation using the
notions of {\it assignment} and {\it break-ups} from subsection \ref{algorithm}. 

We recall the following fact regarding the ambient
metric: Consider any tensor $T=\tilde{\nabla}^{(m)}_{\infty
r_2\dots r_m}\tilde{R}_{ijkl}$, where all the indices
${}_{r_2},\dots ,{}_{r_m},{}_i,{}_j,{}_k,{}_l$ have values between
${}_1$ and ${}_n$. It then follows from the formula $\partial_\infty \tilde{g}_{ab}=2P_{ab}$
in \cite{fg:latest} that:

\begin{equation}
\label{pritzeto} T=-\frac{1}{n-2}\Delta\nabla^{(m-1)}_{r_2\dots
r_m}W_{ijkl}+ \Sum_{h\in H} a_h [T^h(g)(\otimes g\dots \otimes
g)]_{r_2\dots r_mijkl}+Q(R),
\end{equation}
where each tensor $T^h(g)$ is in the
 form $\nabla^{(m')}W_{ijkl}$ with at least {\it two} internal
 contractions.

\par Furthermore, it can be seen from \cite{fg:ci}, \cite{fg:latest} that any component
$T'=\tilde{\nabla}^{(m)}_{r_1\dots r_m}\tilde{R}_{ijkl}$
 where at least {\it two} of the indices ${}_{r_1},\dots,{}_l$
having the value ${}_\infty$ can be expressed as:
\begin{equation} 
\label{pritzeto2}
T'=\Sum_{h\in H} a_h [T^h(g)(\otimes g\dots \otimes
g)]_{r_2\dots r_mijkl}+Q(R),
\end{equation}
 with the same conventions as above.

\par Now, for the purposes of the next definition
we write out $Amb\{Xdiv_i[\nabla_i^\tau
C^v(g)]\}=\sum_{t=1}^{\sigma-1} Amb[C^t(g)]$. We prove our Lemma
by examining the right hand side of the equation:

\begin{equation}
\label{magdalen} Amb\{ Xdiv_i[\nabla_i^\tau
C^v(g)]\}=\sum_{t=1}^{\sigma-1}\Sum_{ass\in ASSIGN^t} \sum_{\gamma\in BR_{ass}} \gamma\{ass[
Amb[C^t(g)]]\}.
\end{equation}

\par It will again be useful to break the right hand side of the
above into sublinear combinations. We will now call the factor
$F_\tau$ (to which $\nabla_i$ belongs) the {\it
important factor}. We will
also call the index $\nabla_i$ in
that factor the ``free index''.

\begin{definition}
\label{kolligiannis} Refer to (\ref{magdalen}). Consider any fixed
$Amb[C^t(g)]$. We write $ASSIGN$ instead of $ASSIGN^t$ for simplicity. 

 We define $ASSIGN^{-}\subset ASSIGN$ to stand for
the index set of assignments where at least two different indices in $Amb[C^t(g)]$ 
are assigned the value ${}_\infty$.

\par We define $ASSIGN^\sharp \subset (ASSIGN\setminus ASSIGN^{-})$ to stand for the index set
of assignments where the pair of values $({}_\infty,{}_0)$ is
assigned to at most one  particular contraction, and moreover if such a
pair is assigned to a particular contraction $({}_a,{}_b)$ then
neither ${}_a$ nor ${}_b$ belongs to the important factor.

\par We define $ASSIGN^1\subset (ASSIGN\setminus ASSIGN^{-})$ to stand for
 the set of assignments in $(ASSIGN\setminus ASSIGN^{-})$ that assign 
the pair of values $({}_0,{}_\infty)$ to exactly one particular contraction $({}_a,{}_b)$ 
 where in addition ${}_a$ belongs to the
important factor.

\par We define $ASSIGN^{*}\subset (ASSIGN\setminus ASSIGN^{-})$ to stand for
 the set of assignments in $(ASSIGN\setminus ASSIGN^{-})$ that 
assign the values $({}_\infty,{}_0)$ 
to exactly one particular contraction $({}_a,{}_b)$ 
 where in addition ${}_a$ belongs to the
important factor and is not the free index $\nabla_i$.

\par We define $ASSIGN^{+}\subset (ASSIGN\setminus ASSIGN^{-})$ to stand for
 the assignment in $(ASSIGN\setminus ASSIGN^{-})$ that
 assigns the value ${}_\infty$ to the free index $\nabla_i$, 
(and the value ${}_0$ to the index against which it contracts) and does not assign 
this value to any other index. 
\end{definition}

\par Clearly:

\begin{equation}
\label{kviri} ASSIGN=ASSIGN^{-}\bigcup ASSIGN^\sharp\bigcup
ASSIGN^{1}\bigcup ASSIGN^{*}\bigcup ASSIGN^{+};
\end{equation}
(where in the above $\bigcup$ stands for a disjoint union).

\par Firstly, by virtue of (\ref{pritzeto}), (\ref{pritzeto2}) we observe that:

\begin{equation}
\label{darouinos} \sum_{ass\in ASSIGN^{-}}
ass[Amb[Xdiv_i[\nabla_i^\tau C^v(g)]]]=\Sum_{f\in F} a_f C^f(g).
\end{equation}

By repeating the cancellation argument from equation 
(\ref{imamis1}) we derive that:

\begin{equation}
\label{glasnost} \begin{split} & \sum_{ass\in ASSIGN^\sharp}
 ass\{Amb[Xdiv_i[\nabla_i^\tau
C^v(g)]]\}=Xdiv_i[\nabla_i^\tau C^v(g)]
+\frac{m_\tau+2}{n-2}\Delta_r^\tau C^v(g)
\\&+\sum_{t\in T^\sharp}
a_t Xdiv_{i_1}\dots Xdiv_{b_t} C^{t,i_1\dots
i_{b_t}}(g)+\sum_{f\in F} a_f C^f(g).
\end{split}
\end{equation}
Note: The expression $\frac{m_\tau+2}{n-2}\Delta_r^\tau C^v(g)$
arises from the assignments where no index is assigned a value
${}_\infty$, by applying $\tilde{\Gamma}_{ia}^\infty=-g_{ia}$ to
the pairs of indices ${}_i,{}_a$ (${}_i$ is the free index) in the
important factor.

\par Now, by virtoe of the Christoffel symbol $\tilde{\Gamma}_{0i}^k=\delta_i^k$,
 we straightforwardly derive:

\begin{equation}
\label{agapil3}
\begin{split}
&\Sum_{ass\in ASSIGN^1} ass[Amb[Xdiv_i[\nabla_i^\tau C^v(g)]]]=
\\&\frac{m_\tau+2}{n-2}\Sum_{k=1,k\ne\tau}^\sigma\Delta_r^k
[C^v(g)]+\Sum_{t\in T^\sharp} a_t Xdiv_{i_1} \dots
Xdiv_{i_{a_t}}C^{t,i_1\dots i_{a_t}}(g)+ \Sum_{f\in F} a_f C^f(g).
\end{split}
\end{equation}

\par Furthermore, we derive:

\begin{equation}
\label{kestek3'}
\begin{split}
& \Sum_{ass\in ASSIGN^{*}} ass\{ Amb [Xdiv_i[\nabla_i^\tau C^v(g)]]\}=
\\&\frac{m_\tau+2}{n-2}\Delta_r^\tau C^v(g)+
\Sum_{t\in T^\sharp} a_t Xdiv_{i_1}\dots Xdiv_{i_{a_t}}
C^{t,i_1\dots i_{a_t}}(g)+\Sum_{f\in F} a_f C^f(g).
\end{split}
\end{equation}
{\it Note:} The sublinear combination $\frac{m_\tau+2}{n-2}\Delta_r^\tau
C^v(g)$ arises as follows: Recall that $ass\in ASSIGN^*$ assigns one pair
of values $({}_\infty,{}_0)$ to any particular contraction
$({}_a,{}_b)$ where ${}_a$ belongs to the important factor and
${}_a$ is {\it not} the free index ${}_i$.  The coefficient then
arises when we use the Christoffel symbol
$\tilde{\Gamma}_{i0}^k=\delta^k_i$ to the pair
$\tilde{\nabla}_iX_b$ (where $X_b$ stands for the index that is
contracting against the index ${}_a$ in the important factor-hence
${}_b$ has been assigned the value ${}_0$). The second Bianchi identity 
is also employed here.

Finally we calculate:

\begin{equation}
\label{kestek3''}
\begin{split}
&\Sum_{ass\in ASSIGN^{+}} ass [Amb Xdiv_i[\nabla_i^\tau C^v(g)]=\frac{n-2-(m_\tau+2)}{n-2}\Delta_r^\tau C^v(g)+
\\&\Sum_{t\in T^\sharp} a_t Xdiv_{i_1}\dots Xdiv_{i_{a_t}}
C^{t,i_1\dots i_{a_t}}(g)+\Sum_{f\in F} a_f C^f(g).
\end{split}
\end{equation}
$\frac{n-2-(m_\tau+2)}{n-2}\Delta_r^\tau C^v(g)$ arises by
 assigning the pair $({}_\infty,{}_0)$ to the divergence pair
$({}_i,\nabla^i)$. The coefficient $\frac{n-2-(m_\tau+2)}{n-2}$ arises by virtue of
 the formula $\tilde{\Gamma}_{i0}^k$ applied to any pair
 $(\tilde{\nabla}_i,X_k)$.(Recall that there the diregence index $\tilde{\nabla}_i$ is
assigned the value ${}_0$).

\par Plugging the equations (\ref{darouinos})--(\ref{kestek3''})
into (\ref{magdalen})
we obtain our Lemma \ref{definit}. $\Box$

\section{Proof of Lemma \ref{nabothewri2}: The divergence ``piece'' in $P(g)$.}
\label{puredivergencesec}

The aim of this section is to prove that when $P(g)$ is 
as in the hypothesis of Lemma \ref{nabothewri2},
then the sublinear combination $\sum_{l\in L^j} a_l C^l(g)$ in $P(g)$\footnote{Recall that 
$\sum_{l\in L^j} a_l C^l(g)$ stands for the sublinear combination of  terms with
 length $\sigma$ and $j$ internal contractions; recall 
that $j\ge 2$ and all other terms of length $\sigma$ in $P(g)$ are 
assumed to be in the form (\ref{karami}) with at least $j+1$ internal contractions.}
can be cancelled out (modulo 
introducing new ``better'' terms as in the statement of Lemma \ref{nabothewri2}) 
by subtracting a divergence of a vector field, 
as allowed by the Deser-Schwimmer conjecture. 

Again, we rely on the only tool we have at our disposal, which 
allows us to pass from the invariance under integration enjoyed by $P(g)$ 
to a {\it local equation}: 
We consider the linear operator $I^1_g(\phi):=[\frac{d}{dt}|_{t=0}[e^{nt\phi}P(e^{2t\phi}g)]$, 
for which $\int_{M^n}I^1_g(\phi)dV_g=0$. 
The aim is to invoke the super divergence formula for $I^1_g(\phi)$,\footnote{Which 
expresses $I^1_g(\phi)$ as a 
divergence, $I^1_g(\phi)=div_iT^i_g(\phi)$} to derive 
the claim in the previous paragraph. The are two challenges one must address: 
Firstly, how can one recover $\sum_{l\in L^j} a_l C^l(g)$ by examining $I^1_g(\phi)$? 
Secondly, how can one use the super divergence formula 
for $I^1_g(\phi)$ to derive the claim on $\sum_{l\in L^j} a_l C^l(g)$?

The first challenge is not hard; it follows by studying the transformation 
law of covariant derivatives of the Weyl curvature, paying 
special attention to {\it internal contractions}.\footnote{See the
 definition above Proposition \ref{bigprop}.} 
This is carried out in subsection \ref{prelimwork}. 
The second challenge is not straightforward; one difficulty, already discussed in \cite{alexakis1},  is 
again that a direct application of the super divergence formula to $I^1_g(\phi)$ 
{\it does not} imply Lemma \ref{nabothewri2}. For that reason, we must formulate 
certain ``main algebraic propositions'', see Propositions 
\ref{pregiade2} and \ref{tool'} in subsection \ref{imp.tools} below, which 
will allow us to derive Lemma \ref{nabothewri2} from the super 
divergence formula applied to $I^1_g(\phi)$. 
These propositions can be viewed as 
analogues of the ``main algebraic Proposition'' in subsection 5.3 \cite{alexakis1}; all three of these 
algebraic propositions will be proven in the series 
of papers \cite{alexakis4, alexakis5, alexakis6}.

However, there is an additional important challenge in this case, briefly 
discussed at the end of subsection \ref{outlofarg}: In a nutshell, the problem 
concerns a certain ``loss of information'' which occurs when we pass from $I^1_g(\phi)$ 
(which satisfies the integral equation  
 $\int_{M^n}I^1_g(\phi)dV_g=0$) to the super divergence 
 formula applied  to $I^1_g(\phi)$, $supdiv[I^1_g(\phi)]=0$. 
As explained at the end of subsection \ref{outlofarg}, 
 upon examining the terms with $\sigma+1$ factors in $I^1_g(\phi)$
 which also have a factor $\nabla\phi$,\footnote{In \ref{outlofarg} we denoted the
  sublinear combination of those terms by $(I^1_g(\phi))_{\nabla\phi}$--see (\ref{bo'az13}).} 
  we can {\it recover} the sublinear combination 
$\sum_{l\in L^j} a_l C^l(g)$ in $P(g)$. However, when we examine the 
terms of length $\sigma+1$  in $supdiv[I^1_g(\phi)]=0$ with a factor $\nabla\phi$, we find that  
those terms arise {\it both} from the sublinear combination $(I^1_g(\phi))_{\nabla\phi}$ but {\it also}
from the sublinear combination of terms $\sum_{q\in Q} a_q C^q(g)\cdot\Delta\phi$ in (\ref{bo'az13}) (which contain 
a factor $\Delta\phi$; this factor gives rise to a factor $\nabla\phi$ 
in the super divergence formula, due to integrations by parts). Therefore
 the additional challenge in this case is that we must somehow ``get rid'' of the
 potentialy  harmful terms in $\sum_{q\in Q} a_q C^q(g)\cdot \Delta(\phi)$; 
 we note that apriori we have {\it no} information on this sublinear combination. 
 We achieve this goal by virtue of the fact that the terms in 
 $\sum_{l\in L^j} a_l C^l(g)$ in the hypothesis of Lemma \ref{nabothewri2}) have at least two
  internal contractions belonging to {\it different factors}.\footnote{Recall
 that we were able to reduce ourselves to 
  this setting by virtue of our elaborate constructions in the ambient metric in the previous section.} 
  Using this fact, we are able to use the super divergence formula
  $supdiv[I^1_g(\phi)]=0$  (and the ``main algebraic propositions'')  in 
  two different ways, in order to first ``get rid'' of the potentially harmful terms in 
  $\sum_{q\in Q} a_q C^q(g)\cdot \Delta(\phi)$ and {\it then} 
 to derive Lemma \ref{nabothewri2}.

\subsection{Preliminary work: How can one recover 
$P(g)$ from $I^1_g(\phi)$?}
\label{prelimwork}

\par As preliminary work for this section, we will
study the images $Image^1_{\phi}[C(g)]$
 of complete contractions $C(g)$
in the form (\ref{karami}) with weight $-n$ and length $\sigma$.\footnote{Recall 
that $Image^1_\phi[C(g)]=\frac{d}{dt}|_{t=0}[e^{nt\phi}C^{2t\phi}g)]$.}

\par In order to do this, we will study the transformation laws of tensors
$T=\nabla^{r_{a_1}\dots r_{a_s}}\nabla^{(m)}_{r_1\dots r_m}
W_{r_{m+1}r_{m+2}r_{m+3}r_{m+4}}$, where each of the indices
 ${}^{r_{a_x}}$ is contracting against an index
${}_{r_{a_x}}$ and all the other indices in $T_h(g)$ are free. For each
such tensor $T_h(g)$, we define $Image^1_{\phi}[T_h]$ to stand for
the sublinear combination in

$$\frac{d}{d\lambda}{|}_{\lambda=0}T_h(e^{2\lambda\phi}g)$$
that involves factors $\nabla^{(p)}\phi$ with $p>0$ (i.e. we are
excluding the partial contractions with a factor $\phi$, without
derivatives).

\par Clearly, for each complete
contraction $C(g)$ in the form (\ref{karami}):

\begin{equation}
\label{anatrans} Image^1_\phi[C(g)]=\Sum_{h=1}^\sigma
C^{sub(T_h)}(g),
\end{equation}
where $C^{sub(T_h)}(g)$ stands for the complete contraction
 that arises from $C(g)$ by replacing the factor $T_h$ by
$Image^1_{\phi}[T_h]$ and then contracting indices according to the same pattern as for $C(g)$. 

{\it Transformation laws:} Recall that the Weyl tensor is conformally invariant, 
i.e.~$W_{ijkl}(e^{2\phi}g)e^{2\phi}W_{ijkl}(g)$. 
Recall also the transformation law of the Levi-Civita connection:
\begin{equation}
\label{levicivita} {\nabla}_k {\eta}_l(e^{2\phi}g)=
{\nabla}_k{\eta}_l(g) -\nabla_k{\phi} {\eta}_l -\nabla_l{\phi} {\eta}_k
+\nabla^s{\phi} {\eta}_s g_{kl}.
\end{equation}

\par Now, in order to perform our calculations, we will have to introduce
some notational conventions. For the purposes of this subsection,
when we write $\nabla^{(m)}R_{ijkl}, \nabla^{(m)}Ric$ or
$\nabla^{(m)}R$ ($R$ is the scalar curvature) we will mean a
tensor in those forms, possibly with some internal contractions.
 Also as usual,
 whenever we write $\nabla^{(m)}R_{ijkl}$ we will mean
 that no two of the indices ${}_i,{}_j,{}_k,{}_l$ are contracting between
themselves and, when we write $\nabla^{(p)}Ric_{ij}$ we
 will mean that no two of the indices ${}_i,{}_j$ are contracting
 between themselves.

{\it Notation:}
We will denote by $T^\alpha(\phi)$ a partial contraction of the
form
 \begin{equation}
 \label{william1}
\nabla^{(m)}_{r_1\dots r_m}R_{r_{m+1}r_{m+2}r_{m+3}r_{m+4}}\nabla^{r_a}\phi,
 \end{equation}
where the factor $\nabla\phi$ is contracting against one of the
 indices ${}_{r_1},\dots, {}_{r_{m+4}}$ in the first factor. 
 $T^\beta(\phi)$ will stand for a partial contraction in
the form:

 \begin{equation}
 \label{william3}
\nabla^{(m)}_{r_1\dots r_m}R_{ijkl}\Delta\phi.
 \end{equation}
$T^\gamma(\phi)$ will stand for a partial contraction in the form:

 \begin{equation}
 \label{william4}
\nabla^{(m)}_{r_1\dots r_m}R_{ijkl}\nabla_x\phi,
 \end{equation}
where ${}_x$ is a free index. $T^\delta(\phi)$ will stand for a
 partial contraction in the form:

 \begin{equation}
 \label{william5}
\nabla^{(m)}_{r_1\dots r_m}R_{ijkl}\nabla^{(2)}_{yx}\phi,
 \end{equation}
where both ${}_y$ and ${}_x$ are free indices; finally,
$T^\epsilon(\phi)$ will be a generic tensor product in the form:

 \begin{equation}
 \label{william6}
\nabla^{(m)}_{r_1\dots r_m}R_{ijkl}\nabla^{(p)}_{s_1\dots s_p}\phi
 \end{equation}
where either $p\ge 3$ or $p=2$ and at least one of the indices
${}_{s_1},{}_{s_2}$ is contracting against the factor $\nabla^{(m)}R_{ijkl}$.

\par We also generically denote by
$T_{Ric}^\alpha(\phi), T_{Ric}^\beta(\phi), T_{Ric}^\gamma(\phi)$,
$T_{Ric}^\delta(\phi)$, $T_{Ric}^\epsilon(\phi)$,
$T_{Ric}^\zeta(\phi)$ tensor fields that are as above, but only with the
 factor $\nabla^{(m)}R_{ijkl}$ formally replaced by an expression
$\nabla^{(m)}[Ric\otimes g]_{ijkl}$ (the indices
${}_i,{}_j,{}_k,{}_l$ belonging to the factors $Ric,g$). We also
denote by $T_{R}^\alpha(\phi), T_{R}^\beta(\phi),
T_{R}^\gamma(\phi)$, $T_{R}^\delta(\phi)$, $T_{R}^\epsilon(\phi)$,
$T_{R}^\zeta(\phi)$ tensor fields that are as above, but
only with the factor $\nabla^{(m)}R_{ijkl}$ replaced by an expression
$\nabla^{(m)}R\otimes [g\otimes g]_{ijkl}$.

\par Now, for each factor in one of the 15 forms above, we
 will denote by $\delta$ the total number of internal
 contractions in the curvature term (i.e.~in the
 factors $Ric_{ab}={R^k}_{akb}$ and $R={R^{ab}}_{ab}$ we 
 {\it also count} the internal contractions in $R$
 itself)--notice that in our definition any internal contractions in the factor
 $\nabla^{(p)}\phi$ do {\it not} count towards $\delta$.
 We will then denote by $\Sum_{\delta \ge t} T$ a
 {\it generic} linear combination  of tensor fields
in any of the forms above other than (\ref{william6}) (and
 its analogues when we replace $\nabla^{(m)}R_{ijkl}$ by
$\nabla^{(m)}Ric\otimes g$ or $\nabla^{(m)}R\otimes g$), with
 at least $t$ internal
contractions in the curvature factor. We will also denote by
$\Sum_{irrelevant} T$ a generic linear combination of tensor
fields in the form (\ref{william6}) or its analogues when we
replace $\nabla^{(m)}R_{ijkl}$ by $\nabla^{(m)}Ric\otimes g$ or
$\nabla^{(m)}R\otimes g\otimes g$.
\newline

\par We will first use the notation above in studying the
 transformation law of normalized factors
$\nabla^{r_{a_1}\dots r_{a_t}}\nabla^{(m)}_{r_1\dots r_m}W_{ijkl}$,
where
 each of the upper indices ${}^{r_{a_h}}$ contracts against a lower
derivative index ${}_{r_{a_h}}$.\footnote{Thus the above factor has $t$ internal
contractions in total.} Using the transformation law for the Weyl curvature and the
 Levi-Civita connection (see the subsection 2.3 in \cite{alexakis1}), 
we calculate that modulo partial
contractions of length $\ge 3$:

\begin{equation}
\label{lahague} \begin{split} &Image^1_\phi [\nabla^{r_{a_1}\dots
r_{a_t}}\nabla^{(m)}_{r_1\dots r_m}W_{ijkl}] =(t\cdot
(n-2)-4{t\choose{2}}) \nabla^{a_2\dots a_t}\nabla^{(m)}_{r_1 \dots
r_m}R_{ijkl}\nabla^{a_1}\phi
\\& +\Sum_{\delta=t-1,\Delta} T
+\Sum_{\delta \ge t} T+\Sum_{irrelevant} T.
\end{split}
\end{equation}
The sublinear combination $\Sum_{\delta=t-1,\Delta} T$ stands for
a generic
 linear combination of partial contractions in the
 form (\ref{william3}) with $t-1$ internal contractions.

\par We now consider a factor
$T=\nabla^{r_{a_1}\dots r_{a_t}}\nabla^{(m)}_{r_1\dots r_m}
W_{r_{m+1}r_{m+2}r_{m+3}r_{m+4}}$ where exactly one of the indices
${}^{r_{a_1}},\dots ,{}^{r_{a_t}}$ is contracting against one of
the internal indices in the factor
$W_{r_{m+1}r_{m+2}r_{m+3}r_{m+4}}$. With no loss of
 generality we assume the index ${}^{r_{a_1}}$ is contracting
 against the index ${}_i={}_{r_{m+1}}$, and we write the
 factor $T$  in the form:

$$T=\nabla^{ir_{a_2}\dots r_{a_t}}\nabla^{(m)}_{r_1\dots r_m}
W_{ijkl}.$$

\par In that case we apply the transformation law of the Weyl curvature and the Leci-Civita connection 
to derive that:

\begin{equation}
\label{lahague2} \begin{split} &Image^1_\phi
[\nabla^{ir_{a_2}\dots r_{a_t}}\nabla^{(m)}_{r_1\dots r_m} W_{ijkl}]=
(n-3)\nabla^{r_{a_2}\dots r_{a_t}}
\nabla^{(m)}_{r_1\dots r_m}R_{ijkl}\nabla^i\phi
\\&+[(t-1)(n-3)-4{t\choose{2}}\frac{n-3}{n-2}]
\nabla^{ir_{a_3}\dots r_{a_t}}\nabla^{(m)}_{r_1\dots r_m} R_{ijkl}\nabla^{r_{a_2}}\phi+
\\& \Sum_{\delta=t-1,negligible} T+\Sum_{\delta=t-1,\Delta} T+\Sum_{\delta \ge t} T+\Sum_{irrelevant} T,
\end{split}
\end{equation}
where the sublinear combination $\Sum_{\delta=t-1,negligible} T$
stands for a generic linear combination of contractions in the
form $\nabla^{ir_{a_3}\dots r_{a_t}}\nabla^{(m)}_{r_1\dots r_m} R_{ijkl}\nabla^j\phi$ 
with $t-1$ internal contractions (notice that $\nabla\phi$ is contracting against the 
index ${}_j$ and the index ${}_i$ is also involved in an internal contraction).  
 The sublinear
combination $\Sum_{\delta=t-1,\Delta} T$ stands for a generic
 linear combination of partial contractions in the
 form (\ref{william3}) with $t-1$ internal contractions.

\par Finally, we consider the transformation law of
factors \\$T=\nabla^{r_{a_1}\dots r_{a_t}}\nabla^{(m)}_{r_1\dots r_m}
W_{r_{m+1}r_{m+2}r_{m+3}r_{m+4}}$ where more than one of the
 indices ${}^{r_{a_1}},\dots ,{}^{r_{a_t}}$ is contracting
 against one of the internal indices in 
$W_{r_{m+1}r_{m+2}r_{m+3}r_{m+4}}$. In view of the anti-symmetry
of the indices ${}_i,{}_j$ and ${}_k,{}_l$, we may assume that
there are
 exactly two of the indices ${}^{r_{a_x}}$ contracting against
  internal indices in the factor $W_{ijkl}$ (modulo introducing quadratic correction terms);
   moreover (for the same reasons),
   we may assume with no loss of generality that
the indices ${}^{r_{a_1}}, {}^{r_{a_2}}$ are contracting
against the indices ${}_i,{}_k$ and also that the tensor $T$ is
symmetric with respect to the indices ${}_j,{}_l$ (this is because
of the first
 Bianchi and the antisymmetry of ${}_i,{}_j$).
 Thus
in this setting we write the factor $T$ in the form:
$$T=\nabla^{ikr_{a_3}\dots r_{a_t}}\nabla^{(m)}_{r_1\dots r_m}
W_{ijkl}.$$

We then calculate:
\begin{equation}
\label{lahague3} \begin{split} &Image^1_\phi
[\nabla^{ikr_{a_3}\dots r_{a_t}}\nabla^{(m)}_{r_1\dots r_m}
W_{ijkl}]= (n-4)(\frac{n-3}{n-2})\nabla^{kr_{a_3}\dots r_{a_t}}\nabla^{(m)}_{r_1\dots r_m}
R_{ijkl}\nabla^i\phi
\\&+(n-4)(\frac{n-3}{n-2})\nabla^{ir_{a_3}\dots r_{a_t}}\nabla^{(m)}_{r_1\dots r_m}
R_{ijkl}\nabla^k\phi
\\&+(t-2)(n-3)\nabla^{ikr_{a_4}\dots r_{a_t}}\nabla^{(m)}_{r_1\dots r_m}
R_{ijkl}\nabla^{r_{a_3}}\phi
\\&-4(\frac{n-3}{n-2})({{\delta-2}\choose{2}}+2{{\delta-1}\choose2})
\nabla^{ikr_{a_4}\dots r_{a_t}}\nabla^{(m)}_{r_1\dots r_m}
R_{ijkl}\nabla^{r_{a_3}}\phi
\\&+\Sum_{\delta=t-1,\Delta} T +\Sum_{\delta=t-1,negligible} T+\Sum_{\delta=t-1,*} T
+\Sum_{\delta \ge t}T+\Sum_{irrelevant} T.
\end{split}
\end{equation}
Here $\Sum_{\delta=t-1,*} T$ stands for a generic linear combination of contractions in the form 
$\nabla^{s_a}\phi\nabla^{(m+t+1)}_{s_1\dots s_{m+t+1}}R$ where $\nabla^{(m+t+1)}_{s_1\dots s_{m+t+1}}R$
is a factor of the scalar curvature, the factor $\nabla^{s_a}\phi$ contracts against that factor 
and the total number of  
internal contractions in this term\footnote{Including the two 
in $R={R^{ab}}_{ab}$ itself.} is $\delta-1$. 
The other linear combinations in the last line of (\ref{lahague3}) follow 
 the same notational conventions as for the RHS of
(\ref{lahague2})).

Let us finally recall from the subsection ``Technical tools'' in \cite{alexakis1} 
the formulas concerning the decomposition of iterated covariant derivatives of the Weyl tensor.

\par We make two notes regarding (\ref{lahague}) and
(\ref{lahague2}): Firstly the coefficients $(t\cdot
(n-2)-4{t\choose{2}})$, $(n-3)$, etc.
 are independent of
$m$. Secondly, since $n>4$ they are nonzero:  this is clear
because by the weight restriction we must have $t\le
\frac{n-4}{2}$.

\par These identities will be used later on in this section.

\subsection{The main algebraic Propositions.}
\label{imp.tools}

Propositions \ref{pregiade2} and \ref{tool'} are the main  tools
we will need for this section. (As explained in the introduction,
 these two Propositions, together with
the main algebraic Proposition 5.1 in \cite{alexakis1} will be proven in the series of papers 
\cite{alexakis4,alexakis5,alexakis6}).
 We will be considering tensor fields
$C^{i_1\dots i_\alpha}_{g}(\Omega_1,\dots , \Omega_p,\phi)$ of length $\sigma+1$
(with no internal contractions) in
the  form:

\begin{equation}
\label{preform2}
pcontr(\nabla^{(m_1)}R_{ijkl}\otimes\dots\otimes\nabla^{(m_s)}R_{ijkl}
\otimes
 \nabla^{(b_1)}\Omega_1\otimes\dots\otimes \nabla^{(b_p)}\Omega_p\otimes\nabla\phi).
\end{equation}
Here $\sigma=s+p$ and ${}_{i_1},\dots ,{}_{i_\alpha}$ are the free indices.

\begin{definition}
\label{preaccept}
A complete or partial contraction in the form (\ref{preform2}) will be called
 {\it acceptable} if the following conditions hold:

 \begin{enumerate}
 \item{Each $b_i\ge 2$, (in other forms each function $\Omega_h$
 is differentiated at least twice).}
 \item{No factors have internal contractions.\footnote{I.e.~no two indices 
in the same factor contract against each other.}}
 \item{The factor $\nabla\phi$ is contracting against 
some other factor (in other words the index ${}_a$ in
$\nabla_a\phi$ is not a free index).}
 \end{enumerate}
\end{definition}

\par The above definition (and also a generalized version of it)
will be used very often both in this paper and in the following ones.
For the purposes of the next Proposition we will introduce two more
definitions which will be used only in the present paper. (For the reader's
convenience we specify that the notion of {\it type} that we introduce below is a 
special case of the notion of {\it weak character} which we introduce in the paper
\cite{alexakis4}).

\begin{definition}
\label{types}
A complete or partial contraction in the form (\ref{preform2})
will have type A if the
factor $\nabla\phi$ is contracting against a factor $\nabla^{(m)}R_{ijkl}$.
It will have type B
if the factor $\nabla\phi$ is contracting against a factor $\nabla^{(B)}\Omega_h$.

\par We will divide the type A contractions into two subcategories: We will 
say a (complete or partial) contraction has type A1 if the 
factor $\nabla\phi$ is contracting against an internal index in a factor
$\nabla^{(m)}R_{ijkl}$,\footnote{(In other words, against one
of the indices ${}_i,{}_j,{}_k,{}_l$.)} has type A2 if it is contracting against a derivative index.

\par We will also say that a (complete or partial) contraction of type
 A1 has a removable index if $\sum m_i+\sum b_i>2p$ in the notation of (\ref{preform2}).
\end{definition}

Given an acceptable partial contraction with $\sigma+1$ factors and 
$a$ free indices, 
$C^{i_1\dots i_a}_g(\Omega_1,\dots,\Omega_p,\phi)$
in the form (\ref{preform2}) and any of its free indices ${}_{i_s}$, we recall that 
$div_{i_s}C^{i_1\dots i_a}_g(\Omega_1,\dots,\Omega_p,\phi)$ 
can be written as a sum of $\sigma+1$
partial contractions with $(a-1)$ free indices each; 
the $t^{th}$ summand corresponds to the 
$(a-1)$-partial contraction that arises when the 
derivative $\nabla^{i_s}$ hits the $t^{th}$ factor.

\begin{definition}
\label{Xdiv} 
Given $C^{i_1\dots i_a}_g(\Omega_1,\dots,\Omega_p,\phi)$ a partial 
contraction in the form (\ref{preform2}) as above, we define 
$Xdiv_{i_s}C^{i_1\dots i_a}_g(\Omega_1,\dots,\Omega_p,\phi)$ 
to stand for the sublinear combination in 
$div_{i_s}C^{i_1\dots i_a}_g(\Omega_1,\dots,\Omega_p,\phi)$ which arises
when the derivative $\nabla^{i_s}$ is not allowed to hit the 
factor $\nabla\phi$, nor the factor to which the index ${}_{i_s}$ belongs. 
\end{definition}

 \par Before stating our Proposition we 
 make two notes regarding the notion of {\it type}, for the reader's convenience.
 Firstly, observe that if a tensor field 
$C^{i_1\dots i_a}_g(\Omega_1,\dots,\Omega_p,\phi)$ in the 
form (\ref{preform2}) is of type A or B then all the complete contractions in the sum
 $Xdiv_{i_1}\dots Xdiv_{i_a}C^{i_1\dots i_a}_g(\Omega_1,\dots,\Omega_p,\phi)$
will also be of type A or B. Secondly: Consider a set of complete
contractions, $\{ C^l_g(\Omega_1,\dots,\Omega_p)\}_{l\in L}$, where all the complete
 contractions are in the form (\ref{preform2}) with a given
 number $\sigma_1$ of factors
 $\nabla^{(m)}R_{ijkl}$ and a given number $P$ or
 factors $\nabla^{(B)}\Omega_h$. Let $L^A\subset L$ stand for the index set
of complete contractions of type A and $L^B\subset L$ be 
the index set of complete contractions of type B. Assume an equation:

$$\sum_{l\in L} a_l C^l_g(\Omega_1,\dots,\Omega_p)=0,$$
which is assumed to hold modulo complete contractions of length $\ge\sigma+2$.
Then, using the fact that the above holds formally, we can easily derive that:
$$\sum_{l\in L^A} a_l C^l_g(\Omega_1,\dots,\Omega_p)=0,$$
$$\sum_{l\in L^B} a_l C^l_g(\Omega_1,\dots,\Omega_p)=0.$$
(Both the above equations hold modulo complete
contractions of length $\ge\sigma+2$).
\newline

\par We are now ready to state our two main algebraic propositions:

\begin{proposition}
\label{pregiade2}
 Consider two linear combinations of acceptable tensor fields in
the form (\ref{preform2}), either all in the form B or all in the form A2,

$$\Sum_{l\in L_1} a_l
C^{l,i_1\dots i_{\alpha}}_{g} (\Omega_1,\dots
,\Omega_p,\phi),$$

$$\Sum_{l\in L_2} a_l
C^{l,i_1\dots i_{\beta_l}}_{g} (\Omega_1,\dots
,\Omega_p,\phi),$$
 where each tensor field above has length $\sigma+1\ge 4$ and a given
number $\sigma_1$  of factors in the form $\nabla^{(m)}R_{ijkl}$.
 We assume  that for each
$l\in L_2$,  $\beta_l\ge \alpha+1$.  We assume that modulo complete contractions
of length $\ge\sigma+2$:

\begin{equation}
\label{hypothese2}
\begin{split}
&\Sum_{l\in L_1} a_l Xdiv_{i_1}\dots Xdiv_{i_{\alpha}}
C^{l,i_1\dots i_{\alpha}}_{g} (\Omega_1,\dots
,\Omega_p,\phi)+
\\&\Sum_{l\in L_2} a_l Xdiv_{i_1}\dots Xdiv_{i_{\beta_l}}
C^{l,i_1\dots i_{\beta_l}}_{g} (\Omega_1,\dots
,\Omega_p,\phi)=0.
\end{split}
\end{equation}

We claim that there is a linear
combination, $\sum_{r\in R} a_r
C^{r,i_1\dots i_{\alpha+1}}_g(\Omega_1,\dots,\Omega_p)$,
 of acceptable $(\alpha+1)$-tensor fields
in the form (\ref{preform2}), all with length $\sigma+1$
 and $\sigma_1$ factors $\nabla^{(m)}R_{ijkl}$ and
all of type B or type A2 respectively 
 so that:\footnote{Recall that given
a $\beta$-tensor field $T^{i_1,\dots, i_\alpha\dots i_\beta}$,
$T^{(i_1\dots i_\alpha)\dots i_\beta}$
stands for a new tensor field that arises from
$T^{i_1,\dots, i_\alpha\dots i_\beta}$ by
symmetrizing over the indices ${}^{i_1},\dots,{}^{i_\alpha}$.}

\begin{equation}
\label{concln}
\begin{split}
&\Sum_{l\in L_1} a_l C^{l,(i_1\dots i_{\alpha})}_{g} (\Omega_1,\dots
,\Omega_p)=
\sum_{r\in R} a_r Xdiv_{i_{\alpha+1}}C^{r,(i_1\dots i_\alpha)i_{\alpha+1}}_g(\Omega_1,\dots,\Omega_p),
\end{split}
\end{equation}
modulo terms of length $\ge\sigma+2$.
\end{proposition}

Second main algebraic Proposition:

\begin{proposition}
\label{tool'}
Consider two linear combinations of acceptable tensor fields in
the form (\ref{preform2}),  all in the form A1,

$$\Sum_{l\in L_1} a_l
C^{l,i_1\dots i_{\alpha}}_{g} (\Omega_1,\dots
,\Omega_p,\phi),$$

$$\Sum_{l\in L_2} a_l
C^{l,i_1\dots i_{\beta_l}}_{g} (\Omega_1,\dots
,\Omega_p,\phi),$$
 where each tensor field above has length $\sigma+1\ge 4$ and a given
number $\sigma_1$  of factors in the form $\nabla^{(m)}R_{ijkl}$.
 We assume  that for each $l\in L_2$,  $\beta_l\ge \alpha+1$. Let $\sum_{j\in J} a_j C^j_g
(\Omega_1,\dots,\Omega_p,\phi)$ stand for a generic linear combination of
complete contractions in the form (\ref{preform2}) of type A2.

  We assume that modulo complete contractions
of length $\ge\sigma+2$:

\begin{equation}
\label{hypothese2'}
\begin{split}
&\Sum_{l\in L_1} a_l Xdiv_{i_1}\dots Xdiv_{i_{\alpha}}
C^{l,i_1\dots i_{\alpha}}_{g} (\Omega_1,\dots
,\Omega_p,\phi)+
\\&\Sum_{l\in L_2} a_l Xdiv_{i_1}\dots Xdiv_{i_{\beta_l}}
C^{l,i_1\dots i_{\beta_l}}_{g} (\Omega_1,\dots
,\Omega_p,\phi)+\sum_{j\in J} a_j C^j_g(\Omega_1,\dots
,\Omega_p,\phi)=0.
\end{split}
\end{equation}
In the one case where $\alpha=1$ and the tensor fields indexed in
$L_1$ have no removable free index,\footnote{Notice that by 
definition and by weight considerations, if one of the tensor 
fields indexed in $L_1$ has this propoerty then {\it all} 
the tensor fields indexed in $L_1$ must have this property.}
 we impose the additional restriction that all the
tensor fields indexed in $L_1$ must have the free index ${}_{i_1}$ not belonging to the factor
against which $\nabla\phi$ is contracting.

\par Our claim is 
 that there is a linear
combination of acceptable $(\alpha+1)$-tensor fields
in the form (\ref{preform2}), all with length $\sigma+1$ and
$\sigma_1$ factors $\nabla^{(m)}R_{ijkl}$ and
all of type  type A1, $\sum_{r\in R} a_r
C^{r,i_1\dots i_{\alpha+1}}_g(\Omega_1,\dots,\Omega_p)$,
with length $\sigma$ so that:\footnote{Recall that given
a $\beta$-tensor field $T^{i_1,\dots, i_\alpha\dots i_\beta}$,
$T^{(i_1\dots i_\alpha)\dots i_\beta}$
stands for a new tensor field that arises from
$T^{i_1,\dots, i_\alpha\dots i_\beta}$ by
symmetrizing over the indices ${}^{i_1},\dots,{}^{i_\alpha}$.}

\begin{equation}
\label{concln2}
\begin{split}
&\Sum_{l\in L_1} a_l C^{l,(i_1\dots i_{\alpha})}_{g} (\Omega_1,\dots
,\Omega_p,\phi)=
\\&\sum_{r\in R} a_r Xdiv_{i_{\alpha+1}}C^{r,(i_1\dots i_\alpha)i_{\alpha+1}}_g
(\Omega_1,\dots,\Omega_p,\phi)+\Sum_{j\in J} a_j C^{j,(i_1\dots i_{\alpha})}_{g} (\Omega_1,\dots
,\Omega_p,\phi),
\end{split}
\end{equation}
modulo terms of length $\ge\sigma+2$. Here the tensor 
fields $C^{j,(i_1\dots i_{\alpha})}_{g}$ are all acceptable and of type A2.
\end{proposition}

{\it Note:} The conclusions of the two Propositions 
above involve (symmetric) {\it tensor fields} 
of rank $\alpha$. Clearly, if we just introduce  new
 factors $\nabla\upsilon$,\footnote{Here $\nabla\upsilon$ is some scalar-valued function.} 
and contract it against the indices ${}^{i_1},\dots,{}^{i_\alpha}$ we see that equations 
(\ref{concln}) implies:

\begin{equation}
\label{concln'}
\begin{split}
&\Sum_{l\in L_1} a_l C^{l,i_1\dots i_{\alpha}}_{g} (\Omega_1,\dots
,\Omega_p)\nabla_{i_1}\upsilon\dots\nabla_{i_\alpha}\upsilon=
\\&\sum_{r\in R} a_r Xdiv_{i_{\alpha+1}}C^{r,i_1\dots i_\alpha i_{\alpha+1}}_g
(\Omega_1,\dots,\Omega_p)\nabla_{i_1}\upsilon\dots\nabla_{i_\alpha}\upsilon,
\end{split}
\end{equation}
modulo terms of length $\ge\sigma+\alpha+2$.\footnote{Recall 
from \cite{alexakis1} that for a 
tensor field $T^{j_1\dots j_\beta}_g(\phi,\upsilon)$ involving factors 
$\nabla\phi,\nabla\upsilon$, $Xdiv_{j_h}T^{j_1\dots j_\beta}_g(\phi,\upsilon)$
stands for the sublinear combination in $div_{j_h}T^{j_1\dots j_\beta}_g(\phi,\upsilon)$ 
where the derivative $\nabla^{j_h}$ is not allowed to hit to which the 
index ${}^{j_h}$ belongs, {\it nor} any of the factors $\nabla\upsilon,\nabla\phi$.} 
Similarly for (\ref{concln2}). In fact, (\ref{concln'}) and (\ref{concln}) are equivalent, 
since (\ref{concln'}) holds formally. We will find it more 
convenient below to use (\ref{concln'}) instead of (\ref{concln}) 
(and similarly, we will use the equation that arises 
from (\ref{concln2}) by contracting 
the free indices against factors $\nabla\upsilon$).

\par Let us also state a Corollary of these two Propositions:

\begin{corollary}
\label{precorgiade}

Consider two linear combinations of acceptable tensor fields in
the form (\ref{preform2}), (with no restriction on the ``type'' of the tensor fields):

$$\Sum_{l\in L_1} a_l
C^{l,i_1\dots i_{\alpha}}_{g} (\Omega_1,\dots
,\Omega_p,\phi),$$

$$\Sum_{l\in L_2} a_l
C^{l,i_1\dots i_{\beta_l}}_{g} (\Omega_1,\dots
,\Omega_p,\phi),$$
 where each tensor field above has length $\sigma+1\ge 4$ and a given
number $\sigma_1$  of factors in the form $\nabla^{(m)}R_{ijkl}$.
 We assume  that for each $l\in L_2$,  $\beta_l\ge \alpha+1$.

  We assume that modulo complete contractions
of length $\ge\sigma+2$:

\begin{equation}
\label{hypothese2''}
\begin{split}
&\Sum_{l\in L_1} a_l Xdiv_{i_1}\dots Xdiv_{i_{\alpha}}
C^{l,i_1\dots i_{\alpha}}_{g} (\Omega_1,\dots
,\Omega_p,\phi)+
\\&\Sum_{l\in L_2} a_l Xdiv_{i_1}\dots Xdiv_{i_{\beta_l}}
C^{l,i_1\dots i_{\beta_l}}_{g} (\Omega_1,\dots
,\Omega_p,\phi)=0.
\end{split}
\end{equation}
In the one case where $\alpha=1$ and the tensor fields indexed in
$L_1$ have no removable free index, we impose the additional restriction that all the
tensor fields indexed in $L_1$ must have the free index ${}_{i_1}$ not belonging to the factor
against which $\nabla\phi$ is contracting if the tensor field
$C^{l,i_1}_g(\Omega_1,\dots,\Omega_p,\phi)$ is of type A1.

\par Our claim is then that there is a linear
combination of acceptable $(\alpha+1)$-tensor fields
in the form (\ref{preform2}), all with length $\sigma+1$ and
$\sigma_1$ factors $\nabla^{(m)}R_{ijkl}$, $\sum_{r\in R} a_r
C^{r,i_1\dots i_{\alpha+1}}_g(\Omega_1,\dots,\Omega_p)$,
with length $\sigma$ so that:

\begin{equation}
\label{concln}
\begin{split}
&\Sum_{l\in L_1} a_l C^{l,(i_1\dots i_{\alpha})}_{g} (\Omega_1,\dots
,\Omega_p,\phi)=
\sum_{r\in R} a_r Xdiv_{i_{\alpha+1}}C^{r,(i_1\dots i_\alpha) i_{\alpha+1}}_g,
\end{split}
\end{equation}
modulo terms of length $\ge\sigma+2$.
\end{corollary}

{\it Corollary \ref{precorgiade} follows from Propositons \ref{pregiade2}, \ref{tool'}:}
We prove this as follows: Divide the index sets $L_1,L_2$ into 
subsets $L_1^A,L_1^{B}$ and $L_2^A,L_2^B$, according to the following rule: 
We say $l\in L_1^{A}/l\in L_2^{A}$  if the tensor field 
$C^{l,i_1\dots i_\alpha}_g(\Omega_1,\dots,\Omega_p,\phi)$/
$C^{l,i_1\dots i_{\beta_l}}_g(\Omega_1,\dots,\Omega_p,\phi)$ 
is of type A. 
We say $l\in  L_1^{B}/l\in L_2^{B}$  if the tensor field 
$C^{l,i_1\dots i_\alpha}_g(\Omega_1,\dots,\Omega_p,\phi)$
/$C^{l,i_1\dots i_{\beta_l}}_g(\Omega_1,\dots,\Omega_p,\phi)$ (in the form (\ref{preform2}))
is of type B.\footnote{See Definition \ref{types}.} 
Observe that since (\ref{hypothese2})  holds formally, then:

\begin{equation}
\label{queenraniaA}
\begin{split}
&\Sum_{l\in L^A_1} a_l Xdiv_{i_1}\dots Xdiv_{i_{\alpha}}
C^{l,i_1\dots i_{\alpha}}_{g} (\Omega_1,\dots
,\Omega_p,\phi)+
\\&\Sum_{l\in L^A_2} a_l Xdiv_{i_1}\dots Xdiv_{i_{\beta_l}}
C^{l,i_1\dots i_{\beta_l}}_{g} (\Omega_1,\dots
,\Omega_p,\phi)=0,
\end{split}
\end{equation}

\begin{equation}
\label{queenraniaB}
\begin{split}
&\Sum_{l\in L^B_1} a_l Xdiv_{i_1}\dots Xdiv_{i_{\alpha}}
C^{l,i_1\dots i_{\alpha}}_{g} (\Omega_1,\dots
,\Omega_p,\phi)+
\\&\Sum_{l\in L^B_2} a_l Xdiv_{i_1}\dots Xdiv_{i_{\beta_l}}
C^{l,i_1\dots i_{\beta_l}}_{g} (\Omega_1,\dots
,\Omega_p,\phi)=0,
\end{split}
\end{equation}
where both the above equations hold moldulo terms of length $\ge\sigma+2$.
We will prove the claim of our Corollary 
{\it separately} for the two sublinear combinations indexed in $L_1^A, L_2^A$. 
We first prove (\ref{concln})  for the sublinear combination indexed in $L_1^A$:

\par We just appy Proposition \ref{pregiade2} to (\ref{queenraniaB}) to derive that 
there exists a tensor field $\sum_{r\in R^A} a_r
C^{r,i_1\dots i_{\alpha+1}}_g(\Omega_1,\dots,\Omega_p,\phi)$, 
as required by our Corollary, so that:

\begin{equation}
\label{kythira}
\begin{split}
&\Sum_{l\in L^B_1} a_l 
C^{l(,i_1\dots i_{\alpha})}_{g} (\Omega_1,\dots
,\Omega_p,\phi)=\sum_{r\in R} a_r Xdiv_{i_{\alpha+1}}C^{r,(i_1\dots i_\alpha) i_{\alpha+1}}_g,
\end{split}
\end{equation}
modulo terms of length $\ge\sigma+1$. This proves (\ref{concln}) for the sublinear combination 
$\Sum_{l\in L^B_1} a_l 
C^{l(,i_1\dots i_{\alpha})}_{g}$.
\newline

\par Now, we prove (\ref{concln}) for the sublinear combination indexed in $L_1^A$. We claim:

\begin{lemma}
\label{cruch}
Assume (\ref{queenraniaA}) and Proposition \ref{tool'}.
Denote by $L_1^{A,*}\subset L_1^A$, $L_2^{A,*}\subset L_2^A$ 
the index set of tensor fields of type $A_1$.
We claim: 

\begin{enumerate}
 \item There exists a   tensor field $\sum_{r\in R^A} a_r
C^{r,i_1\dots i_{\alpha+1}}_g(\Omega_1,\dots,\Omega_p,\phi)$, 
as required by  Corollary \ref{precorgiade}, so that:

\begin{equation}
\label{kythira2}
\begin{split}
&\Sum_{l\in L^{A,*}_1} a_l 
C^{l(,i_1\dots i_{\alpha})}_{g} (\Omega_1,\dots
,\Omega_p,\phi)=
\sum_{r\in R^A} a_r Xdiv_{i_{\alpha+1}}C^{r,(i_1\dots i_\alpha) i_{\alpha+1}}_g
\\&+\sum_{t\in T} a_t C^{l(,i_1\dots i_{\alpha})}_{g} (\Omega_1,\dots
,\Omega_p,\phi),
\end{split}
\end{equation}
where the partial contractions indexed in $T$ in the RHS are all acceptable  in the form 
(\ref{preform2}) with length $\sigma,\sigma_1$ factors 
$\nabla^{(m)}R_{ijkl}$ and have type $A2$. The above
holds modulo partial contractions of length $\ge\sigma+2$.

\item Assume (\ref{hypothese2}), with the additional assumption that 
$L_1^{A,*}=\emptyset$ we claim that we can write:
\begin{equation}
\label{kythira3}
\begin{split}
&\Sum_{l\in L^{A,*}_2} a_l 
Xdiv_{i_1}\dots Xdiv_{i_{\beta_l}}C^{l,i_1\dots i_{\beta_l}}_{g} (\Omega_1,\dots
,\Omega_p,\phi)=
\\&\sum_{t\in T'} a_r Xdiv_{i_1}\dots Xdiv_{i_{\beta_l}} 
C^{r,i_1\dots i_\alpha i_{\alpha+1}}_g,
\end{split}
\end{equation}
where the complete contractions indexed in $T'$ in the RHS are all in the form 
(\ref{preform2}) with length $\sigma,\sigma_1$ factors $\nabla^{(m)}R_{ijkl}$, rank $\beta_l\ge \alpha+1$ 
 and have type $A2$.The above
holds modulo complete contractions of length $\ge\sigma+2$.
\end{enumerate}
\end{lemma}
Observe that proving the above Lemma would imply Corollary \ref{precorgiade}: 

We first prove (\ref{kythira2}); we invoke the last Lemma 
in the Appendix in \cite{alexakis1} to derive:

\begin{equation}
\label{kythira2'}
\begin{split}
&\Sum_{l\in L^{A,*}_1} a_l 
Xdiv_{i_1}\dots Xdiv_{i_\alpha}C^{l,i_1\dots i_{\alpha}}_{g} (\Omega_1,\dots
,\Omega_p,\phi)=
\\&\sum_{r\in R^A} a_r Xdiv_{i_1}\dots Xdiv_{i_\alpha} Xdiv_{i_{\alpha+1}}
C^{r,i_1\dots i_\alpha i_{\alpha+1}}_g(\Omega_1,\dots
,\Omega_p,\phi)
\\&+\sum_{t\in T} a_t Xdiv_{i_1}\dots Xdiv_{i_\alpha}C^{l,i_1\dots i_{\alpha}}_{g} (\Omega_1,\dots
,\Omega_p,\phi),
\end{split}
\end{equation}
modulo complete contractions of length $\ge\sigma+2$. Now, replacing 
the above into (\ref{hypothese2}), we derive
a new local equation:

\begin{equation}
\label{kappos}
\begin{split}
& \Sum_{l\in L^{A,*}_2} a_l 
Xdiv_{i_1}\dots Xdiv_{i_{\beta_l}}C^{l,i_1\dots i_{\beta_l}}_{g} (\Omega_1,\dots
,\Omega_p,\phi)
\\&+\sum_{r\in R^A} a_r Xdiv_{i_1}\dots Xdiv_{i_\alpha} Xdiv_{i_{\alpha+1}}
C^{r,i_1\dots i_\alpha i_{\alpha+1}}_g(\Omega_1,\dots
,\Omega_p,\phi)
\\&+\sum_{t\in T} a_t Xdiv_{i_1}\dots 
Xdiv_{i_\alpha}C^{l,i_1\dots i_{\alpha}}_{g} 
(\Omega_1,\dots,\Omega_p,\phi)
\\&+\Sum_{l\in L^A_1\setminus L_1^{A,*}} a_l 
Xdiv_{i_1}\dots Xdiv_{i_{\alpha}}C^{l,i_1\dots i_{\alpha}}_{g} 
(\Omega_1,\dots,\Omega_p,\phi)=0,
\end{split}
\end{equation}
which holds modulo terms of length $\ge\sigma+2$. 
Now, applying (\ref{kythira3}) to the 
above,\footnote{We treat the sum of $Xdiv$'s 
indexed in $T, L^A_1\setminus L_1^{A,*}$ as 
a generic linear combination $\sum_{j\in J}$ 
(see the statement of Proposition \ref{tool'}).} 
we derive that we can write: 

\begin{equation}
\label{kappos2}
\begin{split}
& \Sum_{l\in L^{A,*}_2} a_l 
Xdiv_{i_1}\dots Xdiv_{i_{\beta_l}}C^{l,i_1\dots i_{\beta_l}}_{g} (\Omega_1,\dots
,\Omega_p,\phi)
\\&+\sum_{r\in R^A} a_r Xdiv_{i_1}\dots Xdiv_{i_{\alpha+1}}
C^{r,i_1\dots i_\alpha i_{\alpha+1}}_g(\Omega_1,\dots
,\Omega_p,\phi)
\\&=\sum_{t\in \tilde{T}} a_t Xdiv_{i_1}\dots Xdiv_{i_{\beta_l}}
C^{t,i_1\dots i_{\beta_l}}_{g} (\Omega_1,\dots
,\Omega_p,\phi),
\end{split}
\end{equation}
which holds modulo complete contractions of length
$\sigma+2$. The complete contractions are as the ones indexed in $T'$
in the RHS of (\ref{kythira3}). 

\par Now, replacing the above into (\ref{kappos}) we derive  a new local equation, 
where all the tensor fields of length $\sigma+1$ are of type $A2$. Thus, we 
are in a position to apply Proposition \ref{pregiade2}. 
We derive that there exists a linear combination of $(\alpha+1)$-tensor fields 
 as in the claim of Proposition \ref{pregiade2}, indexed im $R^+$ below, so that:

\begin{equation}
\label{kappos4}
\begin{split}
& \sum_{L_1^A\setminus L_1^{A,*}} a_l C^{l,i_1\dots i_{\beta_l}}_{g} (\Omega_1,\dots
,\Omega_p,\phi)+\sum_{t\in T} a_t C^{t,(i_1\dots i_{\alpha})}_{g} 
(\Omega_1,\dots,\Omega_p,\phi)=
\\&\sum_{r\in R^+} a_r Xdiv_{i_{\alpha+1}}
C^{r,i_1\dots i_{\alpha}i_{\alpha+1}}_{g} 
(\Omega_1,\dots,\Omega_p,\phi),
\end{split}
\end{equation}
modulo terms of length $\ge\sigma+2$. Now, adding (\ref{kythira2}) and (\ref{kappos4})
we derive Corollary (\ref{precorgiade}). $\Box$ 
\newline

{\it Proof of Lemma \ref{cruch}:} The first claim 
follows by applying Proposition 
\ref{pregiade2} to (\ref{queenraniaA}).\footnote{We treat the sum 
of $Xdiv$'s indexed in $L_2$ as a sum $\sum_{j\in J}\dots$,
as in the statement of \ref{tool'}.} 

To derive the second claim, we proceed by induction: 
Let $\beta_{min}>\alpha$ be the minimum 
rank among the tensor fields indexed in $L_1^{A,*}$; denote the corresponding index set by 
$L_1^{A,*,\alpha,\beta_{min}}$.
We apply Proposition \ref{tool'} to (\ref{hypothese2})
(where the minimum rank among the tensor fields 
indexed in $L_1^{A,*}$ is now $\beta_{min}$);
we derive that there exists a linear combination 
of partial contractions in the form (\ref{preform2}) 
of type $A_1$, each with length $\sigma$, 
with $\sigma_1$ factors $\nabla^{(m)}R_{ijkl}$ and 
rank $\beta_{min}+1$, indexed in $R^A$ below, such that: 

\begin{equation}
\label{kythira2'}
\begin{split}
&\Sum_{l\in L^{A,*,\beta_{min}}_1} a_l 
C^{l,(i_1\dots i_{\beta_{min}})}_{g} (\Omega_1,\dots
,\Omega_p,\phi)=
\\&\sum_{r\in R^A} a_r  Xdiv_{i_{\beta_{min}+1}}
C^{r,(i_1\dots i_{\beta_{min}}) i_{\beta_{min}+1}}_g+
\sum_{t\in T^A} a_t C^{l,(i_1\dots i_{\beta_{min}})}_{g} (\Omega_1,\dots
,\Omega_p,\phi),
\end{split}
\end{equation}
where the terms indexed in $T^A$ have all the features of 
the tensor fields indexed in $L^{A,*,\beta_{min}}_1$, 
{\it but} they are of type $A2$ rather than $A1$. Iteratively applying the above step, 
we derive (\ref{kythira3}).\footnote{Notice that since all local equations 
we deal with involve complete contractions of a fixed weight $-n$, the 
maximum possible number of free indiexes among tensor fields in the 
form (\ref{preform2}) is $\frac{n}{2}$. Therefore we derive 
(\ref{kythira3}) after a {\it finite} number of iterations.}  $\Box$

\subsection{Proof of Lemma \ref{nabothewri2}.}

\par Recall that we are assuming that all the complete contractions
$C^l(g), l\in L_\sigma^j$ in $P(g)$\footnote{Recall that
$L_\sigma\subset L$ stands for the index set of complete
contractions in $P(g)$ with (minimum) length $\sigma$. Recall that
$L^j_\sigma\subset L_\sigma$ stands for the index set of complete
contractions in $\sum_{l\in L_\sigma} a_l C^l(g)$ with the (minimum) number $j$ of
internal contractions.}  have  at least
two of their internal contractions belonging to different factors. Observe
 that since $j\ge 2$ we must have $\sigma\le\frac{n}{2}-2$.

\par We will again consider $I^1_g(\phi)(:=\frac{d}{dt}|_{t=0}[e^{nt\phi}P(e^{2t\phi}g)])$. 
Recall that using (\ref{anatrans}), (\ref{lahague}), (\ref{lahague2}), 
(\ref{lahague3}) and the first  two formulas from 
subsection 5.1 in \cite{alexakis1}
we may express
$I^1_g(\phi)$ as a linear combination of complete contractions in
the form:

\begin{equation}
\label{skoteinh}
contr(\nabla^{(m_1)}R_{ijkl}\otimes\dots\otimes\nabla^{(m_s)}R_{i'j'k'l'}
\otimes\nabla^{(p_1)}Ric_{ab}\otimes\dots
\otimes\nabla^{(p_q)}Ric_{a'b'}\otimes\nabla^{(\nu)}\phi)
\end{equation}
 Recall
that for each $C^l(g)$ of length $\sigma$ in $P(g)$,
$Image^1_{\nabla\phi}[C^l(g)]$
stands for the sublinear combinations in 
$Image^1_\phi[C^l(g)]$\footnote{Recall that 
$Image^1_\phi[C^l(g)]=e^{n\phi}C^l(e^{2\phi}g)-C^l(g)$.} 
 of contractions (in the form
(\ref{skoteinh})) with $\nabla^{(\nu)}\phi=\nabla\phi$ or
$\nabla^{(\nu)}\phi=\Delta\phi$, respectively.

\par By virtue of the formulas in subsection \ref{prelimwork}
(see especially (\ref{anatrans}), (\ref{lahague}), (\ref{lahague2}), 
(\ref{lahague3}) and also the first  two formulas from 
subsection ``Technical Tools'' in \cite{alexakis1})
we observe that all the complete contractions
in $Image^1_{\nabla\phi}[\Sum_{l\in L_\sigma} a_l C^l(g)]$ will
have $\delta_R\ge j-1$; furthermore, (for the same reasons)  it folows that 
any complete contraction  with $\delta=j-1$ and with $\beta>0$ 
factors $\nabla^{(p)}Ric_{ab}$, then the 
indices ${}_a,{}_b$ in each such factor must contract against each other. 
We calculate:
\begin{equation}
\label{biasthka} 
Image^1_{\phi}[\Sum_{l\in L_\sigma} a_l C^l(g)]=Image^1_{\nabla\phi}[\Sum_{l\in L_\sigma} a_l C^l(g)]+
\sum_{u\in U_{\Delta\phi}} a_u C^u_g(\phi)+\Sum_{x\in X} a_x C^x_{g}(\phi),
\end{equation}
where $\sum_{u\in U_{\Delta\phi}} a_u C^u_g(\phi)$ is a generic 
linear combination of terms in the form (\ref{skoteinh}) with 
length $\sigma+1$ and a factor $\nabla^{(\nu)}\phi=\Delta\phi$. $\Sum_{x\in X}
a_x C^x_{g}(\phi)$stands for a generic linear combination
of complete contractions in the form (\ref{skoteinh})
with either length $\ge\sigma+2$ or with length
$\sigma+1$ and a factor $\nabla^{(\nu)}\phi\ne\Delta\phi$.

 On the other hand, we see that for any complete
contraction $C^l(g)$ of length $\ge\sigma+1$ in $P(g)$:
\begin{equation}
\label{dakry}
Image^1_\phi[C^l(g)]=\sum_{u\in U_{\Delta\phi}} a_u C^u_g(\phi)+\Sum_{x\in X}
a_x C^x_{g}(\phi),
\end{equation}
(with the same conventions as above). 

\par For future reference we denote by
$Image^1_{\nabla\phi,*}[P(g)|_{\sigma}]$ the sublinear combinations
in $Image^1_{\phi}[P(g)|_{\sigma}]$ that consists of complete
contractions with no factors $\nabla^{(p)}Ric$, and with a factor
$\nabla\phi$.
\newline

\par Now, we need two main technical Lemmas for this subsection. We
will be considering a linear combination $Y_{g}(\phi)$ in the
form:

\begin{equation}
\label{aintee}
\begin{split}
Y_{g}(\phi)=\Sum_{v=0}^\alpha \{ \Sum_{u\in U^v_{\nabla\phi}} a_u
C^u_{g}(\phi)+ \Sum_{u\in U^v_{\Delta\phi}} a_u C^u_{g}(\phi)\}+
\Sum_{x\in X} a_x C^x_{g}(\phi).
\end{split}
\end{equation}
All the complete contractions $C^u_{g}(\phi)$ here have $\sigma
+1$ factors. The linear combination $\sum_{x\in X} a_x C^x_g(\phi)$ is a  generic
 linear combination as defined earlier (see discussion under (\ref{biasthka})).
 Each complete contraction indexed in
 $U^v_{\nabla\phi}$, $v>0$  has a factor $\nabla\phi$ and $v$ factors
$\nabla^{(p)}Ric$ and {\it either} has $\delta_R\ge j$ {\it or} 
has $\delta_R=j-1$ but then also has 
the property that in all its $v$  factors $\nabla^{(p)}Ric_{ab}$ the indices ${}_a,{}_b$ are 
contracting against each 
other;\footnote{Equivalently, there are $v$ factors $\nabla^{(p)}R$ of the differentiated scalar curvature.} 
if $v=0$ then $\delta_R\ge
j-1$. Each complete contraction
 indexed in $U^v_{\Delta\phi}$
has a factor $\Delta\phi$ and $v$ factors $\nabla^{(p)}Ric$ and
 there is no restriction on $\delta_R$. Thus $\alpha$
 is an upper bound on the number of factors $\nabla^{(p)}Ric$ in the
complete contractions $C^u_g(\phi)$ 
in $Y_g(\phi)$.

 {\it Note:} Observe that $I^1_g(\phi)$ is of the
 form (\ref{aintee}) above, where

\begin{equation}
\label{identification}
\Sum_{u\in U^v_{\nabla\phi}} a_u C^u_{g}(\phi)=Image^1_{\nabla\phi}[P(g)|_{\sigma}].
 \end{equation}
{\it Remark:} We remark that in (\ref{identification}),
  each $C^u_g(\phi)$, $u\in U^v_{\nabla\phi}$ 
with $\beta>0$ factors $R$ (of the scalar curvature) will satisfy $\delta\ge j-1+2\beta$. 
This follows from the decomposition of the Weyl tensor, since each factor $R$ in 
$C^u_g(\phi)$ must have arisen from an (undifferentiated) factor $W_{ijkl}$ in 
$P(g)|_{\sigma}$; thus a term with no internal contractions 
gives rise to a term with two internal contractions. 
\newline

\par We claim the following:

\begin{lemma}
\label{work}
Consider $Y_{g}(\phi)$  as in (\ref{aintee}).
  Assume that $\int_{M^n}Y_{g}(\phi)
dV_{g}=0$. Let $\delta_{min}$  stand for the
minimum $\delta_R$ among the complete contractions indexed in
$U^\alpha_{\Delta\phi}\bigcup U^\alpha_{\nabla\phi}$.\footnote{In other words $\delta_{min}$ 
stands for the minimum number of internal contractions 
in curvature factors (i.e. factors in the form 
$\nabla^{(m)}R_{ijkl}, \nabla^{(p)}Ric_{ab}$),
 among the complete contractions with 
the maximum number $\alpha$ of factors in the form $\nabla^{(p)}Ric$, 
and having afactor $\Delta\phi$ or $\nabla\phi$.
Observe by definition that if $\delta_{min}<j-1$ then $U^\alpha_{\nabla\phi}=\emptyset$.}
We denote the corresponding index sets by 
$U^\alpha_{\Delta\phi,\delta_{min}}, U^\alpha_{\nabla\phi,\delta_{min}}$.
\newline

\par In the case where $\delta_{min}<j-1$,\footnote{(In which case 
$U^\alpha_{\nabla\phi}=\emptyset$ as noted in the previous footnote).} 
and in the case $\delta_{min}=j-1$ 
under the additionnal assumption that $U^\alpha_{\nabla\phi,j-1}=\emptyset$, 
we claim that
 there is a linear combination of vector fields,
$\Sum_{h\in H} a_h C^h_{g}(\Delta\phi)$, each in the
 form (\ref{skoteinh}) with $\nabla^{(\nu)}\phi=\Delta\phi$ 
and with one free index, so that:

\begin{equation}
\label{liapkin}
\begin{split}
& \Sum_{u\in U^\alpha_{\Delta\phi,\delta_{min}}} a_u
C^u_{g}(\phi)- div_i\Sum_{h\in H} a_h C^{h,i}_{g}(\Delta\phi)=
\Sum_{u\in U^\alpha_{\Delta\phi,\delta_{min}+1}} a_u
C^u_{g}(\phi)+
\\& \Sum_{u\in U^{\alpha-1}_{\Delta\phi,\delta_{min}}} a_u
C^u_{g}(\phi)+\Sum_{x\in X} a_x C^x_{g}(\phi),
\end{split}
\end{equation}
where the linear combinations in the RHS are generic 
linear combinations with the following properties:
Each contraction indexed in
$U^\alpha_{\Delta\phi,\delta_{min}+1}$ is in the form (\ref{skoteinh}), has $\alpha$ factors
$\nabla^{(p)}Ric$ and a factor $\Delta\phi$ and
$\delta_R=\delta_{min}+1$. Each complete contraction indexed
 in $U^{\alpha-1}_{\Delta\phi,\delta_{min}}$ has
$\alpha -1$ factors $\nabla^{(p)}Ric$, $\delta_R=\delta_{min}$  and a
factor $\Delta\phi$. 

\par On the other hand, if $\delta_{min}\ge j-1$, and
$\alpha>0$, we  
claim that there is a linear combination of
vector fields, $\Sum_{h\in H} a_h C^{h,i}_{g}(\phi)$, so
that:

\begin{equation}
\label{sabatto}
\begin{split}&\Sum_{u\in U^\alpha_{\nabla\phi,\delta_{min}}\bigcup
 U^\alpha_{\Delta\phi,\delta_{min}}} a_u C^u_{g}(\phi)-
div_i\Sum_{h\in H} a_h C^{h,i}_{g}(\phi)= 
\\&\Sum_{u\in U^\alpha_{\nabla\phi,\delta_{min}+1}\bigcup
 U^\alpha_{\Delta\phi,\delta_{min}+1}} a_u C^u_{g}(\phi)+
\Sum_{u\in U^{\alpha-1}_{\nabla\phi}\bigcup U^{\alpha-1}_{\Delta\phi}} 
a_u C^u_{g}(\phi)+\Sum_{x\in X} a_x
C^x_{g}(\phi),
\end{split}
 \end{equation}
where the linear combinations in the RHS are generic 
linear combinations with the following properties: The terms indexed in
$U^\alpha_{\nabla\phi,\delta_{min}+1}\bigcup
 U^\alpha_{\Delta\phi,\delta_{min}+1}$ still have $\alpha$ factors $\nabla^{(p)}Ric$ and $\delta_{min}+1$
internal contractions in total, while the ones indexed in  
$U^{\alpha-1}_{\nabla\phi}$ have $\alpha-1$ factors
$\nabla^{(p)}Ric$ and $\delta_R\ge j$ and the ones in $U^{\alpha-1}_{\Delta\phi}$ 
have $\alpha-1$ factors
$\nabla^{(p)}Ric$ (and there is no restriction on $\delta_R$).
\end{lemma}

\par We will prove Lemma \ref{work} further down. For the
 time being, we will show how it implies Lemma
\ref{nabothewri2}.
\newline

{\it Lemma \ref{work} implies Lemma \ref{nabothewri2}:}

\par Observe that iteratively applying the above Lemma, 
 starting with $Y_g(\phi)=I^1_{g}(\phi)$, we derive that there is a
 vector field, $\Sum_{h\in H} a_h
C^{h,i}_{g}(\phi)$ so that:

\begin{equation}
\label{synepws}
\begin{split}
&I^1_{g}(\phi)-div_i\Sum_{h\in H} a_h C^{h,i}_{g}(\phi)=
Image^1_{\nabla\phi,*}[P(g)|_{\sigma}]+
\\&\Sum_{u\in U^0_{\Delta\phi,\delta\ge j-1}} a_u
C^u_{g} (\phi)+\Sum_{u\in U^0_{\nabla\phi,\delta\ge j}}
 a_u C^u_{g}(\phi)+\Sum_{x\in X} a_x C^x_{g}(\phi),
\end{split}
\end{equation}
where we recall that $Image^1_{\nabla\phi,*}[P(g)|_{\sigma}]$ stands for the sublinear
combinations in $Image^1_{\phi}[P(g)|_{\sigma}]$ that consists of
complete contractions with no factors $\nabla^{(p)}Ric$, and with
a factor $\nabla\phi$.
 Also, $\Sum_{u\in U^0_{\Delta\phi,
\delta\ge j-1}} a_u C^u_{g} (\phi)$, $\Sum_{u\in
U^0_{\nabla\phi,\delta\ge j}}
 a_u C^u_{g}(\phi)$ stand for generic linear combinations of
 complete contractions in the form (\ref{skoteinh}) with a factor $\nabla\phi,\Delta\phi$
respectively $\delta_R\ge j-1$, $\delta_R\ge j$ respectively.
\newline

Now, in order
to prove Lemma \ref{nabothewri2} we examine the
 sublinear combination $Image^1_{\nabla\phi,*}
[P(g)|_{\sigma}]$. We recall that $L_\sigma^j\subset L$ stands for
the sublinear
 combination in $P(g)$ of complete contractions with length
$\sigma$ and $j$ internal contractions. We recall that for
 each $C^l(g)$, $l\in L^j_\sigma$ there must be at least two
 internal contractions belonging to different factors.
 For every
 $l\in L^j_\sigma$ we denote by
 $Image^1_{\nabla\phi,*,+}[P(g)|_\sigma]$ the sublinear
 combination in $Image^1_{\nabla\phi,*}[P(g)|_\sigma]$ with
 $\delta_R=j-1$.

\par Thus, we derive:

\begin{equation}
\label{kineziko}
\begin{split}
&Image^1_{\nabla\phi,*}[P(g)|_{\sigma}] =
Image^1_{\nabla\phi,*,+}[P(g)|_\sigma]+ \Sum_{h\in H} a_h
C^h_{g}(\phi),
\end{split}
\end{equation}
where each $C^h_{g}(\phi)$ is in the form (\ref{skoteinh}), has $\delta\ge
j$, no factors $\nabla^{(p)}Ric$ and a factor $\nabla\phi$ 
(notice this fits into the generic linear combination 
$\sum_{u\in U^0_{\nabla\phi,\delta\ge j}} a_u C^u_g(\phi)$).

\par In view of the above, we may integrate (\ref{synepws}) to derive 
a global equation of the form:

\begin{equation}
\label{boh9i}
\begin{split}
&\int_{M^n} Image^1_{\nabla\phi,*,+}[P(g)|_\sigma]+
\Sum_{u\in U^0_{\Delta\phi,\delta\ge j-1}} a_u
C^u(g) \Delta\phi
\\&+\Sum_{u\in U^0_{\nabla\phi,\delta\ge j}}
 a_u C^u_{g}(\phi)+\Sum_{x\in X} a_x C^x_{g}(\phi)
dV_{g}=0.
\end{split}
\end{equation}

\par We now show how to derive Lemma \ref{nabothewri2}
from the above. In order to do so, we will introduce some more
notational conventions. 

\begin{definition}
\label{character}
For each tensor field
 $C^{i_1\dots i_\alpha}_g$ in the form 
$$pcontr(\nabla^{(m_1)}W_{ijkl}\otimes\dots\otimes\nabla^{(m_\sigma)}W_{i'j'k'l'})$$
 we will associate a list of numbers,
which we will call the {\it character} of
$C^{l,i_1\dots i_\alpha}_g$: We list out the factors $T_1,\dots T_\sigma$ of
 $C^{l,i_1\dots i_\alpha}_g$ and we define $List(l)=(L_1,\dots ,L_\sigma)$
as follows: $L_i$ will stand for the number of free indices that
 belong to the factor $T_i$. We then define $\vec{\kappa}(l)$
 to stand for the decreasing rearrangement
 of $List(l)$.\footnote{Thus the character is ultimately a list of numbers.
 We can also refer to an abstract ``character'' $\vec{\kappa}$, which
does not necessarily need to correspond to a particular tensor field.
 We also note that two different tensor fields
  can have the same character.} With a slight abuse of notation,
 we will also consider the complete
 contractions in the index set $L_\sigma^j$
 in $P(g)=\sum_{l\in L} a_L C^l(g)$ and define their {\it character} to be
 the character of the {\it tensor field} $C^{l,i_1\dots i_j}_g$
 that arises from $C^l_g$ by replacing all the internal contractions by free indices.
\end{definition}

  We can thus group
  up the different $j$-tensor fields indexed in
  $L_\sigma^j$ according to their double characters: Let $\vec{K}$
 stand for the set of all the different
 characters appearing among the tensor fields indexed in $L_\sigma^j$.
 We write $L_\sigma^j=\bigcup_{\vec{\kappa}\in
 \vec{K}}L_\sigma^{j,\vec{\kappa}}$; here
 $L_\sigma^{j,\vec{\kappa}}\subset L_\sigma^j$
 stands for the index set of $j$-tensor fields
 with a character $\vec{\kappa}$.

 We also introduce a partial
 ordering among characters according to  lexicographic comparison:
 a character $\vec{\kappa}_2$
 is subsequent to  $\vec{\kappa}_1$ if $\vec{\kappa}_1$
 is lexicographically greater than $\vec{\kappa}_2$.

\par  Notice that any two
different characters are {\it comparable}, in the sense
 that one will be subsequent to the other.

\par Now, we return to the derivation of Lemma \ref{nabothewri2}:
We will again proceed by induction. We break up
$L^j_{\sigma}$ into subsets $L^{j,\vec{\kappa}}_\sigma$ that index
complete
 contractions $C^l(g)$ with the same character, $\vec{\kappa}$.\footnote{Recall that
the character of a complete contraction $C(g)$
(with factors $\nabla^{(m)}W_{ijkl}$)
was defined to be the character of the tensor
field that formally arises from $C(g)$ by making all
internal contractions into free indices.} There is a finite list
 of possible characters; we denote
the set of characters of the complete
contractions in $L^j_\sigma$ by $\vec{K}$.

If $\vec{K}=\emptyset$ there is nothing to prove. If
$\vec{K}\ne\emptyset$, we pick out the maximal character
$\vec{\kappa}_{max}\in \vec{K}$. We will show the claim of
 Lemma \ref{nabothewri2} for the sublinear combination
$\Sum_{l\in L^{j,\vec{\kappa}_{max}}_\sigma} a_l C^l(g)$.
In other words, we will prove that there exists a linear
combination of partial contractions, $\Sum_{h\in H} a_h
C^{h,i}(g)$, where each $C^{h,i}(g)$ is in the form (\ref{karami2})
with weight $-n+1$ and $\delta_W=j$, so that:

\begin{equation}
\label{salvationstep} \Sum_{l\in L^{j,\vec{\kappa}_{max}}_\sigma} a_l C^l(g)-div_i
\Sum_{h\in H} a_h C^{h,i}(g)=\Sum_{v\in V} a_v C^v(g)
\end{equation}
where each $C^v(g)$ is in the form (\ref{karami}) with $\delta_W\ge
j+1$. Moreover, each $C^v(g)$ will have at least two internal contractions
 belonging to different factors. The above equation holds modulo
 complete contractions of length $\ge\sigma +1$.

 Clearly,
if we can show that claim then Lemma \ref{nabothewri2} will follow
by induction.
\newline

{\it Proof of (\ref{salvationstep}):} 
We  introduce some notational normalizations.
We write out $\vec{\kappa}_{max}=(k_1,\dots ,k_a)$, where if $i<j$
then $k_i\ge k_j$.\footnote{Recall that each number $k_s$ stands for the number
of internal contractions in some given
factor in $\vec{\kappa}_{Max}$. The second restriction
 reflects the fact that the list is taken in decreasing rearrangement}
Note that $k_1$ is the maximum number of internal contractions 
that can belong to a given
factor among all complete contractions in $L^j_\sigma$.
This follows by the definition of ordering among different refined double characters.
 By the hypothesis of Lemma \ref{nabothewri2} we have that $a\ge 2$ (this reflects the fact that
not all internal contractions can belong to the same factor). Now,
let $b$ be the largest number for which $k_{a-b+1},\dots ,k_a=k_{min}$ have the same value.  
In other words, for each $l\in
L^{j,\vec{\kappa}_{max}}_\sigma$ there are $b$ different factors $\nabla^{(m)}W_{ijkl}$, 
 with $k_{min}$ internal contractions each.
  For notational convenience, we will assume that
 the last $b\cdot k_{min}$ free indices in
$C^{l,i_1\dots i_j}(g)$ correspond to those internal
 contractions (i.e. they arise from those internal contractions
  when we make the internal contractions in $C^l(g)$ into
   free indices in $C^{l,i_1\dots i_j}(g)$).

\par Now, we will say that a $(j-1)$-tensor field in the form 
$$pcontr(\nabla^{(m_1)}W_{ijkl}\otimes\dots\otimes\nabla^{(m_\sigma)}W_{i'j'k'l'}\otimes\nabla\phi)$$
has double-character $\vec{\kappa}'_{max}$ if its $j-1$ free indices are distributed among its
factors according to the pattern: $(k_1,k_2,\dots, k_{a-1}, k_{min}-1)$ (in other words one factor $T_1$
has $k_1$ free indices, the next factor $T_2$ has $k_2$ free indices etc)
{\it and} the factor with $k_a-1$ free indices is also
contracting against the factor $\nabla\phi$. As above, we will extend this definition to
complete contractions in the form 
$$contr(\nabla^{(m_1)}W_{ijkl}\otimes\dots\otimes\nabla^{(m_\sigma)}
W_{i'j'k'l'}\otimes\nabla\phi),$$ with $j-1$ internal contractions:
We will say that such a complete contraction $C_g(\phi)$ has  a double-character
 $\vec{\kappa}_{max}'$ if  the tensor field
$C^{i_1\dots i_{j-1}}_g(\phi)$ which arises from
 $C_g(\phi)$ by making all the internal contractions into free indices
has a double-character $\vec{\kappa}_{max}'$.

\par We then examine the 
complete contractions in $Image^1_{\nabla\phi,*,+}[P(g)|_\sigma]$; in particular 
 we  seek to understand the sublinear combination with double character $\vec{\kappa}_{max}'$.
Denote this 
sublinear combination by $\{Image^1_{\nabla\phi,*,+}[P(g)|_\sigma]\}_{\vec{\kappa}_{max}'}$. 
Using formulas (\ref{lahague}), (\ref{lahague2}), (\ref{lahague3}) we can straightforwardly descibe 
$Image^1_{\nabla\phi,*,+}[P(g)|_\sigma]$: 

\par This sublinear combination arises 
{\it exclusively} from the sublinear combination $\Sum_{l\in L^{j,\vec{\kappa}_{max}}_\sigma} a_l C^l(g)$ 
in $P(g)$ via the following process: 
Consider each $C^l(g)$, $l\in L^{j,\vec{\kappa}_{max}}_\sigma$ and list out its $\sigma$ factors
 $T^l_1,\dots T^l_\sigma$; for each factor $T^l_a$ we let $Subst^{\nabla\phi}_a[C^l(g)]$ 
stand for the (linear combination of) complete contractions 
that arises from $C^l(g)$ by replacing $T^l_a$ by one of the terms explicitly written 
out in the RHSs of (\ref{lahague}), (\ref{lahague2}), (\ref{lahague3})\footnote{(Which of the 
equations we use depends on how many internal contractions in $T^l_a$ 
involve internal indices).} and then replacing all the other factors $T^l_b=\nabla^{(m)}W_{ijkl}$, $b=1,\dots,a-1,a+1,\dots\sigma$,
 by either $\frac{n-3}{n-2}\nabla^{(m)}R_{ijkl}$ or $\nabla^{(m)}R_{ijkl}$ (depending on whether 
$T^l_b$ has an internal contraction involving an internal index or not). 
Then, just by virtue of the formulas (\ref{lahague}), 
(\ref{lahague2}), (\ref{lahague3}) and the first two 
formulas in the subsection 5.1 in \cite{alexakis1}, we calculate: 

\begin{equation}
\begin{split}
\label{peaceatlast} 
&Image^1_{\nabla\phi,*,+}[P(g)|_\sigma]=
\sum_{l\in L^{j,\vec{\kappa}_{max}}_\sigma} a_l \sum_{a=1}^\sigma Subst^{\nabla\phi}_a[C^l(g)]+
\sum_{negligible} C_g(\nabla\phi),
\end{split}
\end{equation}
where $\sum_{negligible} C_g(\nabla\phi)$ stands for a generic linear combination of complete
contractions in the form 
$contr(\nabla^{(m_1}R_{ijkl}\otimes\dots\otimes\nabla^{(m_\sigma)}R_{i'j'k'l'}\otimes\nabla\phi)$ with $j-1$ 
internal contractions in total and with the factor $\nabla\phi$ contracting 
against the index ${}_j$ in a factor $T=\nabla^{(m)}R_{ijkl}$ for which the index ${}_i$ 
is internally contracting in $T$.

{\it Note:} Notice that in the first linear combination in the RHS 
there can be at most $k_1$ internal contractions 
in any of the factors in any $C_g(\nabla\phi)$. 

Given (\ref{peaceatlast}) we are able 
 to explicitly write out the sublinear combination 
$\{Image^1_{\nabla\phi,*,+}[P(g)|_\sigma]\}_{\vec{\kappa}_{max}'}$ in  
$Image^1_{\nabla\phi,*,+}[P(g)|_\sigma]$ of terms with a 
double character $\vec{\kappa}_{max}'$:

\begin{equation}
\label{jeep}
\begin{split}
 &\{Image^1_{\nabla\phi,*,+}[P(g)|_\sigma]\}_{\vec{\kappa}_{max}'}=
\sum_{l\in L^{j,\vec{\kappa}_{max}}_\sigma} a_l \sum_{f=\sigma-b\cdot k_{min}+1}^\sigma 
(Const)_{l,f} C^{l,\iota|i_f}(g)\nabla_{i_f}
\\&+\sum_{negligible} C_g(\nabla\phi).
\end{split}
\end{equation}
Here $C^{l,\iota|i_f}(g)$ stands for the vector field that arises from $C(g)$ by formally replacing all
 factors $\nabla^{(m)}W_{ijkl}$ by $\nabla^{(m)}R_{ijkl}$  and then making the $f^{th}$  internal 
contraction into a free index ${}_{i_f}$. 
\newline

\par Although the full description of the exact form of 
$\{Image^1_{\nabla\phi,*,+}[P(g)|_\sigma]\}_{\vec{\kappa}_{max}'}$ 
is rather messy, all that is really important here is to note 
that by applying the operation $Weylify$ (see subsection 5.1 in \cite{alexakis1}) to this 
sublinear combination, one can {\it recover} the sublinear combination 
$\sum_{l\in L^{j,\vec{\kappa}_{max}}_\sigma} a_l  C^l(g)$ in $P(g)$: 

\par Let us define a formal  opeartion $Op$ which acts on the complete contractions above 
by replacing all internal contractions by factors $\nabla\upsilon$, setting $\phi=\upsilon$ and then 
acting on the resulting complete contractions by the operation $Weylify$. 
It follows from the discussion above (\ref{peaceatlast}), (\ref{jeep}) 
and  the definition of the operation $Weylify$ in \cite{alexakis1} that:

\begin{equation}
\label{themainthing} 
Op[\{Image^1_{\nabla\phi,*,+}[P(g)|_\sigma]\}_{\vec{\kappa}_{max}'}]=
\Sum_{l\in L^{j,\vec{\kappa}_{max}}_\sigma} a_l\cdot 
[b\cdot(k_{min}\cdot ((n-2)-4{{k_{min}\choose{2}}})] C^l(g).
\end{equation}
Two notes are in order: Firstly 
 that the coefficient $[b\cdot(k_{min}\cdot ((n-2)-4{{k_{min}\choose{2}}})]$ is {\it universal}, 
i.e. it depends only on the character $\vec{\kappa}_{max}$ 
(thus it can be factored out in the RHS of 
(\ref{themainthing})), and secondly that it is non-zero: 
This follows from the simple observation that $b\ge 1$ by 
definition and $l\le \frac{n}{2}-3$, since the complete 
contractions $C^l(g)$  have weight $-n$ and  involve 
factors $\nabla^{(m)}W_{ijkl}$ and $\sigma\ge 3$. 
\newline

\par We now consider the linear combination
$\Sum_{u\in U^0_{\Delta\phi,\delta\ge j-1}} a_u C^u_{g} (\phi)$.
We denote by $\Sum_{u\in U^0_{\Delta\phi,\delta= j-1}} a_u C^u_{g}
(\phi)$ the sublinear combination that
 consists of complete contractions with $\delta_R=j-1$.

\par Now, let us also denote by $A=k_1$, where $k_1$ is the
 first number on the list $(k_1,\dots ,k_a)$ for $\vec{\kappa}_{max}$ above.
We will pay special attention to the complete contractions in
$\Sum_{u\in U^0_{\Delta\phi,\delta= j-1}} a_u C^u_{g} (\phi)$
which have a factor $\nabla^{(m)}R_{ijkl}$ with $A$ internal
 contractions. We denote the index set of those complete
 contractions by $SU^0_{\Delta\phi,\delta= j-1}\subset U^0_{\Delta\phi,\delta= j-1}$.

\par We claim that there is a linear combination of vector
fields $C^{h,i}(g)\Delta\phi$ (in the form (\ref{skoteinh}) without factors $\nabla^{(p)}Ric$ and with $\nabla^{(\nu)}\phi=\Delta\phi$), so that:

\begin{equation}
\label{elaela} \Sum_{u\in SU^0_{\Delta\phi,\delta= j-1}} a_u
C^u_{g}\cdot \Delta\phi-div_i \Sum_{h\in H} a_h
C^{h,i}(g)\Delta\phi= \Sum_{u\in U^0_{\Delta\phi,\delta= j}} a_u
C^u_{g}\cdot \Delta\phi+ \Sum_{x\in X} a_x C^x_{g}(\phi),
\end{equation}
where $\Sum_{u\in U^0_{\Delta\phi,\delta= j}} a_u C^u_{g}\cdot
\Delta\phi$ stands for a generic linear
 combination of complete contractions in the form
(\ref{skoteinh}) with $\nabla^{(\nu)}\phi=\Delta\phi$ and with $\delta_R=j$.

\par We prove (\ref{elaela}) after showing that it implies Lemma
\ref{nabothewri2}.
\newline

{\it (\ref{elaela}) implies Lemma \ref{nabothewri2}:} Plugging 
(\ref{elaela}) into (\ref{boh9i}), we may assume that no complete contractions
$C^u_{g}\cdot \Delta\phi$ in (\ref{boh9i}) with $\delta_R=j-1$
have a factor $\nabla^{(m)}R_{ijkl}$ with $A$ internal
contractions. In other words 
(\ref{boh9i})\footnote{See also the notational convention introduced in (\ref{jeep})} 
can now be rewritten in the form:

\begin{equation}
\label{lambetaki} \begin{split} & \int_{M^n}\sum_{l\in L^{j,\vec{\kappa}_{max}}_\sigma} a_l \sum_{f=\sigma-b\cdot k_{min}+1}^\sigma 
(Const)_{l,f} C^{l,\iota|i_f}(g)\nabla_{i_f}+\sum_{negligible} C_g(\nabla\phi)+
\\&\Sum_{l\in L^{j,\sharp}_\sigma}
a_l  C^{l,\iota}(g)\nabla_{i_j}\phi+ \Sum_{h\in H} a_h
C^h_{g}(\nabla\phi)+ \Sum_{u\in U^{0,*}_{\Delta\phi,\delta= j-1}}
a_u C^u_{g}\cdot \Delta\phi+
\\&\Sum_{u\in U^0_{\Delta\phi,\delta\ge j}} a_u C^u_{g}\cdot
\Delta\phi +\Sum_{x\in X} a_x C^x_{g}(\phi) dV_{g}=0.
\end{split}
\end{equation}
Here each compete contraction indexed in $L^{j,\sharp}_\sigma$ has $j-1$
internal contractions and a factor $\nabla\phi$, but its double character
is {\it not} $\vec{\kappa}'_{max}$. The complete contractions
indexed in $H$ are in the form (\ref{skoteinh}) with
$\delta_R\ge j$ and no factors $\nabla^{(p)}Ric$.
 The complete contractions in
$U^0_{\Delta\phi,\delta\ge j}$ have a factor $\Delta\phi$ and
$\delta\ge j$, while the ones indexed in
$U^{0,*}_{\Delta\phi,\delta= j-1}$ have $\delta=j-1$ but also have
no factor $\nabla^{(m)}R_{ijkl}$ with $A$ internal contractions.

\par Integrating by parts with respect to the factor $\Delta\phi$ in the second line of 
(\ref{lambetaki}) we derive:

\begin{equation}
\label{lambetaki2} \begin{split} & \int_{M^n}\sum_{l\in L^{j,\vec{\kappa}_{max}}_\sigma} 
a_l \sum_{f=\sigma-b\cdot k_{min}+1}^\sigma 
(Const)_{l,f} C^{l,\iota|i_f}(g)\nabla_{i_f}+\sum_{negligible} C_g(\nabla\phi)+
\\&\Sum_{l\in L^{j,\sharp}_\sigma}
a_l  C^{l,\iota}(g)\nabla_{i_a}\phi+\Sum_{h\in H} a_h
C^h_{g}(\nabla\phi)- \Sum_{u\in U^{0,*}_{\Delta\phi,\delta= j-1}}
a_u \nabla^i[C^u_{g}]\cdot\nabla_i\phi-
\\&\Sum_{u\in U^0_{\Delta\phi,\delta\ge j}} a_u \nabla^i[C^u_{g}]\cdot
\nabla_i\phi +\Sum_{x\in X} a_x C^x_{g}(\phi) dV_{g}=0.
\end{split}
\end{equation}
(Note that the above equation has been derived independently of (\ref{elaela})).

 We denote the integrand of the above by $M_{g}(\phi)$.
We apply the ``main conclusion'' of the super divergence formula to 
the above integral equation (see subsection 2.2 in \cite{alexakis1})
 and we pick out the sublinear combination $supdiv_{+}[M_{g}
(\phi)]$ of terms with length $\sigma+1$, a factor $\nabla\phi$
 and no internal contractions. Since the super divergence formula holds formally we derive that 
$supdiv_{+}[M_{g}(\phi)]=0$ modulo a linear combination of terms of length $\ge\sigma+2$. 
Thus we derive an equation: 

\begin{equation}
\label{lambetaki3} \begin{split} & \Sum_{l\in
L^{j,\vec{\kappa}_{max}}_\sigma} a_l \Sum_{f=0}^{b\cdot k_{min}-1}
(Const)_{l,f}Xdiv_{i_1}\dots \hat{Xdiv}_{i_{j-f}}\dots
Xdiv_{i_a}C^{l,\iota|i_1\dots i_j}(g)\nabla_{i_{j-f}}\phi+
\\&\sum_{negligible} Xdiv_{i_1}\dots Xdiv_{i_{j-1}}C^{i_1\dots i_{j-1}}_g(\nabla\phi)+
\Sum_{l\in L^{j,\sharp}_\sigma}
a_l Xdiv_{i_1}\dots Xdiv_{i_{j-1}}\\& C^{l,\iota|i_1\dots
i_j}(g)\nabla_{i_j}\phi+\Sum_{h\in H} a_h Xdiv_{i_1}\dots
Xdiv_{i_{z_h}}C^{h,i_1\dots i_{z_h}}_{g}(\nabla\phi)
\\& -\Sum_{u\in U^{0,*}_{\Delta\phi,\delta= j-1}} a_u
Xdiv_{i_1}\dots Xdiv_{i_{j-1}}\nabla^i[C^{u,i_1\dots
i_{j-1}}_{g}]\cdot \nabla_i\phi
\\&-\Sum_{u\in U^0_{\Delta\phi,\delta\ge j}} a_u Xdiv_{i_1}\dots
Xdiv_{i_{u_h}}\nabla^i [C^u_{g}]^{i_1\dots i_{u_h}}\cdot
\nabla_i\phi =0;
\end{split}
\end{equation}
(the linear combination $\sum_{negligible} Xdiv_{i_1}\dots Xdiv_{i_{j-1}}C^{i_1\dots i_{j-1}}_g(\nabla\phi)$ 
just arises from the linear combination of complete contractions $\sum_{negligible} C_g(\nabla\phi)$ 
by making all the internal contractions into free 
indices and then taking $Xdiv$ of the resulting free indices).

We now
apply Corollary \ref{precorgiade} to the above; we pick out the sublinear
combination of $(j-1)$-tensor fields with $A$ factors
$\nabla\upsilon$ contracting against some factor
$\nabla^{(m)}R_{ijkl}$ (this sublinear combination vanishes
separately). We then pick out the sublinear combination of complete contractions with
the  property that the factors $\nabla\upsilon$ are contracting according to the following
pattern:  $k_1$ factors $\nabla\upsilon$ must contract against one factor $T_1$; 
 $k_2$ factors $\nabla\upsilon$ must contract against some other factor $T_2$;
 $k_3$ factors $\nabla\upsilon$ must contract against a third factor $T_3$;
$\dots$; $k_a-1$ factors $\nabla\upsilon$ must contract
against an $a^{th}$ factor $T_a$ {\it and, in this last case}
the factor $\nabla\phi$ must also contract against this factor $T_a$.\footnote{Recall that the
double-character $\vec{\kappa}_{max}'$ was defined to be the
character $(k_1,\dots,k_2,k_3,\dots, k_{a-1},k_a-1)$.} Observe that the sublinear
combination of those terms {\it must vanish separately}, since the equation to
which this sublinear combination belongs holds formally.
We thus obtain a new
true equation which is in the form:

\begin{equation}
\label{georgia}
\begin{split}
&\Sum_{l\in L^{j,\vec{\kappa}_{max}}_\sigma} a_l
\Sum_{f=0}^{b\cdot k_{min}-1} C^{l,\iota|i_1\dots
i_j}(g)\nabla_{i_{j-f}}\phi\nabla_{i_1}\upsilon\dots\hat{\nabla}_{i_{j-f}}\upsilon
\dots \nabla_{i_j}\upsilon+
\\&\sum_{negligible} C^{i_1\dots i_{j-1}}_g(\nabla\phi)\nabla_{i_1}\upsilon\dots \nabla_{i_{j-1}}\upsilon
 \\&-Xdiv_{i_{j+1}}\Sum_{h\in H} a_h C^{h,i_1\dots
i_{j+1}}_{g}\nabla_{i_1}\phi\nabla_{i_2}\upsilon\dots\nabla_{i_j}\upsilon=0;
\end{split}
\end{equation}
(all tensor fields indexed in $H$ are acceptable).

\par Then, setting $\upsilon=\phi$\footnote{Notice that this operation 
will kill the sublinear combination $\sum_{negligible}\dots$.} 
and then applying the operation
$Weylify$ (see subsection 5.1 in  \cite{alexakis1} and also (\ref{themainthing})), 
we derive that Lemma \ref{nabothewri2}
follows from (\ref{elaela}). $\Box$
\newline

{\it Proof of (\ref{elaela}):} We again refer to
(\ref{lambetaki2}),\footnote{Recall the note after that equation,
which shows that it has been derived independently of (\ref{elaela}).} 
(denote the integrand by $M_g(\phi)$) and apply the super
divergence formula to this equation; denote the 
resulting local equation by $supdiv[M_g(\phi)]=0$.
 We pick out the sublinear combination $supdiv[M_g(\phi)]_{+}$ of terms of
length $\sigma+1$ with a factor $\nabla\phi$ and with no internal contractions.
This sublinear combination must vanish separately and we thus derive an equation:

\begin{equation}
\label{lambetaki3'} \begin{split} & \Sum_{l\in
L^{j-1}} a_l
Xdiv_{i_1}\dots
Xdiv_{i_{j-1}}C^{l|i_1\dots i_{j}}(g)\nabla_{i_{j}}\phi
\\&+\Sum_{h\in H} a_h Xdiv_{i_1}\dots
Xdiv_{i_{z_h}}C^{h,i_1\dots i_{z_h}}_{g}(\nabla\phi)
\\& -\Sum_{u\in U^{0,*}_{\Delta\phi,\delta= j-1}} a_u
Xdiv_{i_1}\dots Xdiv_{i_{j-1}}\nabla^i[C^{u,i_1\dots
i_{j-1}}_{g}]\cdot \nabla_i\phi
\\&-\Sum_{u\in U^0_{\Delta\phi,\delta\ge j}} a_u Xdiv_{i_1}\dots
Xdiv_{i_{u_h}}\nabla^i [C^u_{g}]^{i_1\dots i_{u_h}}\cdot
\nabla_i\phi =0;
\end{split}
\end{equation}
here the terms indexed in $L^{j-1}$ are in the form (\ref{preform2}), have rank $j-1$
{\it and also have at most $A-1$ free indices belonging to
the factor against which $\nabla\phi$ contracts}; 
they have arisen from the complete contractions
in $Image^1_{\nabla\phi,*,+}[P(g)|_\sigma]$ 
in (\ref{boh9i}). The tensor fields 
$C^u_{g}]^{i_1\dots i_{u_h}}$, $C^{u,i_1\dots i_{j-1}}_{g}$
arise from the complete contractions $C^h_g(\phi)$, $C^u(g)$
 in (\ref{boh9i}) by replacing
each internal contraction by a free index.

We denote by $H_{*}\subset H$ the index set of tensor 
fields in the above equation where
the factor $\nabla\phi$ is contracting against an 
internal index. We also denote by $L^{j-1,*}\subset L^{j-1}$
the index set of $(j-1)$-tensor fields for which the factor $\nabla\phi$ is 
contracting agianst an internal index in some factor 
$\nabla^{(m)}R_{abcd}$.\footnote{Recall that for each of these last tensor fields
 there can be at most $A-1$ free indices in any given factor.} 
Now, by applying Lemma \ref{cruch} to the above equation, we derive
that we can write: 

\begin{equation}
\label{avaggelsamiou}
\begin{split} 
&\Sum_{l\in
L^{j-1,*}} a_l
Xdiv_{i_1}\dots
Xdiv_{i_{j-1}}C^{l|i_1\dots i_{j}}(g)\nabla_{i_{j}}\phi
\\&+\sum_{h\in H_{*}} a_h Xdiv_{i_1}\dots Xdiv_{i_{z_h}} C^{h,i_1\dots i_{z_h}}_g(\nabla\phi)=
\\&\Sum_{l\in
\tilde{L}^{j-1}} a_lXdiv_{i_1}\dots
Xdiv_{i_{j-1}}C^{l|i_1\dots i_{j}}(g)\nabla_{i_{j}}\phi
\\&+\sum_{h\in H^\sharp} a_h Xdiv_{i_1}\dots Xdiv_{i_{z_h}} C^{h,i_1\dots i_{z_h}}_g(\nabla\phi),
\end{split}
\end{equation}
where the terms indexed in $\tilde{L}^{j-1}$ the RHS stand for a generic 
linear combination of tensor fields in the form (\ref{preform2}),
each with rank $j-1$ {\it and with the factor $\nabla\phi$ contracting against a 
derivative index  and at most $A-1$ free indices in any given factor}; $H^\sharp$ 
indexes a generic linear combination of tensor fields in the form (\ref{preform2}),
each with rank $>j-1$ {\it and with the factor $\nabla\phi$ contracting against a 
derivative index}. Therefore, replacing the above into (\ref{lambetaki3'}) 
we derive a new equation in the form:

 \begin{equation}
\label{lambetaki3'b} \begin{split} & \Sum_{l\in
\tilde{L}^{j-1}} a_l
Xdiv_{i_1}\dots
Xdiv_{i_{j-1}}C^{l|i_1\dots i_{j}}(g)\nabla_{i_{j}}\phi
\\&+\Sum_{h\in H^\sharp} a_h Xdiv_{i_1}\dots
Xdiv_{i_{z_h}}C^{h,i_1\dots i_{z_h}}_{g}(\nabla\phi)-
\\& \Sum_{u\in U^{0,*}_{\Delta\phi,\delta= j-1}} a_u
Xdiv_{i_1}\dots Xdiv_{i_{j-1}}\nabla^i[C^{u,i_1\dots
i_{j-1}}_{g}]\cdot \nabla_i\phi-
\\&\Sum_{u\in U^0_{\Delta\phi,\delta\ge j}} a_u Xdiv_{i_1}\dots
Xdiv_{i_{u_h}}\nabla^i [C^u_{g}]^{i_1\dots i_{u_h}}\cdot
\nabla_i\phi =0,
\end{split}
\end{equation}
where the linear combinations indexed in $\tilde{L}^{j-1}$ and $H^\sharp$ are {\it generic} 
linear combinations of the forms described in the above paragraph.

 Now, for each
$u\in SU^0_{\Delta\phi,\delta=j-1}$, we denote by
$\nabla_{@}^i[C^{u,i_1\dots i_{j-1}}(g)]\nabla_i\phi$ the
 sublinear combination in $\nabla^i[C^{u,i_1\dots i_{j-1}}(g)]\nabla_i\phi$
where $\nabla^i$ is {\it forced} to hit a factor with $A$ free indices.\footnote{Recall that
$SU^0_{\Delta\phi,\delta=j-1}\subset U^0_{\Delta\phi,\delta=j-1}$ is
 exactly the index set of $(j-1)$-tensor fields which have at least one such factor.}
Now, we break up the index set
$SU^0_{\Delta\phi,\delta=j-1}$ into subsets
$SU^{0,\vec{\kappa}}_{\Delta\phi,\delta=j-1}$
 that index tensor fields with the same
 character.\footnote{In other words with the same pattern of distribution of 
free indices among the different factors in $C^{u,i_1\dots i_{j-1}}(g)$.} We denote by $K$ the index
 set of those subsets.
We claim that for every $\vec{\kappa}\in K$ there is a linear combination of $(j+1)$-vector fields $\Sum_{h\in
H^{\vec{\kappa}'}} a_h C^{i_1\dots i_{j+1}}(g)\nabla_{i_1}\phi$
with a  character $\vec{\kappa}$ (also recall that $\nabla\phi$ is contracting against a derivative index)
so that:

\begin{equation}
\label{panagiot}
\begin{split}
&\Sum_{u\in SU^{0,\vec{\kappa}}_{\Delta\phi,\delta=j-1}} a_u
\nabla_{@}^i[C^{u,i_2\dots i_j}(g)]\nabla_i\phi\nabla_{i_2}
\upsilon\dots \nabla_{i_j}\upsilon
\\&- Xdiv_{i_{j+1}}\Sum_{h\in
H^{\vec{\kappa}'}} a_h C^{i_1\dots
i_{j+1}}(g)\nabla_{i_1}\phi\nabla_{i_2} \upsilon\dots
\nabla_{i_j}\upsilon=0.
\end{split}
\end{equation}

Let us check how (\ref{panagiot}) implies (\ref{lambetaki2}):

 We observe that for each $\vec{\kappa}\in K$, there is a
 nonzero combinatorial constant $(Const)_{\vec{\kappa}}$ so
 that for any $u\in SU^{0,\vec{\kappa}}_{\Delta\phi,\delta=j-1}$:

\begin{equation}
\label{triantaoxtw}
\begin{split}
Erase_{\nabla\phi}\{ \nabla_{@}^i[C^{u,i_1\dots i_{j-1}}(g)]\nabla_i\phi\}=
(Const)_{\vec{\kappa}}C^{u,i_1\dots i_{j-1}}(g).
\end{split}
\end{equation}
(In fact $(Const)_{\vec{\kappa}}$ just stands for the number
 of factors in $\vec{\kappa}$ that have $A$ free indices).
To see this clearly observe that if $\vec{\kappa}=(s_1,\dots s_a)$ (where $s_i\ge s_{i+1}$
and $s_1=A$) that means that each tensor field $C^{u,i_1\dots i_{j-1}}(g)$
will have a factor $F_1$ with $s_1=A$ free indices, a second factor $F_2$ with
$s_2$ free indices, $\dots,$ and an $a^{th}$ factor $F_a$ with $s_a$ free indices.
Then if $A=s_1=s_2=\dots s_c$ and $s_{c+1}<A$,
$\nabla^i_{@}[C^{u,i_1\dots i_{j-1}}(g)]$ stands for the sublinear combintaion in
$\nabla^i[C^{u,i_1\dots i_{j-1}}(g)]$ where $\nabla^i$ is forced to
hit one of the factors $F_1,\dots, F_c$. 

Now, applying the operation $Erase_{\phi}$ (see the Appendix in \cite{alexakis1})
to the above and
 then making the factors $\nabla\upsilon$ into internal
 contractions  and finally multiplying by $\Delta\phi$, we
derive (\ref{elaela}).
\newline

{\it Proof of (\ref{panagiot}):}
We apply Proposition \ref{pregiade2} to (\ref{lambetaki3'b}), deriving that
 there is a linear combination of tensor fields of type A1, $\sum_{t\in T} a_t C^{t,(i_1\dots 
i_{j-1}) i_ji_{j+1}}(g)\nabla_{i_j}\phi$,  
(this means that the factors $\nabla\phi$ is contracting against a 
derivative index in some factor $\nabla^{(m)}R_{ijkl}$) so that: 

 \begin{equation}
\label{steepest} \begin{split} & 
\Sum_{l\in \tilde{L}^{j-1}} a_l
C^{l|i_1\dots i_{j-1}i_j}(g)\nabla_{i_{j}}\phi\nabla_{i_1}\upsilon\dots\nabla_{i_{j-1}}\upsilon
 \\&-\Sum_{u\in U^{0,*}_{\Delta\phi,\delta= j-1}} a_u
\nabla^i[C^{u,i_1\dots
i_{j-1}}_{g}]\cdot \nabla_i\phi\nabla_{i_1}\upsilon\dots\nabla_{i_{j-1}}\upsilon
\\&=\sum_{t\in T} a_t Xdiv_{i_{j+1}}C^{t,i_1\dots 
i_{j-1} i_ji_{j+1}}(g)\nabla_{i_j}\phi.\nabla_{i_1}\upsilon\dots\nabla_{i_{j-1}}\upsilon.
\end{split}
\end{equation}
Now, in the above equations we pick out the sublinear combination of terms 
 where the factor $\nabla\phi$ and the factors
$\nabla\upsilon$ are contracting according to the following pattern:
The factor $\nabla\phi$ must contract against a
factor $F_1$ which is also contracting against
another $A$ factors $\nabla\upsilon$. A
second factor $F_2$ contracts against $s_2$
factors $\nabla\upsilon$, $\dots$ and
an $a^{th}$ factor $F_a$ contracts against
$s_a$ factors $\nabla\upsilon$. {\it This sublinear combination must vanish 
separately since (\ref{steepest}) holds formally}. Thus 
we derive (\ref{panagiot}).

\par We have proven (\ref{elaela}).
 $\Box$
\newline

{\bf Proof of Lemma \ref{work}:}
\newline

\par We  prove (\ref{liapkin}) and (\ref{sabatto}).
We will start with  (\ref{liapkin}), and we will first prove this
equation under a simplifying assumption;  our simplifying assumption is that no 
complete contraction in $\Sum_{u\in U^\alpha_{\nabla\phi}\bigcup
U^\alpha_{\Delta\phi}} a_u C^u_{g}(\phi)$ has a factor $R$  of the
scalar curvature. After we complete this proof under the 
simplifying assumption, we will explain how the general case can be derived by 
fitting the argument below to a new downward induction on the maximum 
number of factors $R$ of the scalar curvature, precisely 
as in subsection 5.4 in \cite{alexakis1}.

{\it Proof of (\ref{liapkin}) under the simplifying assumption:}
 Denote by $C^u_{g}(\phi,\Omega^\alpha)$ the
 complete contraction that arises
 from each $C^u_{g}(\phi)$, $u\in U^\alpha$,
by formally replacing the $\alpha$ factors $\nabla^{(p)}Ric$ by factors
$\nabla^{(p+2)}\Omega$.  Thus, $C^u_g(\phi)$, $u\in U_{\nabla\phi}^\alpha,u\in U^\alpha_{\Delta\phi}$  
will be in one of the forms:

\begin{equation}
\label{worthwhilec}
\begin{split}
&contr(\nabla^{a_1\dots
a_t}\nabla^{(m_1)}R_{ijkl}\otimes\dots\nabla^{b_1\dots
b_v}\nabla^{(m_s)}R_{i'j'k'l'}
\\& \otimes\nabla^{y_1\dots y_w} \nabla^{(p_1+2)}\Omega
\otimes\dots\nabla^{x_1\dots
x_o}\nabla^{(p_\alpha+2)}\Omega\otimes\nabla\phi),
\end{split}
\end{equation}

\begin{equation}
\label{Laplac}
\begin{split}
&contr(\nabla^{a_1\dots
a_t}\nabla^{(m_1)}R_{ijkl}\otimes\dots\nabla^{b_1\dots
b_v}\nabla^{(m_s)}R_{i'j'k'l'}
\\& \otimes\nabla^{y_1\dots y_w}
\nabla^{(p_1+2)}\Omega \otimes\dots\nabla^{x_1\dots
x_o}\nabla^{(p_\alpha+2)}\Omega\otimes\Delta\phi).
\end{split}
\end{equation}

For brevity, we write out:

\begin{equation}
\label{meson} \Sum_{u\in U^{\alpha,\delta<j}_{\Delta\phi}} a_u
C^u_{g}(\phi,\Omega^\alpha)=\Sum_{q\in Q} a_q C^q_{g}(\Omega^\alpha)\Delta\phi,
\end{equation}
where each $C^q_{g}(\Omega^\alpha)\Delta\phi$ is in the form
(\ref{Laplac}) with length $\sigma$ and
$\delta_R+\delta_\Omega+\alpha<j$. We denote by $C^{q,i_1\dots
i_a}_{g}(\Omega^\alpha)$ the tensor field that arises from each
$C^q_{g}(\Omega^\alpha)$
 by making all the internal contractions into free indices.

\par Applying the ``main conclusion'' of the super divergence formula 
(see \cite{alexakis1}) to
$\int_M Y_{g}(\phi)dV_g=0$, and picking out the sublinear 
combination of length $\sigma+1$ witha factor $\nabla\phi$, we obtain a
local equation in the form:

\begin{equation}
\label{eotvos}
\begin{split} &\Sum_{q\in Q} a_q Xdiv_{i_1}\dots Xdiv_{i_a}
\nabla^i[C^{q,i_1\dots i_a}_{g}(\Omega^\alpha)]\nabla_i
\phi
\\&+\Sum_{h\in H} a_h Xdiv_{i_1}\dots Xdiv_{i_b} C^{h,i_1\dots
i_b}_{g}(\Omega^\alpha,\phi)=0,
\end{split}
\end{equation}
modulo complete contractions of length $\ge\sigma +2$.\footnote{Recall
that $Xdiv_i[\dots]$ in this context stands for the 
sublinear combination in $div_i[\dots]$ where the derivative $\nabla^i$ 
is not allowed to hit $\nabla\phi$, {\it nor} the factor to 
which the free index ${}_i$ belongs.} All the
tensor fields above are acceptable without internal
 contractions. Moreover, each $b\ge j-1$ (and $b>j-1$ if we are additionnaly 
have $U^\alpha_{\nabla\phi,j-1=\emptyset}$) 
while each $a\le j-1$.

\par We pick out the minimum rank $a$ appearing 
above (we have denoted it by $\delta_{min}$, where 
$\delta_{min} < j-1$, or $\delta_{min}=j-1$ under the additionnal 
assumption that $U^\alpha_{\nabla\phi,j-1}=\emptyset$.) and
 denote by $Q^{\delta_{min}}\subset Q$ the corresponding index set.
We will show that there is a linear combination of acceptable
 $(\delta_{min}+1)$-tensor fields, $\Sum_{h\in H} a_h
C^{h,i_1\dots i_{\delta_{min}+1}}_{g}(\Omega^\alpha)$ so that:

\begin{equation}
\label{ciragan}
\begin{split}& \Sum_{q\in Q^{\delta_{min}}} a_q C^{q,i_1\dots
i_a}_{g}(\Omega^\alpha)
\nabla_{i_1}\upsilon\dots\nabla_{i_{\delta_{min}}}\upsilon
\\&- \Sum_{h\in H}
a_h Xdiv_{i_{min+1}} C^{h,i_1\dots i_{\delta_{min}+1}}_{g}(\Omega^\alpha)
\nabla_{i_1}\upsilon\dots\nabla_{i_{\delta_{min}}}\upsilon=0,
\end{split}
\end{equation}
modulo complete contractions of length $\ge\sigma +1$. Let us
assume we have proven this; we will then show how to derive
(\ref{liapkin}).

{\it Proof that (\ref{ciragan}) implies (\ref{liapkin}):}
 We define an operation $Op$ that replaces each
factor $\nabla^{(p+2)}\Omega$ by a factor $\nabla^{(p)}Ric$ and
replaces each factor $\nabla\upsilon$ by an internal contraction
and finally multiplies the contraction we obtain by $\Delta\phi$
(notice that the operation $Op$ is almost exactly the same as the
operation $Riccify$, defined in the appendix of \cite{alexakis1}).
Applying $Op$ to (\ref{ciragan}) and using the fact that
(\ref{ciragan}) holds formally, we derive (\ref{liapkin}) (the proof of this 
fact is precisely the same as the operation ``Riccify'' in subsection 5.1 in \cite{alexakis1}).
\newline

{\it Proof of (\ref{ciragan}):}  Now, in order to derive
(\ref{ciragan}) from (\ref{eotvos}), we will distinguish two
subcases: Either $\alpha=\sigma$ or $\alpha<\sigma$.

\par We start with the second case. We pick out the sublinear
 combination in $\Sum_{h\in H} \dots$ in (\ref{eotvos})
where the factor $\nabla\phi$
 is contracting against a factor $\nabla^{(m)}R_{ijkl}$ (the condition $\alpha<\sigma$ 
guarantees that there is such a factor).
We denote by $H_{*}\subset H$ the index set of that sublinear
combination. We also denote by
$$\nabla^i_{*}[C^{q,i_1\dots i_a}_{g}(\Omega^\alpha)]\nabla_i\phi$$
the sublinear combination in
$$\nabla^i[C^{q,i_1\dots i_a}_{g}(\Omega^\alpha)]\nabla_i\phi$$
where $\nabla^i$ is only allowed to hit factors
$\nabla^{(m)}R_{ijkl}$.

\par Since (\ref{eotvos}) holds formally, we derive that:

\begin{equation}
\label{eotvos2} 
\begin{split}
 &\Sum_{q\in Q} a_q Xdiv_{i_1}\dots Xdiv_{i_a}
\nabla^i_{*}[C^{q,i_1\dots i_a}_{g}(\Omega^\alpha)]\nabla_i
\phi\\&+\Sum_{h\in H_{*}} a_h Xdiv_{i_1}\dots Xdiv_{i_b}
C^{h,i_1\dots i_b}_{g}(\Omega^\alpha,\phi)=0.
\end{split}
\end{equation}

\par We then denote by $H^\sharp_{*}\subset H_{*}$ the index
set of those tensor fields for which the factor $\nabla\phi$ is
contracting against an internal index of the factor
$\nabla^{(m)}R_{ijkl}$. By 
applying the second claim in Lemma \ref{cruch} we deduce as before
that we can write:

\begin{equation}
\label{tennis} 
\begin{split}
&\Sum_{h\in H^\sharp_{*}} a_h Xdiv_{i_1}\dots
Xdiv_{i_b} C^{h,i_1\dots i_b}_{g}(\Omega^\alpha,\phi)
\\&= \Sum_{h\in H_{**}} a_h Xdiv_{i_1}\dots Xdiv_{i_b} C^{h,i_1\dots
i_b}_{g}(\Omega^\alpha,\phi),
\end{split}
\end{equation}
where each tensor field on the right hand side has $b> j-1$ ($b>j$ 
under the additonnal auumption that $U^\alpha_{\nabla\phi,j-1}=\emptyset$) and
the factor $\nabla\phi$ is contracting against the derivative
index of some factor $\nabla^{(m)}R_{ijkl}$. Substituting the above in
(\ref{eotvos2}) and applying
Proposition \ref{pregiade2} along with the operation
$Erase_{\nabla\phi}$, we derive (\ref{ciragan}).
\newline

\par Now, the case $\alpha=\sigma$: We recall equation (\ref{eotvos}) where
we polarize the function $\Omega$
 into $\Omega_1,\dots ,\Omega_\sigma$, obtaining  a new true equation. We denote the
 corresponding tensor fields by
$C^{q,i_1\dots i_a}_{g}(\Omega_1,\dots ,\Omega_\sigma, \phi)$,
$C^{h,i_1\dots i_b}_{g}(\Omega_1,\dots , \Omega_\sigma,\phi)$. We
denote by $H_{*}\subset H$ the index set of those tensor fields
for
 which $\nabla\phi$ is contracting against the factor
$\nabla^{(a)}\Omega_1$. We also denote by $\nabla^i_{*}$ the
sublinear combination in $\nabla^i$ where $\nabla^i$ is only
allowed to hit the factor $\nabla^{(a)}\Omega_1$.

\par Thus, we derive that:

\begin{equation}
\label{eotvos3}
\begin{split}
&\Sum_{q\in Q} a_q Xdiv_{i_1}\dots Xdiv_{i_a}
\nabla^i_{*}[C^{q,i_1\dots i_a}_{g}(\Omega_1,\dots ,
\Omega_\sigma)]\nabla_i\phi+
\\& \Sum_{h\in H_{*}} a_h Xdiv_{i_1}\dots Xdiv_{i_b}
C^{h,i_1\dots i_b}_{g}(\Omega_1,\dots ,\Omega_\sigma,\phi)=0.
\end{split}
\end{equation}

\par We denote by $H^\sharp_{*}\subset H_{*}$ the index
set of those tensor fields for which the factor $\nabla\phi$ is
contracting against a factor $\nabla^{(2)}\Omega_1$ (with exactly two derivatives). 
Then, by applying  Lemma 4.6 
in \cite{alexakis4},\footnote{As noted in \cite{alexakis4}, that Lemma is a {\it consequence} 
of the ``fundamental Proposition''. There is no logical dependence of that Lemma
on any of the work in this paper or in \cite{alexakis1,alexakis3}.} we can derive:

\begin{equation}
\label{tennis} 
\begin{split}&\Sum_{h\in H^\sharp_{*}} a_h Xdiv_{i_1}\dots
Xdiv_{i_b} C^{h,i_1\dots i_b}_{g}(\Omega^\sigma,\phi)
\\&= \Sum_{h\in
H_{**}} a_h Xdiv_{i_1}\dots Xdiv_{i_b} C^{h,i_1\dots
i_b}_{g}(\Omega^\sigma,\phi),
\end{split}
\end{equation}
where each tensor field on the right hand side has $b> j-1$ and
the factor $\nabla\phi$ is contracting against a factor
$\nabla^{(A)}\Omega$ with $A\ge 3$.\footnote{We only need
 to check that the hypotheses of this Lemma are satisfied
in the case where the minimum rank $b_{min}$ appearing
among the tensor fields in $H^\sharp_{*}$ is $b_{min}=1$, and in addition
there are tensor fields indexed in $H^\sharp_{*}$
 with rank 1 and the one free index ${}_{i_1}$ belonging to an
expression $\nabla^{(2)}_{ii_1}\Omega_1\nabla^i\phi$.
 In that case we see that the extra assumption needed
 for Lemma 4.6 is fulfilled by
 virtue of a weight restriction--this follows since $j\ge 2$ and
  all internal contractions in
  $P(g)|_{L_\sigma^j}$ involve only derivative indices.}
 Substituting the above in (\ref{eotvos2}) and applying
 the operation $Erase_{\nabla\phi}$ from the Appendix in \cite{alexakis1}, we
notice that all resulting complete contractions are 
acceptable (in particular they have each factor
 $\nabla^{(a_s)}\Omega_s$ having $a_s\ge 2$); thus 
we derive our claim (\ref{liapkin}).
\newline

{\it Proof of (\ref{sabatto}) under the simplifying assumption:}

\par The proof of this claim will be slightly more involved in 
the case where $\delta_{min}=j-1$.
We first prove (\ref{sabatto}) in the case $\delta_{min}>j-1$:

{\bf The case $\delta_{min}>j-1$:} For each $u\in U^\alpha_{\Delta\phi}$, we denote
by $C^{u,i_1}_g(\phi)$ the complete contraction that arises from $C^u_g(\phi)$ by
replacing $\Delta\phi$ by $\nabla_{i_1}\phi$. We observe that we can write:

$$\sum_{u\in U^\alpha_{\Delta\phi}} a_u C^u_g(\phi)-div_{i_1} C^{u,i_1}_g(\phi)=
\sum_{u\in U'^\alpha_{\nabla\phi}}C^u_g(\phi)$$
 where the linear combination on the RHS stands for a
 generic linear combination of complete contractions with a
  factor $\nabla\phi$, $\alpha$ factors $\nabla^{(p)}Ric$ and $\delta>j-1$.

  In view of this equation, it suffices to show (\ref{sabatto})
 under the additional assumption that $U^\alpha_{\Delta\phi}=\emptyset$.

Now, to show (\ref{sabatto}) under this extra assumption we consider the equation
$\int_M Y_g(\phi)dV_g=0$ (in the assumption of Lemma \ref{work}).
 We apply the ``main conclusion'' of the super divergence formula 
(see \cite{alexakis1}) and we pick out the sublinear combination of terms with length $\sigma+1$ and
 a factor $\nabla\phi$, thus deriving a new local equation: 

\begin{equation}
\label{robbanks}
\begin{split}
&\sum_{u\in U^\alpha_{\nabla\phi,\delta_{min}}} a_u
Xdiv_{i_1}\dots Xdiv_{i_{\delta_{min}}} C^{u,i_1\dots
i_{\delta_{min}}}_g(\phi,\Omega^\alpha)+
\\& \sum_{u\in U^\alpha_{\nabla\phi,\delta>\delta_{min}}} a_u
Xdiv_{i_1}\dots Xdiv_{i_\delta} C^{u,i_1\dots
i_a}_g(\phi,\Omega^\alpha)=0.
\end{split}
\end{equation}
Then (\ref{sabatto}) follows by applying Corollary \ref{precorgiade}
 to (\ref{robbanks});\footnote{Since $\delta_{min}> j-1\ge 1$ we do not have to worry about 
the extra restrictions in that corollary when $\alpha=1$.}
we derive that there is a linear combination of $(\beta+1)$-vector
fields (indexed in $U^\alpha_{\nabla\phi,\delta=\delta_{min}+1}$ below) so that:

\begin{equation}
\label{robbanks2}
\begin{split}
&\sum_{u\in U^\alpha_{\nabla\phi,\delta_{min}}} a_u
 C^{u,i_1\dots i_{\delta_{min}}}_g(\phi,\Omega^\alpha)
\nabla_{i_1}\dots \nabla_{i_{\delta_{min}}}\upsilon+
\\& \sum_{u\in U^\alpha_{\nabla\phi,\delta_{min}+1}} a_u
 Xdiv_{i_{\delta_{min}+1}} C^{u,i_1\dots i_{\delta_{min}+1}}_g
(\phi,\Omega^\alpha)\nabla_{i_1}\upsilon\dots \nabla_{i_{\delta_{min}}}\upsilon=0.
\end{split}
\end{equation}
By applying the operation $Riccify$ to the above (see subsection 5.1 in \cite{alexakis1})  
we derive our claim, in the case $\delta_{min}>j-1$. 
 \newline

Now, the case $\delta_{min}=j-1$: We 
introduce some notational conventions.  
We break up the linear combination $\sum_{u\in U^\alpha_{\Delta\phi,j-1}}\dots$
 into two pieces (by dividing the index set $U^\alpha_{\Delta\phi,j-1}$ into two subsets): 
We will say that $u\in U^{\alpha|I}_{\Delta\phi,j-1}$ if the complete contraction 
$C^u_g(\phi)$ has the property that all $\alpha$ factors $\nabla^{(p)}Ric$ have 
at least one internal contraction;\footnote{Not counting the one in $Ric_{ab}={R^i}_{aib}$ itself.}
 we  define $U^{\alpha|II}_{\Delta\phi,j-1}=U^{\alpha}_{\Delta\phi,j-1}\setminus U^{\alpha|I}_{\Delta\phi,j-1}$.

By applying the second Bianchi identity we derive that we can write:

\begin{equation}
\label{doukissa} 
\begin{split}
& \sum_{u\in U^{\alpha|I}_{\Delta\phi,j-1}} a_u C^u_g(\phi)=
 \sum_{u\in \tilde{U}^{\alpha|I}_{\Delta\phi,j-1}} a_u C^u_g(\phi)+
\sum_{u\in U^{\alpha-1}_{\Delta\phi}} a_u C^u_g(\phi),
\end{split}
\end{equation}
where the complete contractions indexed in $\tilde{U}^{\alpha|I}_{\Delta\phi,j-1}$ have all the 
generic properties of the complete contractions indexed in $U^{\alpha|I}_{\Delta\phi,j-1}$ {\it and} 
have the feature that the two indices ${}_a,{}_b$ in every factor $\nabla^{(p)}Ric_{ab}$ contract 
against each other;\footnote{In other words that factor is of the form $\nabla^{(p)}R$, 
where $R$ stands for the scalar curvature} the linear combination 
$\sum_{u\in U^{\alpha-1}_{\Delta\phi}} a_u C^u_g(\phi)$ stands for a generic linear combination 
as allowed in the RHS of (\ref{sabatto}). 

\par In view of the above, we may assume with no loss of generality that 
all the $\alpha$ factors $\nabla^{(p)}Ric_{ab}$ in every $C^u_g(\phi), u\in U^{\alpha|I}_{\Delta\phi,j-1}$ 
have the indices ${}_a,{}_b$ contracting against each other.
Now, for each $u\in U^{\alpha|I}_{\Delta\phi,j-1}$, we denote
by $C^{u,i_1}_g(\phi)$ the vector field that arises from $C^u_g(\phi)$ by
replacing $\Delta\phi$ by $\nabla_{i_1}\phi$. We observe that we can write:

\begin{equation}
\label{momentum}
\sum_{u\in U^{\alpha|I}_{\Delta\phi,j-1}} a_u C^u_g(\phi)-div_{i_1} C^{u,i_1}_g(\phi)=
\sum_{u\in U'^\alpha_{\nabla\phi}} a_u C^u_g(\phi),
\end{equation}
 where the linear combination on the RHS stands for a
 generic linear combination of complete contractions with a
  factor $\nabla\phi$, $\alpha$ factors $\nabla^{(p)}Ric$ and $\delta=j-1$ {\it and} with the feature that 
all the $\alpha$ factors $\nabla^{(p)}Ric_{ab}$  
have the indices ${}_a,{}_b$ contracting against each other. Thus, by 
virtue of (\ref{doukissa}), (\ref{momentum})
 we may assume with no loss of generality 
that $U^{\alpha|I}_{\Delta\phi,j-1}=\emptyset$.

\par Now, we again consider the ``main conclusion'' of the super divergence formula 
applied to the integral equation $\int_{M^n}Y_g(\phi)dV_g=0$ (see the statement of Lemma \ref{work})
 and we pick out the sublinear combination of terms with length $\sigma+1$ and
 a factor $\nabla\phi$, thus deriving a new local equation: 

\begin{equation}
\label{robbanks'}
\begin{split}
&\sum_{u\in U^\alpha_{\nabla\phi,j-1}} a_u
Xdiv_{i_1}\dots Xdiv_{i_{j-1}} C^{u,i_1\dots
i_\beta}_g(\phi,\Omega^\alpha)
\\&+\sum_{u\in U'} a_u
Xdiv_{i_1}\dots Xdiv_{i_{j-1}} C^{u,i_1\dots i_\beta}_g
(\phi,\Omega^\alpha)
\\& +\sum_{u\in U^\alpha_{\nabla\phi,\delta>j-1}} a_u
Xdiv_{i_1}\dots Xdiv_{i_{\delta}} C^{u,i_1\dots
i_{\delta}}_g(\phi,\Omega^\alpha)=0,
\end{split}
\end{equation}
where the tensor fields indexed in $U'$ have rank $j-1$ but at 
least one of the $\alpha$ factors 
$\nabla^{(B)}\Omega$ does not contain a free index.\footnote{These terms 
arise from the sublinear combination $\sum_{u\in U^{\alpha|I}_{\Delta\phi,j-1}} a_u C^u_g(\phi).$}
 The terms in $U^\alpha_{\nabla\phi,\delta>j-1}$ have rank $>j-1$. 

Now, we apply Corollary \ref{precorgiade}
 to (\ref{robbanks'}) and pick out the sublinear combination 
of terms where all factors $\nabla^{(p)}\Omega$ contract against at 
least one factor $\nabla\upsilon$; this sublinear combination must vanish separately,
thus we derive  a new equation:

\begin{equation}
\label{robbanks''}
\begin{split}
&\sum_{u\in U^\alpha_{\nabla\phi,j-1}} a_u
 C^{u,i_1\dots
i_\beta}_g(\phi,\Omega^\alpha)\nabla_{i_1}\upsilon\dots\nabla_{i_{j-1}}\upsilon
\\&- \sum_{h\in H} a_h
Xdiv_{i_{j}} C^{u,i_1\dots i_j}_g(\phi,\Omega^\alpha)\nabla_{i_1}\upsilon\dots\nabla_{i_{j-1}}\upsilon=0,
\end{split}
\end{equation}
where the tensor fields indexed in $H$ have each factor $\nabla^{(p)}\Omega$ contracting 
against at least one factor $\nabla\upsilon$. Now, recall that
for each of the terms indexed in $U^\alpha_{\nabla\phi,j-1}$ in the above, the 
last index ${}_{r_p}$ in each factor $\nabla^{(p)}_{r_1\dots r_p}\Omega$ is contracting against 
a factor $\nabla\upsilon$. By just permuting indices  we may assume 
that the same is true  for each of the terms
indexed in $H$. Then, {\it since (\ref{robbanks''}) 
holds formally} we may assume that this last index in 
each of the factors $\nabla^{(p)}\Omega$ is {\it not} permuted in the formal permutations 
of indices by which (\ref{robbanks''}) is made {\it formally} true. 
Then, applying the operation $Riccify$ to the above we derive that there is a vector field 
$\sum_{h\in H'} a_h C^{h,i}_g(\phi)$ so that:

\begin{equation}
\begin{split}
\label{letssee} 
&\sum_{u\in U^\alpha_{\nabla\phi,j-1}} a_u C^{u}_g(\phi)- 
\sum_{h\in H'} a_h C^{u,i_1}_g(\phi)=\sum_{u\in U^\alpha_{\nabla\phi,\delta>j-1}} a_u C^u_g(\phi)+
\sum_{x\in X} a_x C^x_g(\phi),
\end{split}
\end{equation}
(in the notational conventions of (\ref{sabatto})). Thus we are reduced to 
showing our claim in the case  where $U^\alpha_{\nabla\phi,j-1}=U^{\alpha|I}_{\nabla\phi,j-1}=\emptyset$.
Our claim then follows by appealing to (\ref{liapkin}), where we now have the additionnal 
assumption that $U^\alpha_{\nabla\phi,j-1}=\emptyset$. $\Box$
\newline

 {\it Proof of
(\ref{liapkin}), (\ref{sabatto}) without the simplifying
assumption:}
\newline

\par Now, we explain how to prove (\ref{liapkin}), (\ref{sabatto})  in the case
where there are factors $R$ among the contractions indexed in
$U^\alpha_{\nabla\phi}\bigcup U^\alpha_{\Delta\phi}$ (i.e. {\it
without} the simplifying assumption). 
Let $M\epsilon\gamma$ be the maximum number of
factors $R$ among the complete contractions indexed in 
$U^\alpha_{\nabla\phi}\bigcup U^\alpha_{\Delta\phi}$. 

 We then reduce ourselves to
the case where there are no such factors  
by a downward induction on $M\epsilon\gamma$. 
 If  $M\epsilon\gamma\le\sigma -3$ then we can prove 
(\ref{liapkin}), (\ref{sabatto}) by just repeating 
the proof under the simplifying assumption, as in the proof of Lemmas
5.4, 5.5 when $M\le \sigma-3$ in \cite{alexakis1}:
 We prove the claims for the sublinear combinations in the LHS which 
have $M\epsilon\gamma$ factors $R$, allowing terms in the RHSs with 
$M\epsilon\gamma-1$ factors $R$. The argument runs unobstructed when 
$M\epsilon\gamma\ge 3$, since (as in \cite{alexakis1})
 whenever we wish to invoke Corollary \ref{precorgiade} 
or the first technical Proposition from \cite{alexakis4}, the
equations to which we apply them have $\sigma\ge 3$. 
\newline

{\it The case $M\epsilon\gamma\le 2$:} We reduce ourselves to the case 
$M\epsilon\gamma\ge 3$ by an {\it explicit} construction of divergences, 
followed by an application of the ``main conclusion''  of the super divergence formula:

 In the setting of equations (\ref{liapkin}), (\ref{sabatto}) we denote by
 $U_{\Delta\phi,*}, U_{\nabla\phi,*}$ the index sets of 
 complete contractions $C^u(g)\Delta\phi, C^u_g(\phi)$ with at least
$\sigma-2$ factors $R$ (hence for those complete contractions
$\alpha$ can be $\sigma,\sigma-1$ or $\sigma-2$). We claim
 that there exist a vector fields as required in (\ref{liapkin}), (\ref{sabatto}), indexed in
$H,H'$ below, so that:

\begin{equation}
\label{liapkinR}
\begin{split}
& \Sum_{u\in U_{\Delta\phi,*}} a_u C^u_{g}(\phi)- div_i\Sum_{h\in
H} a_h C^{h,i}_{g}(\Delta\phi)= (Const)_* C^*(g)\Delta\phi+
\\& \Sum_{u\in U_{\Delta\phi,-}} a_u
C^u_{g}(\phi) +\Sum_{x\in X} a_x C^x_{g}(\phi),
\end{split}
\end{equation}

\begin{equation}
\label{sabattoR} 
\begin{split}&\Sum_{u\in U_{\nabla\phi,*}} a_u C^u_{g}(\phi)-
div_i\Sum_{h\in H'} a_h C^{h,i}_{g}(\phi)= (Const)_\sharp C^\sharp_{g}(\phi)
\\&+ \Sum_{u\in U_{\nabla\phi,\delta>j-1,-}} a_u C^u_{g}(\phi)+\Sum_{x\in X} a_x
C^x_{g}(\phi);
\end{split}
\end{equation}
Here the complete contractions $C^*, C^\sharp$ are 
{\it explicit} complete contractions in  forms  and respectively
(their form depends on the value fo $\alpha$). 
Moreover the complete contractions indexed in
$U_{\Delta\phi,-}, U_{\nabla\phi,\delta>j-1,-}$ have at least
$j$ internal contractions and also strictly fewer
than $\sigma-2$ factors $R$, and a factor $\Delta\phi,\nabla\phi$ respectively.

\par Finally,  we claim that there exists a vector field 
indexed in $H''$ below  as required in (\ref{sabatto}) so that:

\begin{equation}
\label{lars}
\begin{split}
&-(Const)_*\nabla^i[C^*(g)]\nabla_i\phi+ (Const)_\sharp C^u_{g}(\phi)- div_i\Sum_{h\in H''} a_h
C^{h,i}_{g}(\phi)=\\& \Sum_{u\in U_{\nabla\phi,\delta>j-1,-}} a_u
C^u_{g}(\phi) +\Sum_{x\in X} a_x C^x_{g}(\phi),
\end{split}
\end{equation}
with the same conventions as above in the RHS.

Clearly, if we can show these three Lemmas we will have
 then reduced ourselves to showing Lemma \ref{work} under the
additional assumption that each contraction with a factor
$\nabla\phi$ or $\Delta\phi$ can have at most $\sigma -3$
 factors $R$; this case  has already been settled.

{\it Mini-Proof of (\ref{liapkinR}), (\ref{sabattoR}):}
The divergences needed for these equations are constructed {\it explicitly}.
The (simple) technique we use is explained in great detail in section 3 of \cite{alexakis3}.  
Consider the two factors $T_1,T_2$ which are {\it not} 
in the form $R$ (of the scalar curvature). Consider any particular contractions
 between these two factors. If one of the indices is a derivative index, we explicitly 
subtract the divergence corresponding to that index. The correction terms we 
get are either allowed in the RHSs of our equations, or have 
increased the number of internal contractions.\footnote{We also apply 
the second Bianchi 
identity whenever necessary to create particular 
contractions as described above, whenever possible.} One can see that we can perform this 
explicit construction repeatedly to obtain the terms $C^*,C^\sharp$ in the RHSs. 

 {\it Mini-Proof of (\ref{lars}):} We perform the same explcicit construction as 
described above for the term $-(Const)_*\nabla^i[C^*(g)]\nabla_i\phi$, 
to  derive (\ref{lars}) with an additional term $[-2(Const)_*+(Const)_\sharp]C^u_{g}(\phi)$
in the RHS. Then, substituting this into the integral
 equation $\int_{M^n} Y_g(\phi)dV_g=0$ (see the assumption of Lemma \ref{work})
and applying the main conclusion of the super divergence formula, we derive that 
$-2(Const)_*+(Const)_\sharp=0$. $\Box$

\end{document}